\documentclass[11pt]{amsart}

\renewcommand{\bar}{\overline}

\renewcommand{\phi}{\varphi}

\newcommand{\pz}{\partial_z}
\newcommand{\pzb}{\partial_{\bar z}}

\usepackage{amsmath,amsthm,amscd}
\usepackage{calc}

\title
[]{Canonical Metrics on the Moduli Space of Riemann Surfaces II}
\author{Kefeng Liu}
\address{Department of Mathematics\\
University of California at Los Angeles\\ Los Angeles, CA 90095-1555, USA\\
Center of Mathematical Sciences, Zhejiang University, Hangzhou,
China} \email{liu@math.ucla.edu, liu@cms.zju.edu.cn}
\author{Xiaofeng Sun}
\address{Department of Mathematics\\
Harvard University\\ Cambridge, MA 02138, USA}
\email{xsun@math.harvard.edu}
\author{Shing-Tung Yau}
\address{Department of Mathematics\\
Harvard University\\ Cambridge, MA 02138, USA}
\email{yau@math.harvard.edu}
\date{\today}

\thanks{The authors are supported by the NSF}

\newtheorem{theorem}{Theorem}[section]

\newtheorem{lemma}{Lemma}[section]
\newtheorem{cor}{Corollary}[section]
\newtheorem{prop}{Proposition}[section]
\newtheorem{claim}{Claim}
\newtheorem{definition}{Definition}[section]

\theoremstyle{remark}
\newtheorem{rem}{Remark}[section]

\begin{document}
\maketitle

\numberwithin{equation}{section}

\tableofcontents

\newcommand{\M}{{\mathcal M}}
\section{Introduction}\label{intro}

In \cite{lsy1} we started the project to understand the canonical
metrics on the Teichm\"uller and the moduli spaces of Riemann
surfaces, especially the K\"ahler-Einstein metric. Our goal is to
understand the geometry and topology of the moduli spaces from
understanding those classical metrics, as well as to find new
complete K\"ahler metrics with good curvature property. In
\cite{lsy1} we studied in detail two new complete K\"ahler
metrics, the Ricci and the perturbed Ricci metric. In particular
we proved that the Ricci metric has bounded holomorphic
bisectional curvature, and the perturbed Ricci metric, not only
has bounded holomorphic bisectional curvature, but also has
bounded negative holomorphic sectional curvature, and bounded
negative Ricci curvature. By using the perturbed Ricci metric as a
bridge we were able to prove the equivalence of several classical
complete metrics on the Teichm\"uller and the moduli spaces,
including the Teichm\"uller metric, the Kobayashi metric, the
Cheng-Yau-Mok K\"ahler-Einstein metric, the McMullen metric, as
well as the Ricci and the perturbed Ricci metric. This also solved
an old conjecture of Yau about the equivalence of the
K\"ahler-Einstein metric and the Teichm\"uller metric.

 In this paper we continue our study on these metrics and other classical metrics,
 in particular the K\"ahler-Einstein metric, and the perturbed Ricci metric.
One of the main results is the good understanding of the
K\"ahler-Einstein metric, from which we will
  derive some corollaries about the geometry of the
 moduli spaces. We will first prove the equivalence of the
 Bergman metric and the Carath\'eodory metric to the K\"ahler-Einstein
 metric. This completes our project on comparing all of the known
 complete metrics on the Teichm\"uller and moduli spaces.
 We then prove that the Ricci curvature of the perturbed Ricci metric has
 negative upper and lower bounds, and it also has bounded geometry. Recall that it also has
 bounded negative holomorphic sectional curvature. The perturbed Ricci metric
 is the first known complete
 K\"ahler metric on the Teichm\"uller and the moduli space with such good negative curvature
  property. We then focus on the K\"ahler-Einstein
 metric, study
 in detail its boundary behaviors and prove that not only
 it has bounded
 geometry, but also all of the covariant derivatives of its curvature are
 uniformly bounded. It is natural to expect interesting applications of the good
 properties of the perturbed Ricci and the K\"ahler-Einstein
 metric. For example, as an application of our detailed understanding of these metrics
  we prove that the logarithmic cotangent bundle of the moduli space is stable in the
 sense of Mumford.

This paper is organized as follows. In Section 2, by using the
Bers' embedding theorem, we prove that both the Carath\'eodory
metric and the Bergman metric on the Teichm\"uller space are
equivalent to the Kobayashi metric which is equivalent to the
Ricci metric, perturbed Ricci metric, K\"ahler-Einstein metric,
and the McMullen metric by the work in \cite{lsy1}. The
equivalence between the Bergman metric and the K\"ahler-Einstein
metric was first conjectured by Yau as Problem 44 in his $120$
open problems in geometry \cite{yau2},  \cite{yau4}.

In Section 3, we will prove that both the Ricci metric and the
perturbed Ricci metric have bounded curvature. Especially, with a
suitable choice of the perturbation constant, the holomorphic
sectional curvature and the Ricci curvature of the perturbed Ricci
metric are pinched by negative constants. As a simple corollary,
we immediately see that the dual of the logarithmic cotangent
bundle has no nontrivial holomorphic section. By using Bers'
embedding theorem and minimal surface, we also prove that the
Teichm\"uller space equipped with either of these two metrics has
bounded geometry: bounded curvature and lower bound of injectivity
radius. McMullen proved that the McMullen metric has bounded
geometry.

Having a complete K\"ahler-Einstein metric puts strong
restrictions on the geometric structure of the moduli space. In
Section 4 we will study the cohomology classes defined by the
K\"ahler forms and Ricci forms of the Ricci metric and the
K\"ahler-Einstein metric. As a direct corollary, we will see
easily that the moduli space is of logarithmic general type. One
of the most interesting applications of these study is a proof of
the stability of the logarithmic cotangent bundle of the moduli
space in the sense of Mumford.

Finally in Section 5, we will study the bounded geometry of the
K\"ahler-Einstein metric. We set up the Monge-Amper\'e equation
from a new metric obtained by deforming the Ricci metric along the
K\"ahler-Ricci flow. The work in \cite{lsy1} provides a $C^2$
estimate. We follow the work of Yau in \cite{yau3} and do the
$C^3$ and $C^4$ estimates. These estimates imply that the
curvature of the K\"ahler-Einstein metric is bounded. The same
method proves that all of the covariant derivatives of the
curvature are bounded. This may be used to understand the
complicated boundary of the Teichm\"uller space.

Now we give the precise statements of the main results in this
paper. We fix an integer $g\geq 2$ and denote by $\mathcal
T=\mathcal T_g$ the Teichm\"uller space, and $\M_g$ the moduli
space of closed Riemann surfaces of genus $g$. Our first result is
the following theorem which will be proved in Section 2:

\begin{theorem}
The Bergman metric and the Carath\'eodory metric both are
equivalent to the Kobayashi metric, therefore to the
K\"ahler-Einstein metric on the Teichm\"uller space.

\end{theorem}

Recall that we say two metrics are equivalent if they are
quasi-isometric to each other. We note that the equivalence
between the Bergman metric and the K\"ahler-Einstein metric was
conjectured by Yau in \cite{yau2}, \cite{yau4}. The proof of the
first part of Theorem 1.1 only needs the Bers' embedding and the
most basic properties of the Bergman and the Carath\'eodory
metrics.

Our second main result, proved in Section 3, is about the
curvature properties of the Ricci and the perturbed Ricci metric.
We have two theorems, the first one is about the Ricci metric:

\begin{theorem}
The holomorphic bisectional curvature, the holomorphic sectional
curvature and the Ricci curvature of the Ricci metric $\tau$ on
the moduli space $\M_g$ are bounded.
\end{theorem}

And the second theorem is about the perturbed Ricci metric:
\begin{theorem}
For any constant $C>0$, the bisectional curvature of the perturbed
Ricci metric $\widetilde\tau=\tau+Ch$ is bounded. Furthermore,
with a suitable choice of $C$, the holomorphic sectional curvature
and the Ricci curvature of $\widetilde\tau$ are bounded from above
and below by negative constants.
\end{theorem}

Both of the above theorems are proved by a detailed analysis of
the boundary behavior of the metrics and their curvature. Together
with the following result, which is proved by using minimal
surface theory and the Bers' embedding, they imply that Ricci
metric and the perturbed Ricci metric both have bounded geometry.

\begin{cor}
The injectivity radius of the Teichm\"uller space equipped with
the Ricci metric or the perturbed Ricci metric is bounded from
below.
\end{cor}

Let $\bar\M_g$ be the Deligne-Mumford compactification of the
moduli space of Riemann surfaces. Let $\bar E$ denote the
logarithmic extension of the cotangent bundle of the moduli space.
The next result is an interesting consequence of our detailed
understanding of the K\"ahler-Einstein metric on the moduli
spaces. It is proved in Section 4:

\begin{theorem}
The first Chern class of $\bar E$ is positive and $\bar E$ is
Mumford stable with respect to $c_1(\bar E)$.
\end{theorem}

This theorem also implies that ${\bar\M}_g$ is of logarithmic
general type for any $g\geq 2$. We prove this theorem by using the
good control of the K\"ahler-Einstein metric and the Ricci metric
we obtained near the boundary of the moduli space.

 Our final result, proved in Section 5, is the following theorem about the
K\"ahler-Einstein metric:

\begin{theorem}
The K\"ahler-Einstein metric on the Teichm\"uller space ${\mathcal
T}_g$ has bounded geometry. Furthermore the covariant derivatives
of its curvature are all uniformly bounded.
\end{theorem}

This theorem is proved by using the K\"ahler-Ricci flow and the
method of Yau in his proof of the Calabi conjecture to obtain
higher order estimates of the curvature.

We will continue our study on the geometry of the moduli space and
Teichm\"uller space. The topics will include the goodness on the
moduli space of these metrics and the $L^2$-cohomology of these
metrics and the Weil-Petersson metric on the Teichm\"uller space,
and the convergence of the K\"ahler-Ricci flow starting from the
Ricci and the perturbed Ricci metric.

{\bf Acknowledgement.} We would like to thank H. D. Cao and F. Luo
for interesting discussions and helps. The second author would
like to thank P. Li, Z. Lu, C. McMullen, R. Wentworth and B. Wong
for their help and encouragement. The second author would
especially like to thank Professor R. Schoen for many help and
guidance.

\section{The Carath\'eodory Metric and the Bergman Metric}\label{BC}

In this section we prove that the Carath\'eodory metric and the
Bergman metric on the Teichm\"uller space are equivalent to the
Kobayashi metric by using the Bers' embedding theorem. This
achieves one of our initial goals on the equivalence of all known
complete metrics on the Teichm\"uller space.

We first describe the idea. By the Bers' embedding theorem, we
know that for each point $p$ in the Teichm\"uller space, we can
find an embedding map of the Teichm\"uller space into $\mathbb
C^n$ such that $p$ is mapped to the origin and the image of the
Teichm\"uller space contains the ball of radius $2$ and is
contained inside the ball of radius $6$. The Kobayashi metric and
the Carath\'eodory metric of these balls coincide and can be
computed directly. Also, both of these metrics have restriction
property. Roughly speaking, the metrics on a submanifold are
larger than those on the ambient manifold. We use explicit form of
these metrics on the balls together with this property to estimate
the Kobayashi and the Carath\'eodory metric on the Teichm\"uller
space and compare them on a smaller ball. On the other hand, the
norm defined by the Bergman metric at each point can be estimated
by using the quotient of peak sections at this point. We use upper
and lower bounds of these peak sections to compare the Bergman
metric, the Kobayashi metric and the Euclidean metric on a small
ball in the image under the Bers' embedding of the Teichm\"uller
space.

At first, we briefly recall the definitions of the Carath\'eodory,
Bergman and Kobayashi metric on a complex manifold. Please see
\cite{ko2} for details.

Let $X$ be a complex manifold and of dimension $n$. let $\Delta_R$
be the disk in $\mathbb C$ with radius $R$. Let $\Delta=\Delta_1$
and let $\rho$ be the Poincar\'e metric on $\Delta$. Let $p\in X$
be a point and let $v\in T_p X$ be a holomorphic tangent vector.
Let $\text{Hol}(X,\Delta_R)$ and $\text{Hol}(\Delta_R,X)$ be the
spaces of holomorphic maps from $X$ to $\Delta_R$ and from
$\Delta_R$ to $X$ respectively. The Carath\'eodory norm of the
vector $v$ is defined to be
\[
\Vert v\Vert_C=\sup_{f\in\text{Hol}(X,\Delta)}\Vert f_\ast
v\Vert_{\Delta,\rho}
\]
and the Kobayashi norm of $v$ is defined to be
\[
\Vert v\Vert_K=\inf_{f\in\text{Hol}(\Delta_R,X),\ f(0)=p,\
f'(0)=v}\frac{2}{R}.
\]

Now we define the Bergman metric on $X$. Let $K_X$ be the
canonical bundle of $X$ and let $W$ be the space of $L^2$
holomorphic sections of $K_X$ in the sense that if $\sigma\in W$,
then
\[
\Vert\sigma\Vert_{L^2}^2=\int_X
(\sqrt{-1})^{n^2}\sigma\wedge\bar\sigma<\infty.
\]
The inner product on $W$ is defined to be
\[
(\sigma,\rho)=\int_X (\sqrt{-1})^{n^2}\sigma\wedge\bar\rho
\]
for all $\sigma,\rho\in W$. Let $\sigma_1,\sigma_2,\cdots$ be an
orthonormal basis of $W$. The Bergman kernel form is the
non-negative $(n,n)$-form
\[
B_X=\sum_{j=1}^\infty(\sqrt{-1})^{n^2}\sigma_j\wedge\bar\sigma_j.
\]

With a choice of local coordinates $z_i,\cdots,z_n$, we have
\[
B_X=BE_X(z,\bar z)(\sqrt{-1})^{n^2}dz_1\wedge\cdots\wedge dz_n
\wedge d\bar z_1\wedge\cdots\wedge d\bar z_n
\]
where $BE_X(z,\bar z)$ is called the Bergman kernel function. If
the Bergman kernel $B_X$ is positive, one can define the Bergman
metric
\[
B_{i\bar j}=\frac{\partial^2\log BE_X(z,\bar z)}{\partial z_i
\partial \bar z_j}.
\]
The Bergman metric is well-defined and is nondegenerate if the
elements in $W$ separate points and the first jet of $X$.

We will use the following notations:
\begin{definition}
Let $X$ be a complex space. For each point $p\in X$ and each
holomorphic tangent vector $v\in T_p X$ , we denote by $\Vert
v\Vert_{B,X,p}$, $\Vert v\Vert_{C,X,p}$ and $\Vert v\Vert_{K,X,p}$
the norms of $v$ measured in the Bergman metric, the
Carath\'eodory metric and the Kobayashi metric of the space $X$
respectively.
\end{definition}

Now we fix an integer $g\geq 2$ and denote by $\mathcal T=\mathcal
T_g$ the Teichm\"uller space of closed Riemann surface of genus
$g$. Our main theorem of this section is the following:
\begin{theorem}\label{bck}
Let $\mathcal T$ be the Teichm\"uller space of closed Riemann
surfaces of genus $g$ with $g\geq 2$. Then there is a positive
constant $C$ only depending on $g$ such that for each point
$p\in \mathcal T$ and each vector $v\in T_p\mathcal T$, we have
\[
C^{-1}\Vert v\Vert_{K,\mathcal T,p}\leq
\Vert v\Vert_{B,\mathcal T,p}\leq
C\Vert v\Vert_{K,\mathcal T,p}
\]
and
\[
C^{-1}\Vert v\Vert_{K,\mathcal T,p}\leq
\Vert v\Vert_{C,\mathcal T,p}\leq
C\Vert v\Vert_{K,\mathcal T,p}.
\]
\end{theorem}

{\bf Proof.} We will show that the norms defined by these metrics
are uniformly equivalent at each point of $\mathcal T$. We first
collect some known results in the following lemma.

\begin{lemma}
Let $X$ be a complex space. Then
\begin{enumerate}
\item $\Vert \cdot\Vert_{C,X}\leq \Vert \cdot \Vert_{K,X}$;
\item Let $Y$ be another complex space and $f:X\to Y$ be a
holomorphic map. Let $p\in X$ and $v\in T_p X$. Then
$\Vert f_\ast(v)\Vert_{C,Y,f(p)}\leq \Vert v \Vert_{C,X,p}$
and $\Vert f_\ast(v)\Vert_{K,Y,f(p)}\leq \Vert v \Vert_{K,X,p}$;
\item If $X\subset Y$ is a connected open
subset and $z\in X$ is a point. Then
with any local coordinates we have $BE_Y(z)\leq BE_X(z)$;
\item If the Bergman kernel is positive, then at each point
$z\in X$, a peak section $\sigma$ at $z$ exists. Such a peak section
is unique up to a constant factor $c$ with norm $1$. Furthermore,
with any choice of local coordinates, we have
$BE_X(z)=|\sigma(z)|^2$;
\item If the Bergman kernel of $X$ is positive, then
$\Vert \cdot \Vert_{C,X}\leq 2\Vert \cdot \Vert_{B,X}$;
\item If $X$ is a bounded convex domain in $\mathbb C^n$, then
$\Vert \cdot\Vert_{C,X}= \Vert \cdot \Vert_{K,X}$;
\item Let $|\cdot |$ be the Euclidean norm and
let $B_r$ be the open ball with center $0$ and radius $r$ in
$\mathbb C^n$. Then for any holomorphic tangent vector $v$ at
$0$,
\[
\Vert v\Vert_{C,B_r,0}=\Vert v\Vert_{K,B_r,0}=\frac{2}{r}|v|
\]
where $|v|$ is the Euclidean norm of $v$.
\end{enumerate}
\end{lemma}

{\bf Proof.} The first six claims are Proposition 4.2.4,
Proposition 4.2.3, Proposition 3.5.18, Proposition 4.10.4,
Proposition 4.10.3, Theorem 4.10.18 and Theorem 4.8.13 of
\cite{ko2}.

The last claim follows from the second claim easily. By rotation,
we can assume that $v=b\frac{\partial}{\partial z_1}$. Let
$\Delta_r$ be the disk with radius $r$ in $\mathbb C$ with
standard coordinate $z$ and let $\widetilde
v=b\frac{\partial}{\partial z}$ be the corresponding tangent
vector of $\Delta_r$ at $0$. Now, consider the maps $i:\Delta_r\to
B_r$ and $j:B_r\to \Delta_r$ given by $i(z)=(z,0,\cdots,0)$ and
$j(z_1,\cdots,z_n)=z_1$. We have $i_\ast(\widetilde v)=v$ and
$j_\ast(v)=\widetilde v$. By the Schwarz lemma it is easy to see
that $\Vert\widetilde v\Vert_{C,\Delta_r,0}=\frac{2}{r}|\widetilde
v|$. So we have
\[
\Vert v\Vert_{C,B_r,0}\geq \Vert j_\ast(v)\Vert_{C,\Delta_r,0}
=\Vert\widetilde v\Vert_{C,\Delta_r,0}=\frac{2}{r}|\widetilde v|
=\frac{2}{r}|v|
\]
and
\[
\Vert v\Vert_{C,B_r,0}= \Vert i_\ast(v)\Vert_{C,B_r,0}
\leq\Vert\widetilde v\Vert_{C,\Delta_r,0}=\frac{2}{r}|v|.
\]
This shows that the last claim holds for the Carath\'eodory
metric. By the sixth claim, we know that the last claim also holds
for the Kobayashi metric. This finishes the proof.

\qed

Now we prove the theorem. We first compare the Carath\'eodory
metric and the Kobayashi metric. By the above lemma it is easy to
see that if $X\subset Y$ is a subspace, then $\Vert
\cdot\Vert_{C,Y}\leq \Vert \cdot \Vert_{C,X}$ and $\Vert
\cdot\Vert_{K,Y}\leq \Vert \cdot \Vert_{K,X}$. Let $p\in \mathcal
T$ be an arbitrary point and let $n=3g-3=\dim_{\mathbb C}\mathcal
T$. Let $f_p:\mathcal T\to \mathbb C^n$ be the Bers' embedding map
with $f_p(p)=0$. In the following, we will identify $\mathcal T$
with $f_p(\mathcal T)$ and $T_p\mathcal T$ with $T_0\mathbb C^n$.
We know that
\begin{eqnarray}\label{bc5}
B_2\subset \mathcal T\subset B_6.
\end{eqnarray}
Let $v\in T_0\mathbb C^n$ be a holomorphic tangent vector. By
using the above lemma we have
\begin{eqnarray}\label{bc10}
\Vert v\Vert_{C,\mathcal T,0}\leq \Vert v\Vert_{K,\mathcal T,0}
\end{eqnarray}
and
\begin{eqnarray}\label{bc20}
\Vert v\Vert_{C,\mathcal T,0}\geq
\Vert v\Vert_{C,B_6,0}=\frac{1}{3}|v|=\frac{1}{3}
\Vert v\Vert_{C,B_2,0}=\frac{1}{3}\Vert v\Vert_{K,B_2,0}
\geq \frac{1}{3}\Vert v\Vert_{K,\mathcal T,0}.
\end{eqnarray}
By combining the above two inequalities, we have
\[
 \frac{1}{3}\Vert v\Vert_{K,\mathcal T,0}\leq
\Vert v\Vert_{C,\mathcal T,0}\leq
\Vert v\Vert_{K,\mathcal T,0}.
\]
Since the above constants are independent of the choice of
$p$, we proved the second claim of the theorem.

Now we compare the Bergman metric and the Kobayashi metric. By the
above lemma we know that the Bergman norm is bounded from below by
half of the Carath\'eodory norm provided the Bergman kernel is
non-zero. For each point $p\in \mathcal T_g$, let $f_p$ be the
Bers' embedding map with $f_p(p)=0$. Since $f_p(\mathcal
T_g)\subset B_6$, by the above lemma we know that
$BE_{f_p(\mathcal T_g)}(0)\geq BE_{B_6}(0)$. However, we know that
the Bergman kernel on $B_6$ is positive. This implies that the
Bergman kernel is non-zero at every point of the Teichm\"uller
space.

By the above lemma and the equivalence of the Carath\'eodory metric
and the Kobayashi metric, we know that the Bergman metric is
bounded from below by a constant multiple of the Kobayashi
metric.

When we fix a point $p$ and the Bers' embedding map $f_p$, from
inequality \eqref{bc20} we know that
\begin{eqnarray}\label{bc30}
|v|\leq 3\Vert v\Vert_{C,\mathcal T,0}
\leq 3\Vert v\Vert_{K,\mathcal T,0}.
\end{eqnarray}
Let $z_1,\cdots,z_n$ be the standard coordinates on
$\mathbb C^n$ with $r_i=|z_i|$ and let
$dV=\left (\sqrt{-1}\right )^n dz_1\wedge d\bar z_1
\wedge\cdots\wedge dz_n\wedge d\bar z_n$ be the volume form. Let
$\sigma=\alpha(z)dz_1\wedge\cdots\wedge dz_n$ be a peak section
over $\mathcal T$ at $0$ such that
\begin{eqnarray*}
\int_{\mathcal T}|\alpha|^2\ dV=1.
\end{eqnarray*}
Then we have $BE_{\mathcal T}(0)=|\alpha(0)|^2$. Now we
consider a peak section $\sigma_1=\alpha_1(z)
dz_1\wedge\cdots\wedge dz_n$ over $B_6$ at $0$ with
$\int_{B_6}|\alpha_1|^2\ dV=1$. Similarly we have that
$BE_{B_6}(0)=|\alpha_1(0)|^2$. By the above lemma and \eqref{bc5}
we have
\begin{eqnarray}\label{bc40}
|\alpha(0)|=\left ( BE_{\mathcal T}(0)\right )^{\frac{1}{2}}
\geq \left ( BE_{B_6}(0)\right )^{\frac{1}{2}}=|\alpha_1(0)|.
\end{eqnarray}
Let $v_n=\int_{B_1}dV$ be the volume of the unit ball in
$\mathbb C^n$ and let
\[
w_n=\frac{1}{n}\int_{x_1^2+\cdots+x_n^2\leq 4,\ x_i\geq 0}
(x_1^2+\cdots+x_n^2)x_1\cdots x_n\ dx_1\cdots dx_n
\]
where $x_1,\cdots,x_n$ are real variables. We see that both $v_n$
and $w_n$ are positive constants only depending on $n=3g-3$.

Now we consider the constant section $\sigma_2=a\
dz_1\wedge\cdots\wedge dz_n$ over $B_6$ where
$a=6^{-\frac{n}{2}}v_n^{-\frac{1}{2}}$. we have $\int_{B_6}a^2\
dV=1$. Since $\sigma_1$ is a peak section at $0$, we know that
$|\alpha_1(0)|\geq a$. By using inequality \eqref{bc40} we have
\begin{eqnarray}\label{bc50}
|\alpha(0)|\geq 6^{-\frac{n}{2}}v_n^{-\frac{1}{2}}.
\end{eqnarray}

To estimate the Bergman norm of $v$, by rotation, we may assume
$v=b\frac{\partial}{\partial z_1}$. So $|v|=|b|$. Let
$\tau=f(z)dz_1\wedge\cdots\wedge dz_n$ be an arbitrary section
over $\mathcal T$ with $f(0)=0$ and
$\int_{\mathcal T}|f|^2\ dV=1$. We have
$\int_{B_2}|f|^2\ dV\leq 1$.

Let $I$ be the index set $I=\{ (i_1,\cdots,i_n)\mid i_k\geq 0,\
\sum i_k\geq 1\}$. Since $f(0)=0$ and $f$ is holomorphic, we
can expand $f$ as a power series on $B_2$ as
\[
f(z)=\sum_{(i_1,\cdots,i_n)\in I} a_{i_1\cdots i_n}
z_1^{i_1}\cdots z_n^{i_n}.
\]
This implies $df(v)=a_{10\cdots 0}b$. Since $\int_{B_2}|f|^2\
dV\leq \int_{\mathcal T}|f|^2\ dV=1$, we have
\begin{align*}
\begin{split}
1\geq & \int_{B_2}|f|^2\ dV=\int_{B_2}
\sum_{(i_1,\cdots,i_n)\in I} |a_{i_1\cdots i_n}|^2
r_1^{2i_1}\cdots r_n^{2i_n}\ dV\geq \int_{B_2}
|a_{10\cdots 0}|^2 r_1^2\ dV\\
=& |a_{10\cdots 0}|^2 (4\pi)^n
\int_{r_1^2+\cdots+r_n^2\leq 4}
r_1^3r_2\cdots r_n\ dr_1\cdots dr_n=|a_{10\cdots 0}|^2 (4\pi)^n w_n
\end{split}
\end{align*}
which implies
\begin{eqnarray}\label{bc60}
|a_{10\cdots 0}|\leq (4\pi)^{-\frac{n}{2}}w_n^{-\frac{1}{2}}.
\end{eqnarray}
So we have
\begin{eqnarray}\label{bc70}
|df(v)|=|a_{10\cdots 0}||b|
\leq (4\pi)^{-\frac{n}{2}}w_n^{-\frac{1}{2}}|v|.
\end{eqnarray}
Let $W'$ be the set of sections over $\mathcal T$ such that
\[
W'=\{\tau=f(z)dz_1\wedge\cdots\wedge dz_n\mid f(0)=0,\
\int_{\mathcal T}|f|^2\ dV=1\}.
\]
By combining \eqref{bc30}, \eqref{bc50} and \eqref{bc70} we have
\begin{align}\label{bc80}
\begin{split}
\Vert v\Vert_{B,\mathcal T,0}=&\sup_{\tau\in W'}
\frac{|df(v)|}{|\alpha(0)|}\leq
\frac{(4\pi)^{-\frac{n}{2}}w_n^{-\frac{1}{2}}|v|}
{6^{-\frac{n}{2}}v_n^{-\frac{1}{2}}}
=\left (\frac{3}{2\pi}\right )^{\frac{n}{2}}
\left (\frac{v_n}{w_n}\right )^{\frac{1}{2}}|v|\\
\leq & 3\left (\frac{3}{2\pi}\right )^{\frac{n}{2}}
\left (\frac{v_n}{w_n}\right )^{\frac{1}{2}}
\Vert v\Vert_{K,\mathcal T,0}.
\end{split}
\end{align}
Since the constant in the above inequality only depends on the
dimension $n$, we know that the Bergman metric is uniformly
equivalent to the Kobayashi metric. This finished the proof.

\qed

\begin{rem}
After we proved this theorem, the second author was informed by C.
McMullen that the equivalence of the Carath\'eodory metric and the
Kobayashi metric maybe known. A more interesting question is
whether these two metrics coincide or not. We would like to study
this problem in the future.
\end{rem}
\section{The Negativity of the Ricci Curvature of the Perturbed Ricci Metric}\label{ricciricci}

In this section, we first study the curvature bounds of the Ricci
metric and the perturbed Ricci metric. By using the Bers'
embedding theorem, we show that the injectivity radius of the
Teichm\"uller space equipped with the Ricci metric or the
perturbed Ricci metric is bounded from below. This implies that
both the Ricci metric and the perturbed Ricci metric have bounded
geometry on the Teichm\"uller space.

The boundedness of the curvatures of these metrics was obtained by
analyzing their asymptotic behavior. The proof of the negativity
of the holomorphic sectional curvature and Ricci curvature of the
perturbed Ricci metric is more delicate. Near the boundary of the
moduli space and in the degeneration directions, these curvatures
are dominated by the contribution from the Ricci metric which is
negative. In the nondegeneration directions and in the interior of
the moduli space, these curvatures are dominated by the
contribution of the constant multiple of the Weil-Petersson metric
when the constant is large which is also negative. However we know
that the curvature of the linear combinations of two metrics is
not linear, we need to handle the error terms carefully.

Let $\M_g$ be the moduli space of closed oriented Riemann surfaces
of genus $g$ with $g\geq 2$ and let $\mathcal T_g$ be the
corresponding Teichm\"uller space. Let $\bar\M_g$ be the
Deligne-Mumford compactification of $\M_g$ and let
$D=\bar\M_g\setminus\M_g$ be the compactification divisor. It is
well known that $D$ is a divisor of normal crossings. In
\cite{lsy1} we studied various metrics on $\M_g$ and $\mathcal
T_g$. We briefly recall the results here.

Fix a point $p\in \M_g$. Let $X=X_p$ be a Riemann surface
corresponding to $p$. Let $z$ be the local holomorphic coordinate
on $X$ and let $s_1,\cdots,s_n$ be local holomorphic coordinates
on $\M_g$ where $n=3g-3$ is the complex dimension of $\M_g$. Let
$HB(X)$ and $Q(X)$ be the spaces of harmonic Beltrami
differentials and holomorphic quadratic differentials on $X$
respectively and let $\lambda$ be the hyperbolic metric on $X$.
Namely,
\[
\pz\pzb\log\lambda=\lambda.
\]
By the deformation theory of Kodaira-Spencer, we
know that the tangent space $T_p\M_g$ is identified with $HB(X)$ and the map
$T_p\M_g\to HB(X)$ is given by
\[
\frac{\partial}{\partial s_i}\mapsto A_i\frac{\partial}{\partial z}\otimes d\bar z
\]
where $A_i=-\pzb(\lambda^{-1}\partial_{s_i}\pzb\log\lambda)$.
Similarly the cotangent space $T^\ast_p\M_g$ is identified with
$Q(X)$. For $\mu=\mu(z)\frac{\partial}{\partial z}\otimes d\bar z\in HB(X)$ and
$\phi=\phi(z)dz^2\in Q(X)$, the duality
between them is given by
\[
\langle\mu,\phi\rangle=\int_X \mu(z)\phi(z)\ dzd\bar z
\]
and the Teichm\"uller norm of $\phi$ is defined to be
\[
\Vert\phi\Vert_T=\int_X |\phi(z)|\ dzd\bar z.
\]
By using the above notation, the norm of the Teichm\"uller metric is
given by
\[
\Vert\mu\Vert_T=\sup_{\phi\in Q(X)}
\{\text{Re}\langle\mu,\phi\rangle\mid\Vert\phi\Vert_T=1\}
\]
for all $\mu\in HB(X)\cong T_p\M_g$.

The Weil-Petersson metric on $\M_g$ is defined by
\[
h_{i\bar j}(p)=\int_X A_i\bar{A_j}\ dv
\]
where $dv=\frac{\sqrt{-1}}{2}\lambda dz\wedge d\bar z$ is the volume form of $X$. The Ricci
metric $\tau$ is the negative Ricci curvature of the Weil-Petersson metric
\[
\tau_{i\bar j}=\partial_i\partial_{\bar j}\log\det(h_{k\bar l}).
\]
By the works in \cite{tr1} and \cite{lsy1} we know that the Ricci metric is complete. Now we
take linear combination of the Ricci metric and the Weil-Petersson metric to define the
perturbed Ricci metric
\[
\widetilde\tau_{i\bar j}=\tau_{i\bar j}+Ch_{i\bar j}
\]
where $C>0$. In \cite{lsy1} we proved the following theorem
\begin{theorem}
For suitable choice of large constant $C$, the holomorphic sectional curvature of the
perturbed Ricci metric $\widetilde\tau$ has negative upper
bound. Furthermore, on $\M_g$, the Ricci
metric, the perturbed Ricci metric, the K\"ahler-Einstein metric, the Asymptotic Poincar\'e
metric are equivalent.
\end{theorem}

This theorem was proved by using the curvature properties of the
perturbed Ricci metric and the estimates of its asymptotic
behavior.

Now we prove several claims about the boundedness of the curvature of the Ricci metric
and the perturbed Ricci metric which were stated in \cite{lsy1}.

\begin{theorem}\label{riccibdd}
The holomorphic sectional curvature, bisectional curvature and the Ricci curvature of
the Ricci metric $\tau$ on the moduli space $\M_g$ are bounded.
\end{theorem}

As part of the Theorem 4.4 of \cite{lsy1}, this theorem was
roughly proved in \cite{lsy1}. Here we give a detailed proof since
we need the techniques later.

{\bf Proof.} We follow the notations and computations in \cite{lsy1}. Let $p\in D$ be a
codimension $m$ boundary point and let $(t_1,\cdots,t_m,s_{m+1},\cdots,s_n)$ be the pinching
coordinates of $\M_g$ at $p$ where $t_1,\cdots,t_m$ represent the degeneration directions.
Let $X_{t,s}$ be the Riemann surface corresponding to the point with coordinates
$(t_1,\cdots,s_n)$ and let $l_i$ be the length of the short geodesic loop on the $i$-th collar.
Let $u_i=\frac{l_i}{2\pi}$. We fix $\delta>0$ and assume that $|(t,s)|<\delta$. When
$\delta$ is small enough, from the work of \cite{wol1} and \cite{lsy1} we know that
\[
u_i=-\frac{\pi}{\log|t_i|}\left ( 1+O\left (\left ( \frac{\pi}{\log|t_i|}
\right )^2\right )\right ).
\]
Now we let $u_0=\sum_{i=1}^m u_i+\sum_{j=m+1}^m |s_j|$.

By Corollary 4.2 and Theorem 4.4 of \cite{lsy1}, the work of Masur in \cite{ma1} and Wolpert,
if we use $\widetilde R_{i\bar jk\bar l}$
to denote the curvature tensor of the Ricci metric $\tau$, we have
\begin{enumerate}
\item $\tau_{i\bar i}=\frac{3}{4\pi^2}\frac{u_i^2}{|t_i|^2}(1+O(u_0))$, if $i\leq m$;
\item $\tau_{i\bar j}=O\left (\frac{u_i^2u_j^2}{|t_it_j|}(u_i+u_j)\right )$,
if $i,j\leq m$ and $i\ne j$;
\item $\tau_{i\bar j}=O\left (\frac{u_i^2}{|t_i|}\right )$,
if $i\leq m$ and $j\geq m+1$;
\item $\tau_{i\bar j}=O\left (\frac{u_j^2}{|t_j|}\right )$,
if $j\leq m$ and $i\geq m+1$;
\item $\tau_{i\bar j}=O(1)$, if $i,j\geq m+1$;
\item The matrix $\left ( \tau_{i\bar j}\right )_{i,j\geq m+1}$ is positive definite and has a
positive lower bound depending on $p,n,\delta$;
\item $\widetilde R_{i\bar ii\bar i}=\frac{3u_i^4}{8\pi^4|t_i|^4}(1+O(u_0))$, if $i\leq m$;
\item $\widetilde R_{i\bar ii\bar i}=O(1)$, if $i\geq m+1$.
\end{enumerate}
Now we let
\[
\Lambda_i=
\begin{cases}
\frac{u_i}{|t_i|} & i\leq m\\
1 & i\geq m+1.
\end{cases}
\]
We divide the index set into three parts. Let
\begin{enumerate}
\item $A_1=\{(i,i,i,i)\mid i\leq m\}$;
\item $A_2=\{(i,j,k,l)\mid \text{ at least one of } i,j,k,l \leq m
\text{ and they are not all equal} \}$;
\item $A_3=\{(i,j,k,l)\mid i,j,k,l\geq m+1\}$.
\end{enumerate}
By following the computations of \cite{lsy1} we know that, if
$(i,j,k,l)\in A_2$, then
\[
\widetilde R_{i\bar jk\bar l}=O(\Lambda_i\Lambda_j\Lambda_k\Lambda_l)O(u_0).
\]

Let $v=a_1\frac{\partial}{\partial t_1}+\cdots+a_n\frac{\partial}{\partial s_n}$ and
 $w=b_1\frac{\partial}{\partial t_1}+\cdots+b_n\frac{\partial}{\partial s_n}$ be two
tangent vectors at $(t,s)$. We have
\begin{align*}
\begin{split}
|\widetilde R(v,\bar v,w,\bar w)|=&\left | \sum_{i,j,k,l}a_i\bar a_j b_k\bar b_l
\widetilde R_{i\bar jk\bar l}\right |\leq \sum_{i,j,k,l}
\left | a_i\bar a_j b_k\bar b_l\widetilde R_{i\bar jk\bar l} \right |
=I_1+I_2+I_3
\end{split}
\end{align*}
where $I_\alpha=\sum_{(i,j,k,l)\in A_\alpha}
\left | a_i\bar a_j b_k\bar b_l\widetilde R_{i\bar jk\bar l} \right |$. To estimate the
norms of $v$ and $w$, we have
\[
\tau(v,v)=\sum_{i,j\leq m, i\ne j} a_i\bar a_j\tau_{i\bar j}
+\sum_{i\leq m<j} a_i\bar a_j\tau_{i\bar j}
+\sum_{j\leq m<i} a_i\bar a_j\tau_{i\bar j}
+\sum_{i\leq m} |a_i|^2\tau_{i\bar i}+\sum_{i,j\geq m+1} a_i\bar a_j\tau_{i\bar j}.
\]
By using the asymptotic of $\tau$ and the Schwarz inequality we
have
\[
\left |\sum_{i,j\leq m, i\ne j} a_i\bar a_j\tau_{i\bar j}
+\sum_{i\leq m<j} a_i\bar a_j\tau_{i\bar j}
+\sum_{j\leq m<i} a_i\bar a_j\tau_{i\bar j}\right |
\leq O(u_0)\sum_{i=1}^n |a_i|^2\Lambda_i^2.
\]
Since the matrix $\left ( \tau_{i\bar j}\right )_{i,j\geq m+1}$ has a local
positive lower bound, we know there is a positive constant $c$ depending on $p,n,\delta$
such that
\[
\sum_{i,j\geq m+1} a_i\bar a_j\tau_{i\bar j}\geq c \sum_{i=1}^n |a_i|^2
=c\sum_{i=m+1}^n |a_i|^2\Lambda_i^2.
\]
Finally we have
\[
\sum_{i\leq m} |a_i|^2\tau_{i\bar i}=\frac{3}{4\pi^2}(1+O(u_0))\sum_{i=1}^m |a_i|^2\Lambda_i^2.
\]
By combining the above inequalities we know there is another
positive constant $c_1$ depending on $p,n,\delta$ such that
\begin{eqnarray}\label{taunorm}
\tau(v,v)\geq c_1\sum_{i=1}^n |a_i|^2\Lambda_i^2.
\end{eqnarray}
Similar estimates hold for the Ricci norm of $w$.

Now for each term in $I_2$, by using the Schwarz inequality, we
have
\[
\left | a_i\bar a_j b_k\bar b_l\widetilde R_{i\bar jk\bar l} \right |
=O(u_0)|a_i\bar a_j b_k\bar b_l\Lambda_i\Lambda_j\Lambda_k\Lambda_l|
\leq O(u_0)(|a_i|^2\Lambda_i^2+|a_j|^2\Lambda_j^2)(|b_k|^2\Lambda_k^2+|b_l|^2\Lambda_l^2).
\]
So we have
\begin{eqnarray*}
I_2=O(u_0)\left (\sum_{i=1}^n |a_i|^2\Lambda_i^2\right )
\left (\sum_{i=1}^n |b_i|^2\Lambda_i^2\right )\leq c_0\tau(v,v)\tau(w,w)
\end{eqnarray*}
for some positive constant $c_0$. By enlarging this constant, we also have
\[
I_3\leq O(1)\sum_{(i,j,k,l)\in A_3} | a_i\bar a_j b_k\bar b_l|
\leq O(1) \left (\sum_{i=m+1}^n |a_i|^2\Lambda_i^2\right )
\left (\sum_{i=m+1}^n |b_i|^2\Lambda_i^2\right )\leq c_0\tau(v,v)\tau(w,w)
\]
and
\[
I_1=\frac{3}{8\pi^4}(1+O(u_0))\sum_{i=1}^m |a_i|^2|b_i|^2\Lambda_i^4\leq c_0\tau(v,v)\tau(w,w).
\]
By combining the above inequalities we know that there is a
positive constant $\widetilde c$ depending on $p,n,\delta$ such
that if $\delta$ is small enough, then
\[
|\widetilde R(v,\bar v,w,\bar w)|\leq \widetilde c \tau(v,v)\tau(w,w).
\]
So we have proved that for each point $p\in D$ there is an open
neighborhood $U_p$ such that the bisectional curvature of the
Ricci metric is bounded by a constant which depends on $U_p$.
Since $D$ is compact, we can find a finite cover of $D$ by such
$U_p$. Let $U$ be the union of such $U_p$. Then we can find a
universal constant $c$ which bounds the bisectional curvature at
each point in $U$. Since $\M_g\setminus U$ is a compact set, we
know the bisectional curvature is bounded there. So we proved that
the bisectional curvature of the Ricci metric is bounded.

The boundedness of the holomorphic sectional curvature can be
proved similarly if we replace $w$ by $v$ in the above argument.
Finally since the Ricci curvature is the average of the
bisectional curvature and the holomorphic sectional curvature, it
is bounded. We finish the proof.

\qed

We now investigate the curvatures of the perturbed Ricci metric.
We have
\begin{theorem}\label{prbdd}
For any constant $C>0$, the bisectional curvature of the perturbed Ricci metric
$\widetilde\tau=\tau+Ch$ is bounded. Furthermore, with suitable choice of $C$, the
holomorphic sectional curvature and the Ricci curvature of $\widetilde\tau$ are bounded
from above and below by negative constants.
\end{theorem}

{\bf Proof.} We use $R_{i\bar jk\bar l}$ and $P_{i\bar jk\bar l}$ to denote the curvature
tensor of the Weil-Petersson metric and the perturbed Ricci metric respectively. We
use the same notations as in the proof of the above theorem. Let $p\in D$ be a codimension
$m$ boundary point, let $(t_1,\cdots,s_n)$ be the local pinching coordinates and let
$A_1,A_2,A_3$ be the partition of the index set. Assume $\delta$ is a small positive
constant and $|(t,s)|<\delta$.

By Corollary 4.1, Corollary 4.2 and Theorem 5.2 of \cite{lsy1} we
have
\begin{enumerate}
\item $\widetilde\tau_{i\bar i}=\frac{u_i^2}{|t_i|^2}
\left (\frac{3}{4\pi^2}+\frac{1}{2}Cu_i\right )(1+O(u_0))$, if $i\leq m$;
\item $\widetilde\tau_{i\bar j}=\frac{u_i^2u_j^2}{|t_it_j|}
(O(u_i+u_j)+CO(u_iu_j))$, if $i,j\leq m$ and $i\ne j$;
\item $\widetilde\tau_{i\bar j}=\frac{u_i^2}{|t_i|}(O(1)+CO(u_i) )$,
if $i\leq m$ and $j\geq m+1$;
\item $\widetilde\tau_{i\bar j}=\frac{u_j^2}{|t_j|}(O(1)+CO(u_j) )$,
if $j\leq m$ and $i\geq m+1$;
\item $P_{i\bar ii\bar i}=\left (\left (
\frac{9}{16\pi^4}-\frac{3}{16\pi^4}\left (1+\frac{2\pi^2Cu_i}{3}\right )^{-1}\right )
\frac{u_i^4}{|t_i|^4}
+\frac{3C}{8\pi^2}\frac{u_i^5}{|t_i|^4}\right ) (1+O(u_0))$, if $i\leq m$;
\item $P_{i\bar ii\bar i}=O(1)+CR_{i\bar ii\bar i}$, if $i\geq m+1$;
\item $P_{i\bar jk\bar l}=O(1)+CR_{i\bar jk\bar l}$, if $(i,j,k,l)\in A_3$;
\item $P_{i\bar jk\bar l}=O(\Lambda_i\Lambda_j\Lambda_k\Lambda_l)O(u_0)+CR_{i\bar jk\bar l}$,
if  $(i,j,k,l)\in A_2$
\end{enumerate}
where all the $O$-terms are independent of $C$.

Let $v$ and $w$ be holomorphic vectors as above. To estimate the bisectional curvature, we
have
\begin{align}\label{pr800}
\begin{split}
|P(v,\bar v,w,\bar w)|\leq & \sum_{i=1}^m |a_i|^2|b_i|^2P_{i\bar ii\bar i}
+\sum_{(i,j,k,l)\in A_2}\left | a_i\bar a_j b_k\bar b_l P_{i\bar jk\bar l} \right |
+\sum_{(i,j,k,l)\in A_3}\left | a_i\bar a_j b_k\bar b_l P_{i\bar jk\bar l} \right |\\
\leq & \sum_{i=1}^m |a_i|^2|b_i|^2P_{i\bar ii\bar i}
+O(u_0)\sum_{(i,j,k,l)\in A_2}\left | a_i\bar a_j b_k\bar b_l
\Lambda_i\Lambda_j\Lambda_k\Lambda_l\right |\\
&+C\sum_{(i,j,k,l)\in A_2}\left | a_i\bar a_j b_k\bar b_l R_{i\bar jk\bar l} \right |
+O(1)\sum_{(i,j,k,l)\in A_3}\left | a_i\bar a_j b_k\bar b_l\right |\\
&+C\sum_{(i,j,k,l)\in A_3}\left | a_i\bar a_j b_k\bar b_l R_{i\bar jk\bar l} \right |.
\end{split}
\end{align}
Let $c_i$ denote certain positive constants only depending on
$p,n,\delta$. By the proof of the above theorem, since
$\widetilde\tau\geq \tau$, we have
\begin{eqnarray}\label{pr801}
\left |O(u_0)\sum_{(i,j,k,l)\in A_2}\left | a_i\bar a_j b_k\bar b_l
\Lambda_i\Lambda_j\Lambda_k\Lambda_l\right |\right |\leq c_0 \tau(v,v)\tau(w,w)
\leq c_0\widetilde\tau(v,v)\widetilde\tau(w,w).
\end{eqnarray}
We also have
\begin{eqnarray}\label{pr802}
\left | O(1)\sum_{(i,j,k,l)\in A_3}\left | a_i\bar a_j b_k\bar b_l\right |\right |
\leq c_1 \tau(v,v)\tau(w,w)\leq c_2\widetilde\tau(v,v)\widetilde\tau(w,w).
\end{eqnarray}
Now we estimate $\widetilde\tau(v,v)$. By using similar argument
as in the above proof we know that
\[
\widetilde\tau(v,v)\geq c_2\sum_{i=1}^n |a_i|^2\widetilde\tau_{i\bar i}.
\]
So for $i\leq m$ we have
\begin{align*}
\begin{split}
&|a_i|^2|b_i|^2P_{i\bar ii\bar i}=|a_i|^2|b_i|^2\left (\left (
\frac{9}{16\pi^4}-\frac{3}{16\pi^4}\left (1+\frac{2\pi^2Cu_i}{3}\right )^{-1}\right )
\frac{u_i^4}{|t_i|^4}
+\frac{3C}{8\pi^2}\frac{u_i^5}{|t_i|^4}\right ) (1+O(u_0))\\
\leq &c_3|a_i|^2|b_i|^2\frac{u_i^4}{|t_i|^4}\left (\frac{3}{4\pi^2}+\frac{1}{2}Cu_i\right )^2
\leq c_4 |a_i|^2|b_i|^2\widetilde\tau_{i\bar i}^2
\end{split}
\end{align*}
which implies
\begin{eqnarray}\label{pr803}
\sum_{i=1}^m |a_i|^2|b_i|^2P_{i\bar ii\bar i}\leq c_5
\sum_{i=1}^n |a_i|^2|b_i|^2\widetilde\tau_{i\bar i}^2\leq
c_6\widetilde\tau(v,v)\widetilde\tau(w,w).
\end{eqnarray}
To estimate the rest two terms in the right hand side of
\eqref{pr800} we need the estimates of the curvature tensor of the
Weil-Petersson metric which is done in the proof of Corollary 4.2
of \cite{lsy1}. By collecting the results there we know that
$R_{i\bar jk\bar l}=O(\Lambda_i\Lambda_j\Lambda_k\Lambda_l)O(u_0)$
if $(i,j,k,l)\in A_2$ and $R_{i\bar jk\bar l}=O(1)$ if
$(i,j,k,l)\in A_3$. By using a similar argument as in the above
proof we know that
\begin{eqnarray}\label{pr804}
\left | C\sum_{(i,j,k,l)\in A_2}\left | a_i\bar a_j b_k\bar b_l R_{i\bar jk\bar l} \right |
\right |\leq c_7 C\widetilde\tau(v,v)\widetilde\tau(w,w)
\end{eqnarray}
and
\begin{eqnarray}\label{pr805}
\left |C\sum_{(i,j,k,l)\in A_3}\left | a_i\bar a_j b_k\bar b_l R_{i\bar jk\bar l} \right |
\right |\leq  c_8 C\widetilde\tau(v,v)\widetilde\tau(w,w).
\end{eqnarray}
These imply that
\[
|P(v,\bar v,w,\bar w)|\leq (c_9C+c_{10})\widetilde\tau(v,v)\widetilde\tau(w,w).
\]
By using the compactness argument as above we proved that the bisectional curvature of
$\widetilde\tau$ is bounded. However, the bounds depend on the choice of $C$.

By using a similar method it is easy to see that the holomorphic
sectional curvature is also bounded. However, in \cite{lsy1} we
showed that, for suitable choice of $C$, the holomorphic sectional
curvature has a negative upper bound. So for this $C$, the
holomorphic sectional curvature of $\widetilde\tau$ is pinched
between negative constants.

Finally, we consider the Ricci curvature of
$\widetilde\tau=\tau+C\, h$. We first define two new tensors. Let
$\widehat{R}_{i\bar jk\bar l}=P_{i\bar jk\bar l}-CR_{i\bar jk\bar
l}$ and let $P_{i\bar j}=-Ric(\omega_{\widetilde{\tau}})_{i\bar
j}$. We only need to show that there are positive constants
$\alpha_1$ and $\alpha_2$ which may depend on $C$ such that
\begin{eqnarray}\label{needpr}
\alpha_1\left (\widetilde{\tau}_{i\bar j}\right )\leq \left
(P_{i\bar j}\right )\leq \alpha_2\left (\widetilde{\tau}_{i\bar
j}\right ).
\end{eqnarray}

Based on Lemma 5.2 of \cite{lsy1} and by Corollary 4.1 and 4.2 of
\cite{lsy1} we can estimate the asymptotic of the perturbed Ricci
metric.
\begin{lemma}\label{prest}
Let $p\in \bar\M_g\setminus\M_g$ be a codimension $m$ boundary
point and let $(t,s)=(t_1,\cdots,s_n)$ be the pinching
coordinates. Let $\delta>0$ be a small constant such that
$|(t,s)|<\delta$. Let $C$ be a positive constant. Let $B_1=\left
(\tau_{i\bar j}\right )_{i,j\geq m+1}$, let $B_2=\left (h_{i\bar
j}\right )_{i,j\geq m+1}$ and let $B=B_1+CB_2$. Let $\left
(B^{i\bar j}\right )=\bar{\left (B^{-1}\right )}$ and let
$x_i=2\pi^2Cu_i$ for $i\leq m$. Then we have
\begin{enumerate}
\item $\widetilde\tau_{i\bar i}=
\frac{\Lambda_i^2}{4\pi^{2}}(3+x_i)(1+O(u_0))$ and
$\widetilde\tau^{i\bar
i}=4\pi^2\Lambda_i^{-2}(3+x_i)^{-1}(1+O(u_0))$ if $i \leq m$;
\item $\widetilde\tau_{i\bar
j}=\frac{u_i^2u_j^2}{|t_it_j|}(O(u_i+u_j)+CO(u_iu_j))$ and
$\widetilde\tau^{i\bar
j}=O(|t_it_j|)\min\{(1+x_i)^{-1},(1+x_j)^{-1}\}$ if $i,j\leq m$
and $i\ne j$; \item $\widetilde\tau_{i\bar
j}=\frac{u_i^2}{|t_i|}(O(1)+CO(u_i))$ and $\widetilde\tau^{i\bar
j}=O(|t_i|)(1+x_i)^{-1}$ if $i\leq m$ and $j\geq m+1$; \item
$\widetilde\tau_{i\bar j}=\tau_{i\bar j}+Ch_{i\bar j}$ and
$\widetilde\tau^{i\bar j}=C^{-1}\left ( h^{i\bar
j}+O(C^{-1})+O(u_0)\right )$ if $i,j\geq m+1$; \item
$\widehat{R}_{i\bar jk\bar
l}=O(\Lambda_i\Lambda_j\Lambda_k\Lambda_l)O(u_0)$ if $(i,j,k,l)\in
A_2$; \item $\widehat{R}_{i\bar jk\bar l}=O(1)$ if $(i,j,k,l)\in
A_3$.
\end{enumerate}
\end{lemma}
{\bf Proof.} Let $a=a(t_1,\bar t_1,\cdots,s_n,\bar s_n)$ and
$b=b(t_1,\bar t_1,\cdots,s_n,\bar s_n)$ be any local functions
defined for $|(t,s)|<\delta$. Assume there is local constant $c_1$
depending on $p,\delta$ and $n$ such that
\[
0< c_1\leq a,b
\]
for $|(t,s)|<\delta$. We first realize that there are constants
$\mu_i>0$ depending on $c_1$, $p$, $n$ and $\delta$ such that
\[
1+x_i\leq \mu_i(a+Cb)
\]
for $|(t,s)|<\delta$. In fact, we can pick
$\mu_i=\max\{\frac{1}{c_1},\frac{2\pi^2u_i}{c_1}\}$ since $u_i$ is
small when $\delta$ is small.

The first four claims followed from Corollary 4.1, Corollary 4.2,
Lemma 5.1 and Lemma 5.2 of \cite{lsy1}. By the proof of Lemma 5.2
of \cite{lsy1} we have the linear algebraic formula
\[
\det(\widetilde{\tau})=\left (\prod_{i=1}^m
\frac{\Lambda_i^2}{4\pi^{2}}(3+x_i)\right )\det (B)(1+O(u_0)).
\]
These claims followed from similar computations of the
determinants of the minor matrices.

The last two claims follow from the same techniques and
computations as in the appendix of \cite{lsy1}.

\qed

Now we estimate $P_{i\bar j}$. We first compute $P_{i\bar i}$ with
$i\leq m$. We have
\begin{align*}
\begin{split}
P_{i\bar i}=&\widetilde\tau^{i\bar i}P_{i\bar ii\bar
i}+\sum_{(k,l)\ne (i,i)}\widetilde\tau^{k\bar l}P_{i\bar ik\bar
l}\\
=&\widetilde\tau^{i\bar i}P_{i\bar ii\bar i}+\sum_{(k,l)\ne
(i,i)}\widetilde\tau^{k\bar l}\widehat{R}_{i\bar ik\bar l}
+C\sum_{(k,l)\ne (i,i)}\widetilde\tau^{k\bar l}R_{i\bar ik\bar l}.
\end{split}
\end{align*}
We estimate each term in the right hand side of the above formula.
We have
\begin{align}\label{pr900}
\begin{split}
\widetilde\tau^{i\bar i}P_{i\bar ii\bar
i}=&4\pi^2\Lambda_i^{-2}(3+x_i)^{-1}\left (\left
(\frac{9}{16\pi^4}-\frac{3}{16\pi^4}\left (1+\frac{x_i}{3}\right
)^{-1}\right )\Lambda_i^4+\frac{3x_i}{16\pi^4}\Lambda_i^4\right
)(1+O(u_0))\\
=& \frac{3}{4\pi^2}\Lambda_i^2(3+x_i)^{-1}\left ( 3-\left
(1+\frac{x_i}{3}\right )^{-1}+x_i\right )(1+O(u_0))\\
=& \frac{3}{4\pi^2}\left (1-\frac{3}{(3+x_i)^2}\right
)\Lambda_i^2(1+O(u_0)).
\end{split}
\end{align}
By the fifth claim of the above lemma we have
\begin{align}\label{pr910}
\begin{split}
\left |\sum_{(k,l)\ne (i,i)}\widetilde\tau^{k\bar
l}\widehat{R}_{i\bar ik\bar l}\right |=O(\Lambda_i^2)O(u_0).
\end{split}
\end{align}
We also have
\[
C\sum_{(k,l)\ne (i,i)}\widetilde\tau^{k\bar l}R_{i\bar ik\bar l}=
C\sum_{k\ne i,\ l\ne i}\widetilde\tau^{k\bar l}R_{i\bar ik\bar l}
+C\sum_{k\ne i}\widetilde\tau^{k\bar i}R_{i\bar ik\bar i}
+C\sum_{l\ne i}\widetilde\tau^{i\bar l}R_{i\bar ii\bar l}.
\]
By the proof of Corollary 4.2 of \cite{lsy1} we know that, if
$k\ne i$, then $R_{i\bar ik\bar i}=O\left
(\frac{u_i^5}{|t_i|^3}\right )O(\Lambda_k)$. By combining with the
above lemma, we have
\begin{eqnarray}\label{pr920}
\left |C\sum_{k\ne i}\widetilde\tau^{k\bar i}R_{i\bar ik\bar
i}\right
|=CO(u_i^3)\Lambda_i^2(1+x_i)^{-1}=\frac{x_i}{1+x_i}\Lambda_i^2O(u_i^2)
=\Lambda_i^2O(u_i^2).
\end{eqnarray}
Similarly, we have
\begin{eqnarray}\label{pr930}
\left |C\sum_{l\ne i}\widetilde\tau^{i\bar l}R_{i\bar ii\bar
l}\right |=\Lambda_i^2O(u_i^2).
\end{eqnarray}
Now we fix $i$ and let $v=\frac{\partial}{\partial t_i}$,
$w=b_1\frac{\partial}{\partial
t_1}+\cdots+b_n\frac{\partial}{\partial s_n}$ with $b_i=0$. Since
the bisectional curvature of the Weil-Petersson metric is
non-positive, we have
\[
0\leq R(v,\bar v,w,\bar w)=\sum_{k\ne i,\ l\ne i}b_k\bar
b_lR_{i\bar ik\bar l}.
\]
This implies that the matrix $\left (R_{i\bar ik\bar l}\right
)_{k\ne i,\ l\ne i}$ is semi-positive definite. So we know
\[
\sum_{k\ne i,\ l\ne i}\widetilde\tau^{k\bar l}R_{i\bar ik\bar
l}\geq 0
\]
since it is the trace of the product of a positive definite matrix
and a semi-positive definite matrix. Again, by using the proof of
Lemma 4.2 of \cite{lsy1} we have
\begin{eqnarray}\label{pr940}
0\leq C\sum_{k\ne i,\ l\ne i}\widetilde\tau^{k\bar l}R_{i\bar
ik\bar l}\leq x_i\Lambda_i^2O^+(1)
\end{eqnarray}
where $O^+(1)$ represents a positive bounded term. By combining
formulas \eqref{pr900}, \eqref{pr910}, \eqref{pr920},
\eqref{pr930} and \eqref{pr940} we have
\[
P_{i\bar i}\geq\frac{3}{4\pi^2}\left (1-\frac{3}{(3+x_i)^2}\right
)\Lambda_i^2(1+O(u_0))+\Lambda_i^2O(u_0)+\Lambda_i^2O(u_i^2)
\]
and
\[
P_{i\bar i}\leq \frac{3}{4\pi^2}\left
(1-\frac{3}{(3+x_i)^2}+x_iO^+(1)\right
)\Lambda_i^2(1+O(u_0))+\Lambda_i^2O(u_0)+\Lambda_i^2O(u_i^2)
\]
which imply
\begin{eqnarray}\label{pr950}
\frac{3\Lambda_i^2}{8\pi^2}\left (1-\frac{3}{(3+x_i)^2}\right
)\leq P_{i\bar i}\leq \frac{3\Lambda_i^2}{2\pi^2}\left
(1-\frac{3}{(3+x_i)^2}+x_iO^+(1)\right )
\end{eqnarray}
when $\delta$ is small. The above estimate is independent of the
choice of $C$.

Now we estimate $P_{i\bar j}$ with $i,j\leq m$ and $i\ne j$. We
have
\begin{align}\label{pr960}
\begin{split}
P_{i\bar j}=& \widetilde\tau^{k\bar l}P_{i\bar jk\bar
l}=\widetilde\tau^{k\bar l}\widehat{R}_{i\bar jk\bar l}
+C\widetilde\tau^{k\bar l}R_{i\bar jk\bar l}\\
=&\widetilde\tau^{k\bar l}\widehat{R}_{i\bar jk\bar
l}+C\sum_{k\leq m}\widetilde\tau^{k\bar k}R_{i\bar jk\bar
k}+C\sum_{k,l\leq m,\ k\ne l}\widetilde\tau^{k\bar l}R_{i\bar
jk\bar l}\\
&+C\sum_{k\leq m<l}\widetilde\tau^{k\bar l}R_{i\bar jk\bar
l}+C\sum_{l\leq m<k}\widetilde\tau^{k\bar l}R_{i\bar jk\bar
l}+C\sum_{k,l\geq m+1}\widetilde\tau^{k\bar l}R_{i\bar jk\bar l}.
\end{split}
\end{align}
By the above lemma we have
\begin{eqnarray}\label{pr970}
\left | \widetilde\tau^{k\bar l}\widehat{R}_{i\bar jk\bar l}\right
|=\Lambda_i\Lambda_jO(u_0).
\end{eqnarray}
We also have
\[
\left |C\sum_{k\leq m}\widetilde\tau^{k\bar k}R_{i\bar jk\bar
k}\right |\leq \left |C\sum_{k\leq m,k\ne
i,j}\widetilde\tau^{k\bar k}R_{i\bar jk\bar k}\right |+ \left
|C\widetilde\tau^{i\bar i}R_{i\bar ji\bar i}\right |+\left
|C\widetilde\tau^{j\bar j}R_{i\bar jj\bar j}\right |.
\]
By the proof of Lemma 4.2 of \cite{lsy1} we have $R_{i\bar ji\bar
i}=O\left (\frac{u_i^5u_j^3}{|t_i^3t_j|}\right )$ which implies
\[
\left |C\widetilde\tau^{i\bar i}R_{i\bar ji\bar i}\right
|=4\pi^2\Lambda_i^{-2}(3+x_i)^{-1}CO\left
(\frac{u_i^5u_j^3}{|t_i^3t_j|}\right )=\Lambda_i\Lambda_jO(u_0).
\]
Similarly we have
\[
\left |C\widetilde\tau^{j\bar j}R_{i\bar jj\bar j}\right
|=\Lambda_i\Lambda_jO(u_0).
\]
Again, by the proof of Lemma 4.2 of \cite{lsy1}, for $k\leq m$ and
$k\ne i,j$ we have
\[
R_{i\bar jk\bar k}=O\left
(\frac{u_iu_ju_k^3}{|t_it_jt_k^2|}u_0\right )
\]
which implies
\[
\left |C\sum_{k\leq m,k\ne i,j}\widetilde\tau^{k\bar k}R_{i\bar
jk\bar k}\right |=\Lambda_i\Lambda_jO(u_0).
\]
By combining the above three formulas we have
\[
\left |C\sum_{k\leq m}\widetilde\tau^{k\bar k}R_{i\bar jk\bar
k}\right |=\Lambda_i\Lambda_jO(u_0).
\]
Similarly we can show that
\[
\left | C\sum_{k,l\leq m,\ k\ne l}\widetilde\tau^{k\bar l}R_{i\bar
jk\bar l}\right |=\Lambda_i\Lambda_jO(u_0),
\]
\[
\left | C\sum_{k\leq m<l}\widetilde\tau^{k\bar l}R_{i\bar jk\bar
l}\right |=\Lambda_i\Lambda_jO(u_0)
\]
and
\[
\left | C\sum_{l\leq m<k}\widetilde\tau^{k\bar l}R_{i\bar jk\bar
l}\right |=\Lambda_i\Lambda_jO(u_0).
\]
Finally,
\[
\left |C\sum_{k,l\geq m+1}\widetilde\tau^{k\bar l}R_{i\bar jk\bar
l}\right |\leq \sum_{k,l\geq m+1}\left |C\widetilde\tau^{k\bar
l}\right | \left |R_{i\bar jk\bar l}\right |=\sum_{k,l\geq
m+1}O(1) \left |R_{i\bar jk\bar l}\right |
=\Lambda_i\Lambda_jO(u_0).
\]
By combining the above results we have
\begin{eqnarray}\label{pr980}
P_{i\bar j}=\Lambda_i\Lambda_jO(u_0).
\end{eqnarray}
By using the same method we know that, if $i\leq m < j$, then
\begin{eqnarray}\label{pr990}
P_{i\bar j}=\Lambda_iO(u_0)
\end{eqnarray}
and if $j\leq m < i$, then
\begin{eqnarray}\label{pr1000}
P_{i\bar j}=\Lambda_jO(u_0).
\end{eqnarray}

The next step is to estimate the matrix $\left (P_{i\bar j}\right
)_{i,j\geq m+1}$. We will show that this matrix is bounded from
above and below by positive constant multiples of the matrix $B_1$
defined in the above lemma where the constants depend on $\delta$
and $C$. We first estimate $P_{i\bar j}$ with fixed $i,j\geq m+1$.
We have
\begin{align*}
\begin{split}
P_{i\bar j}=&\widetilde\tau^{k\bar l}P_{i\bar jk\bar
l}=\widetilde\tau^{k\bar l}\widehat{R}_{i\bar jk\bar
l}+\widetilde\tau^{k\bar l}R_{i\bar jk\bar l}\\
=&\sum_{k\leq m}\widetilde\tau^{k\bar k}\widehat{R}_{i\bar jk\bar
k}+\sum_{k,l\leq m,\ k\ne l}\widetilde\tau^{k\bar
l}\widehat{R}_{i\bar jk\bar l}+\sum_{k\leq
m<l}\widetilde\tau^{k\bar l}\widehat{R}_{i\bar jk\bar
l}+\sum_{l\leq m<k}\widetilde\tau^{k\bar l}\widehat{R}_{i\bar
jk\bar l}\\
&+\sum_{k,l\geq m+1}\widetilde\tau^{k\bar l}\widehat{R}_{i\bar
jk\bar l}+C\sum_{k\leq m}\widetilde\tau^{k\bar k}R_{i\bar jk\bar
k}+C\sum_{k,l\leq m,\ k\ne l}\widetilde\tau^{k\bar l}R_{i\bar
jk\bar l}\\
&+C\sum_{k\leq m<l}\widetilde\tau^{k\bar l}R_{i\bar jk\bar
l}+C\sum_{l\leq m<k}\widetilde\tau^{k\bar l}R_{i\bar jk\bar l}
+C\sum_{k,l\geq m+1}\widetilde\tau^{k\bar l}R_{i\bar jk\bar l}.
\end{split}
\end{align*}
By Lemma \ref{prest} and the proof of Corollary 4.2 of \cite{lsy1}
we know that $\left | \sum_{k\leq m}\widetilde\tau^{k\bar
k}\widehat{R}_{i\bar jk\bar k}\right |=O(u_0)$, $\left
|\sum_{k,l\leq m,\ k\ne l}\widetilde\tau^{k\bar
l}\widehat{R}_{i\bar jk\bar l} \right |=O(u_0)$, $\left
|\sum_{k\leq m<l}\widetilde\tau^{k\bar l}\widehat{R}_{i\bar jk\bar
l} \right |=O(u_0)$, $\left |\sum_{l\leq m<k}\widetilde\tau^{k\bar
l}\widehat{R}_{i\bar jk\bar l} \right |=O(u_0)$, \\
$\left |
\sum_{k,l\geq m+1}\widetilde\tau^{k\bar l}\widehat{R}_{i\bar
jk\bar l}\right |=O(C^{-1})$, $\left | C\sum_{k,l\leq m,\ k\ne
l}\widetilde\tau^{k\bar l}R_{i\bar jk\bar l}\right |=O(u_0)$,
$\left | C\sum_{k\leq m<l}\widetilde\tau^{k\bar l}R_{i\bar jk\bar
l}\right |=O(u_0)$ and $\left | C\sum_{l\leq
m<k}\widetilde\tau^{k\bar l}R_{i\bar jk\bar l}\right |=O(u_0)$.
Also, since for $i,j,k,l\geq m+1$, $R_{i\bar jk\bar l}=O(1)$, we
have
\[
C\sum_{k,l\geq m+1}\widetilde\tau^{k\bar l}R_{i\bar jk\bar
l}=\sum_{k,l\geq m+1}h^{k\bar l}R_{i\bar jk\bar
l}+O(C^{-1})+O(u_0).
\]
By combining the above arguments we have
\begin{eqnarray}\label{pr1010}
P_{i\bar j}=C\sum_{k\leq m}\widetilde\tau^{k\bar k}R_{i\bar jk\bar
k}+\sum_{k,l\geq m+1}h^{k\bar l}R_{i\bar jk\bar
l}+O(C^{-1})+O(u_0).
\end{eqnarray}

The matrix $\left ( \sum_{k,l\geq m+1}h^{k\bar l}R_{i\bar jk\bar
l}(0,s)\right )_{i,j\geq m+1}$ is just the negative of the Ricci
curvature matrix of the restriction of the Weil-Petersson metric
to the boundary piece. So we know it is positive definite and is
bounded from below by a constant multiple of $B_2(0,s)$. By
continuity we know that, when $\delta$ is small enough, the matrix
$\left ( \sum_{k,l\geq m+1}h^{k\bar l}R_{i\bar jk\bar
l}(t,s)\right )_{i,j\geq m+1}$ is bounded from below by a constant
multiple of $B_2(t,s)$. Again, since $h^{k\bar l}R_{i\bar jk\bar
l}=O(1)$ when $i,j,k,l\geq m+1$ and the fact that matrices $B_1$
and $B_2$ are locally equivalent, we know that $\left (
\sum_{k,l\geq m+1}h^{k\bar l}R_{i\bar jk\bar l}\right )_{i,j\geq
m+1}$ is locally bounded from above and below by positive
constants multiples of $B_1$.

Finally, by using the fact that the bisectional curvature of the
Weil-Petersson metric is non-positive and $\widetilde\tau^{k\bar
k}>0$, we know that the matrix $\left (\widetilde\tau^{k\bar
k}R_{i\bar jk\bar k}\right )_{i,j\geq m+1}$ is positive
semi-definite. Also, we know that $C\sum_{k\leq
m}\widetilde\tau^{k\bar k}R_{i\bar jk\bar k}=O(1)$.

Now by using formula \eqref{pr1010} we know that there are
positive constants $\beta_1\leq \beta_2$ depending on $\delta$,
the point $p$ and the choice of $C$ such that as long as $\delta$
is small enough and $C$ is large enough,
\begin{eqnarray}\label{pr1020}
\beta_1 B_1\leq \left (P_{i\bar j}\right )_{i,j\geq m+1}\leq
\beta_2 B_1.
\end{eqnarray}

We know that there is a constant $c_0>0$ such that $h\leq c_0\tau$
which implies $\tau\leq\widetilde\tau\leq (1+c_0C)\tau$. By
combining formulas \eqref{pr950}, \eqref{pr980}, \eqref{pr990},
\eqref{pr1000} and \eqref{pr1020} we know that, when $\delta$ is
small enough and $C$ is large enough, there are positive constants
$\alpha_1\leq \alpha_2$ depending on $p$, $\delta$ and $C$ such
that
\[
\alpha_1\widetilde\tau\leq \left (P_{i\bar j}\right )\leq
\alpha_2\widetilde\tau.
\]
Now by using the compactness argument as we did before, we can
find an open neighborhood $U$ of $D$ in $\bar\M_g$ and a $C_0>0$
such that
\[
\alpha_1\widetilde\tau\leq \left (P_{i\bar j}\right )\leq
\alpha_2\widetilde\tau
\]
on $U$ for positive constants $\alpha_1$ and $\alpha_2$ as long as
$C\geq C_0$.

Let $V=\bar\M_g\setminus U$. We know $V$ is compact. We also know
that, for $C$ large enough,
\[
Ric(\widetilde\tau)=Ric(C^{-1}\widetilde\tau)=Ric(h+C^{-1}\tau).
\]
Since the Ricci curvature of the Weil-Petersson metric has a
negative upper bound, a perturbation of the Weil-Petersson metric
with a small error term still has negative Ricci curvature on a
compact set $V$. Also, on $V$ the perturbed Ricci metric is
bounded. So for $C$ large, we know that the Ricci curvature of the
perturbed Ricci metric is pinched between negative constant
multiples of the perturbed Ricci metric. Here the bounds depend on
the choice of $C$. This finished the proof.

\qed

As a direct corollary of the above theorem, we show a vanishing
theorem similar to the work of Faltings \cite{fa}.

\begin{cor}
Let $D$ be the compactification divisor of the Deligne-Mumford
compactification of $\M_g$. Then
\[
H^0\left (\bar\M_g, \Omega\left (\log D\right )^\ast\right )=0.
\]
\end{cor}

{\bf Proof.} We first pick a constant $C>0$ such that the Ricci
curvature of the perturbed Ricci metric $\widetilde\tau=\tau+Ch$
is pinched by negative constants. Let $\sigma$ be a holomorphic
section of $\Omega\left (\log D\right )^\ast$.

Let $p\in \bar{\M_g}\setminus \M_g$ be a codimension $m$ point and
let $(t_1,\cdots,t_m,s_{m+1},\cdots,s_n)$ be local pinching
coordinates. Then locally we have
\[
\sigma=\sum_{i=1}^m a_i(t,s)t_i\frac{\partial}{\partial t_i}
+\sum_{j=m+1}^n a_j(t,s)\frac{\partial}{\partial s_j}
\]
where $a_i$ are bounded local holomorphic functions for $1\leq
i\leq n$. It is clear that, restricted to $\M_g$, $\sigma$ is a
holomorphic vector field. Now we equip the moduli space $\M_g$
with the perturbed Ricci metric $\widetilde\tau$. From the above
expression of $\sigma$, it is easy to see that
\[
\Vert\sigma\Vert_{\widetilde\tau}\in L^2(\M_g,\widetilde\tau)
\]
since $\widetilde\tau$ is equivalent to the asymptotic Poincar\'e
metric. Now we have the Bochner formula
\[
\Delta_{\widetilde\tau}\Vert\sigma\Vert_{\widetilde\tau}^2
=\Vert\nabla\sigma\Vert_{\widetilde\tau}^2
-Ric_{\widetilde\tau}(\sigma,\sigma).
\]
To integrate, we need a special cut-off function. In \cite{ls2}, a
monotone sequence of cut-off functions $\rho_{\epsilon}$ with the
properties that $\Delta_{\widetilde\tau}\rho_\epsilon$ is
uniformly bounded for each $\epsilon$ and the measure of the
support of $\Delta_{\widetilde\tau}\rho_\epsilon$ goes to $0$ as
$\epsilon$ goes to zero. We will recall the construction in the
next section.

By using the cut-off function $\rho_\epsilon$ we have
\[
\lim_{\epsilon\to 0}\int_{\M_g}\rho_\epsilon
\Delta\Vert\sigma\Vert_{\widetilde\tau}^2
dV_{\widetilde\tau}=\lim_{\epsilon\to
0}\int_{\M_g}\Delta\rho_\epsilon\Vert\sigma\Vert_{\widetilde\tau}^2
dV_{\widetilde\tau}=0
\]
since $\sigma$ is an $L^2$ section with respect to
$\widetilde\tau$ and the measure of
$\Delta_{\widetilde\tau}\rho_\epsilon$ goes to $0$. This above
formula implies $Ric(\sigma,\sigma)=0$ since $Ric(\widetilde\tau)$
is negative which implies $\sigma=0$. Thus we have proved the
corollary.

\qed

Finally, we show that the Teichm\"uller space equipped with the
Ricci metric or the perturbed Ricci metric has bounded geometry.

\begin{cor}\label{injectivity}
The injectivity radius of the Teichm\"uller space equipped with
the Ricci metric or the perturbed Ricci metric is bounded from below.
\end{cor}

{\bf Proof.} We only prove that there is a lower bound for the
injectivity radius of the Ricci metric since the case of the
perturbed Ricci metric can be done in the same way.

In Theorem \ref{riccibdd} we showed that the curvature of the
Ricci metric is bounded. We denote the sup of the curvature by
$\delta$. If $\delta\leq 0$ then the injectivity radius is $+\infty$
by the Cartan-Hadamard theorem. Now we assume $\delta>0$.

Assume the injectivity radius of $(\mathcal T_g, \tau)$ is $0$, then
for any $\epsilon>0$, there is a point $p=p_\epsilon$ such that the
injectivity radius at $p$ is less than $\epsilon$.

Let $f_p$ be the Bers' embedding map such that $f_p(p)=0$. By
using a similar argument as in the proof of Theorem \ref{bck}, and
by changing some constants, we know that the Ricci metric and the
Euclidean metric are equivalent on the Euclidean ball
$B_{1}\subset f_p(\mathcal T_g)$. By using the Rauch comparison
theorem to compare the Ricci metric on the ball $B_1$ and the
standard sphere of constant curvature $\delta$, we know that there
is no conjugate point of $p$ within distance $\epsilon$ when
$\epsilon$ is small enough.

So the only case we need to rule out is that there is a closed
geodesic loop $\gamma$ containing $p$ such that
$l_\tau(\gamma)\leq 2\epsilon$. We know that when $\epsilon$ small
enough, $\gamma\subset B_1$ since the Ricci metric and the
Euclidean metric are equivalent on $B_1$. This implies that the
Euclidean length of $\gamma$, denoted by $\widetilde l(\gamma)\leq
c\epsilon$ for some constant $c$ only depending on the comparison
constants of the Ricci metric and the Euclidean metric on $B_1$. It
is clear that $\gamma$ bounds a minimal disk $\widetilde\Sigma$
with respect to the Euclidean metric. By the isoperimetric
inequality, we know that the Euclidean area $A_{\sp
E}(\widetilde\Sigma)$ satisfies
\[
A_{\sp E}(\widetilde\Sigma)\leq c_1\widetilde l(\gamma)^2
\leq c_1c^2\epsilon^2.
\]
By using the equivalence of the metrics, we know the area
$A_\tau(\widetilde\Sigma)$ of the surface $\widetilde\Sigma$ under
the Ricci metric is small if $\epsilon$ is small enough. Thus
$\gamma$ bounds a minimal disk $\Sigma$ with respect to the Ricci
metric. By the Gauss-Codazzi equation we know that the curvature
$R_{\sp\Sigma}$ of the metric on $\Sigma$ induced from the Ricci
metric is bounded above by $\delta$. By using the isoperimetric
inequality, we know that
\[
A_\tau(\Sigma)\leq c_2\epsilon^2.
\]
However, by the Gauss-Bonnet theorem, since the geodesic $\gamma$
has at most one vertex $p$ and the outer angle $\theta$ at $p$ is at
most $\pi$, we have
\[
\int_\Sigma R_{\Sigma}\ dv_\tau+\int_\gamma \kappa_\gamma\ ds+\theta
=2\pi\chi(\Sigma)=2\pi
\]
where $dv_\tau$ is the induced area form from the Ricci metric. Since
$\gamma$ is a geodesic, we see that the second term in the left hand
side of the above formula is $0$. Since $R_\Sigma\leq\delta$ and
$\theta<\pi$, we have
\[
\delta A_\tau(\Sigma)\geq\int_{\Sigma}R_{\sp\Sigma}\ dv_\tau
=2\pi-\theta>\pi
\]
which implies $A_\tau(\Sigma)\geq \frac{\pi}{\delta}$. By
comparing the above two inequalities, we get a contradiction as
long as $\epsilon$ is small enough. This finishes the proof.

\qed

\section{The Stability of the Logarithmic Cotangent Bundle}\label{stab}

In this section we investigate the cohomology classes defined by
the currents $\omega_\tau$ and $\omega_{\sp{}KE}$. Since both of
these K\"ahler forms have Poincar\'e growth, it is natural to identify
them with the first Chern class of the logarithmic cotangent
bundle of $\bar\M_g$. This implies this bundle is positive over
the compactified moduli space which directly implies that the
moduli space is of log general type.

The next step is to show that the restriction of the
K\"ahler-Einstein metric to a subbundle of the logarithmic
cotangent bundle will not have growth worse than Poincar\'e
growth. Then we prove that the logarithmic cotangent bundle $\bar
E$ over $\bar\M_g$ is stable with respect to the first Chern class
of this bundle.

More precisely we have the
following theorem:
\begin{theorem}\label{stablemain}
The first Chern class of $\bar E$ is positive and $\bar E$ is Mumford
stable with respect to $c_1(\bar E)$.
\end{theorem}

We first setup our notation. On $\M_g$, let $g_p$, $\tau$, $\widetilde\tau$, $g_{_{WP}}$
and $g_{_{KE}}$ be the asymptotic Poincar\'e metric, the Ricci metric, the
perturbed Ricci metric, the Weil-Petersson metric and the K\"ahler-Einstein metric
respectively. Let
$\omega_p$, $\omega_\tau$ and $\omega_{_{KE}}$ be the corresponding K\"ahler
forms of these metrics. Let $Ric(\omega_\tau)$ be the Ricci form of the Ricci
metric.

Let $D=\bar\M_g\setminus\M_g$ be the compactification divisor.
In order to prove the stability, we need to control the growth of these
K\"ahler forms near $D$. We fix a cover of $\M_g$ by local charts.

For each point $y \in D$, we can pick local pinching coordinate charts
$U_y\subset\widetilde U_y$ centered at $y$ with
$U_y=(\Delta^\ast_{\delta_y})^{m_y}\times \Delta_{\delta_y}^{n-m_y}$ and
$\widetilde U_y=(\Delta^\ast_{2\delta_y})^{m_y}\times \Delta_{2\delta_y}^{n-m_y}$
such that the estimates in Corollary 4.1, Corollary 4.2, Theorem 4.4 and
Theorem 5.2 of \cite{lsy1} hold on $\widetilde U_y$. Here $\Delta_{\delta_y}$ is
the disk of radius $\delta_y>0$ and  $\Delta_{\delta_y}^\ast$ is the punctured disk
of radius $\delta_y$ and $m_y$ is the codimension of the point $y$ and $n=3g-3$ is
the complex dimension of $\M_g$.

Since $D$ is compact, we can find $p$ such
charts $U_1=U_{y_1},\cdots,U_p=U_{y_p}$ such that there is a neighborhood $V_0$ of $D$ with
$D\subset V_0\subset \bar{V_0}\subset \cup_{i=1}^p U_i$.
Now we choose coordinate charts $V_1,\cdots,V_q$ such that the estimate of Theorem 5.2
of \cite{lsy1} hold and
\begin{enumerate}
\item $\M_g\subset \left ( \cup_{j=1}^p U_j \right ) \cup \left
(\cup_{j=1}^q V_j\right )$; \item $\left (\cup_{j=1}^q V_j \right
)\cap \bar{V_0}=\emptyset$.
\end{enumerate}
Let $\psi_1,\cdots,\psi_{p+q}$ be a partition of unity subordinate to the
cover $U_1,\cdots,U_p,V_1,\cdots,V_q$ such that $supp (\psi_i)\subset U_i$ for
$1\leq i\leq p$ and $supp (\psi_i)\subset V_{i-p}$ for $p+1\leq i\leq p+q$.

Let $\alpha_i=m_{y_i}$ and let
$t_1^i,\cdots,t_{\alpha_i}^i,s_{\alpha_i+1}^i,\cdots,s_n^i$ be the pinching
coordinates on $U_i$ where $t_1^i,\cdots,t_{\alpha_i}^i$ represent the degeneration
directions.

To prove the theorem, we need a special cut-off function. Such function was used in
\cite{ls2}. We include a short proof here since we need to use the construction later.

\begin{lemma}\label{cutoff}
For any small $\epsilon>0$ there is a smooth function
$\rho_{\epsilon}$ such that
\begin{enumerate}
\item $0\leq \rho_{\epsilon}\leq 1$; \item For any open
neighborhood $V$ of $D$ in $\bar\M_g$, there is a $\epsilon>0$
such that $supp (1-\rho_{\epsilon})\subset V$; \item For each
$\epsilon>0$, there is a neighborhood $W$ of $D$ such that
$\rho_{\epsilon}\mid_W\equiv 0$; \item $\rho_{\epsilon'}\geq
\rho_{\epsilon}$ if $\epsilon'\leq \epsilon$; \item There is a
constant $C$ which is independent of $\epsilon$ such that
\[
-C\omega_p\leq \sqrt{-1}\partial\bar\partial\rho_{\epsilon}\leq
C\omega_p.
\]
\end{enumerate}
\end{lemma}

{\bf Proof.}
We fix a smooth function $\phi\in C^{\infty}(\mathbb R)$ with $0\leq \phi\leq 1$ such
that
\begin{eqnarray*}
\phi(x)=
\begin{cases}
0\ \ \  x\geq 1;\\
1\ \ \  x\leq 0.
\end{cases}
\end{eqnarray*}
Now let
\[
\phi_\epsilon (z)=\phi\left (\frac{\left (\log\frac{1}{|z|}\right )^{-1}-\epsilon}
{\epsilon}\right ).
\]
For $1\leq i \leq p$ and $\epsilon>0$ small, we let
\[
\phi_\epsilon^i(t_1^i,\cdots,s_n^i)=\prod_{j=1}^{\alpha_i}(1-\phi_\epsilon(t_j^i)).
\]
The cut-off function is defined by
\[
\rho_\epsilon=1-\sum_{i=1}^p\psi_i\phi_\epsilon^i.
\]
It is easy to check that $\rho_\epsilon$ satisfy all the conditions.

\qed

Now we discuss the logarithmic cotangent bundle. Let
$U_1,\cdots,V_q$ be the cover of $\M_g$ as above. For each
$1\leq i\leq p$, let $W_i=\Delta_{\delta_{y_i}}^{m_{y_i}}\times
\Delta_{\delta_{y_i}}^{n-m_{y_i}}$. Then
$W_1,\cdots,W_p,V_1,\cdots,V_q$ is a cover of $\bar\M_g$. On each
$U_i$, a local holomorphic frame of the holomorphic cotangent
bundle $T^\ast\M_g$ is
$dt_1^i,\cdots,dt_{\alpha_i}^i,ds_{\alpha_i+1}^i,\cdots,ds_n^i$.
Let
\begin{eqnarray}\label{logframe}
e_j^i=
\begin{cases}
\frac{dt_j^i}{t_j^i}& j\leq \alpha_i;\\
ds_j^i & j\geq \alpha_i+1.
\end{cases}
\end{eqnarray}
The logarithmic cotangent bundle $\bar E$ is the extension of
$T^\ast\M_g$ to $\bar\M_g$ such that on each $U_i$,
$e_1^i,\cdots,e_n^i$ is a local holomorphic frame of $\bar E$. It
is very easy to check this fact by writing down the transition
maps. In the following, we will use $g_{_{WP}}^\ast$, $\tau^\ast$
and $g_{_{KE}}^\ast$ ro represent the metrics on $\bar E$ induced
by the Weil-Petersson metric, the Ricci metric and the
K\"ahler-Einstein metric respectively.

To discuss the stability of $\bar E$, we need to fix a K\"ahler class on $\bar\M_g$. It is
natural to use the first Chern class of $\bar E$. we denote this class by $\Phi$.
We first identify the current represented
by the K\"ahler form $\omega_{_{KE}}$ with $\Phi$.

\begin{lemma}
The currents $\omega_\tau$ and $\omega_{_{KE}}$ are positive closed currents.
Furthermore,
\[
[\omega_\tau]=[\omega_{_{KE}}]=c_1(\bar E).
\]
\end{lemma}

{\bf Proof.} It is clear that $\omega_\tau$ and $\omega_{_{KE}}$ are positive currents. Let
$\phi$ be an arbitrary smooth $(2n-3)$-form on $\bar\M_g$.
To show that $\omega_{_{KE}}$ is closed, we only need to show
\begin{eqnarray}\label{c1}
\int_{\bar\M_g}\omega_{_{KE}}\wedge d\phi=0.
\end{eqnarray}
We first check
\begin{eqnarray}\label{c2}
\int_{\bar\M_g}\left |\omega_{_{KE}}\wedge d\phi\right |
=\int_{\M_g}\left |\omega_{_{KE}}\wedge d\phi\right |<\infty.
\end{eqnarray}
To simplify the notations, on each $U_i$, we let $t_j^i=s_j^i$ for $\alpha_i+1\leq j\leq n$.
On each $U_i$ we assume
\[
d\phi=\sum_{\alpha,\beta}a_{\alpha\bar\beta}^i dt_1^i\wedge d\bar{t_1^i}\cdots\wedge
\widehat{dt_\alpha^i}\wedge d\bar{t_\alpha^i}\cdots\wedge
dt_\beta^i\wedge\widehat{d\bar{t_\beta^i}}\cdots\wedge dt_n^i\wedge d\bar{t_n^i}
\]
where $a_{\alpha\bar\beta}^i$ are bounded smooth functions on
$U_i$. We denote $dt_1^i\wedge d\bar{t_1^i}\cdots\wedge
dt_n^i\wedge d\bar{t_n^i}$ by $dt^i\wedge d\bar{t^i}$. By
\cite{lsy1} we know that the K\"ahler-Einstein metric is
equivalent to the Ricci metric and the asymptotic Poincar\'e
metric. By using Corollary 4.2 of \cite{lsy1}, we know that,
restricted to each $U_i$, there is a constant $C$ depending on
$\phi$ such that
\begin{align*}
\begin{split}
|\omega_{_{KE}}\wedge d\phi|\leq &\left (\frac{\sqrt{-1}}{2}\right )^n
\sum_{\alpha,\beta}\left |
\left (g_{_{KE}}\right )_{\alpha\bar\beta}a_{\alpha\bar\beta}^i\right |
dt^i\wedge d\bar{t^i}\\
\leq &  \left (\frac{\sqrt{-1}}{2}\right )^n C
\left (\sum_{j=1}^{\alpha_i}\frac{1}{\left |t_j^i \right |^2
\left (\log\left |t_j^i \right |\right )^2}+1\right )dt^i\wedge d\bar{t^i}
\end{split}
\end{align*}
since $a_{\alpha\bar\beta}^i$ is bounded on $U_i$. This implies
\[
\int_{U_i}\psi_i|\omega_{_{KE}}\wedge d\phi|\leq \int_{U_i}\psi_i
\left (\frac{\sqrt{-1}}{2}\right )^n C
\left (\sum_{j=1}^{\alpha_i}\frac{1}{\left |t_j^i \right |^2
\left (\log\left |t_j^i \right |\right )^2}+1\right )dt^i\wedge d\bar{t^i}
<\infty.
\]
So we have
\[
\int_{\M_g}\left |\omega_{_{KE}}\wedge d\phi\right | \leq
\sum_{j=1}^q\int_{V_j}\psi_{j+p} \left |\omega_{_{KE}}\wedge
d\phi\right | +\sum_{j=1}^p\int_{U_j}\psi_{j} \left
|\omega_{_{KE}}\wedge d\phi\right |<\infty.
\]
Let $\rho_\epsilon$ be the cut-off function constructed above. By
the dominating convergence theorem, we have
\begin{eqnarray}\label{c3}
\int_{\M_g}\omega_{_{KE}}\wedge d\phi=\lim_{\epsilon\to 0}
\int_{\M_g}\rho_\epsilon\omega_{_{KE}}\wedge d\phi.
\end{eqnarray}
Let $h$ be an Hermitian metric on $\bar E$. Let
$Ric(h)=-\partial\bar\partial\log\det(h)$ be its Ricci form.
Clearly,
\[
[Ric(h)]=c_1(\bar E)
\]
and
\[
\lim_{\epsilon\to 0}
\int_{\M_g}\rho_\epsilon Ric(h)\wedge d\phi=\int_{\M_g} Ric(h)\wedge d\phi
=\int_{\bar\M_g} Ric(h)\wedge d\phi=-\int_{\bar\M_g}d(Ric(h))\wedge \phi=0.
\]
Since $\omega_{_{KE}}=-\partial\bar\partial\log\det(g_{_{KE}}^\ast)$, we have
\begin{align}\label{c4}
\begin{split}
\int_{\M_g}\omega_{_{KE}}\wedge d\phi=& \lim_{\epsilon\to 0}
\int_{\M_g}\rho_\epsilon\omega_{_{KE}}\wedge d\phi
=\lim_{\epsilon\to 0}
\int_{\M_g}\rho_\epsilon\omega_{_{KE}}\wedge d\phi-\lim_{\epsilon\to 0}
\int_{\M_g}\rho_\epsilon Ric(h)\wedge d\phi\\
=&\lim_{\epsilon\to 0}\int_{\M_g}\rho_\epsilon \partial\bar\partial\log\left (
\frac{\det(h)}{\det(g_{_{KE}}^\ast)}\right )\wedge d\phi\\
=& -\lim_{\epsilon\to 0}\int_{\M_g}\log\left (
\frac{\det(g_{_{KE}}^\ast)}{\det(h)}\right )
\partial\bar\partial\rho_\epsilon \wedge d\phi.
\end{split}
\end{align}
By using the frame in \eqref{logframe},  by Theorem 1.4 and
Corollary 4.2 of \cite{lsy1} we know that there are positive
constants $C_1$ and $C_2$ which may depend on $\phi$ such that, on
each $U_i$,
\[
C_1\leq \frac{\det(g_{_{KE}}^\ast)}{\det(h)}\leq C_2 \prod_{j=1}^{\alpha_i}
\left (\log\left |t_j^i\right |\right )^2
\]
which implies that there is a constant $C_3$ such that
\begin{eqnarray}\label{c5}
\left |\log\left ( \frac{\det(g_{_{KE}}^\ast)}{\det(h)}\right )\right |
\leq C_3+2\sum_{j=1}^{\alpha_i}\log\log\frac{1}{\left |t_j^i\right |}.
\end{eqnarray}
Now by Lemma \ref{cutoff} we can pick $\epsilon_0$ small enough such that for any
$0<\epsilon<\epsilon_0$, $supp(1-\rho_\epsilon)\subset V_0$. We also know
that $supp(\partial\bar\partial\rho_\epsilon)\subset supp(1-\rho_\epsilon)$. Again, by
Lemma \ref{cutoff}, since
\[
-C\omega_p\leq \partial\bar\partial\rho_\epsilon
\leq C\omega_p
\]
we know that there is a constant $C_4$ which depends on $\phi$ such that
\begin{eqnarray}\label{c6}
\left |\partial\bar\partial\rho_\epsilon\wedge d\phi  \right |\leq
\left (\frac{\sqrt{-1}}{2}\right )^n C_4
\left (\sum_{j=1}^{\alpha_i}\frac{1}{\left |t_j^i \right |^2
\left (\log\left |t_j^i \right |\right )^2}+1\right )dt^i\wedge d\bar{t^i}.
\end{eqnarray}
By combining \eqref{c5} and \eqref{c6}, and a simple computation
we can show that
\begin{eqnarray}\label{c7}
\int_{U_i}\left |\log\left (
\frac{\det(g_{_{KE}}^\ast)}{\det(h)}\right )
\partial\bar\partial\rho_\epsilon \wedge d\phi\right |<\infty.
\end{eqnarray}
From the above argument we know that
\begin{align*}
\begin{split}
&\left |\int_{\M_g}\log\left (
\frac{\det(g_{_{KE}}^\ast)}{\det(h)}\right )
\partial\bar\partial\rho_\epsilon \wedge d\phi\right |
=\left |\int_{\M_g\cap supp(\partial\bar\partial\rho_\epsilon)}\log\left (
\frac{\det(g_{_{KE}}^\ast)}{\det(h)}\right )
\partial\bar\partial\rho_\epsilon \wedge d\phi\right |\\
\leq & \sum_{j=1}^{p}\int_{U_i\cap supp(1-\rho_{\epsilon})}\psi_i
\left |\log\left (
\frac{\det(g_{_{KE}}^\ast)}{\det(h)}\right )
\partial\bar\partial\rho_\epsilon \wedge d\phi\right |\to 0
\end{split}
\end{align*}
as $\epsilon$ goes to $0$ because of \eqref{c7} and the fact that
the Lebesgue measure of $U_i\cap supp(1-\rho_{\epsilon})$ goes to
$0$. By combining with \eqref{c4} we know that
\[
\int_{\M_g}\omega_{_{KE}}\wedge d\phi=0
\]
which implies $\omega_{_{KE}}$ is a closed current. Similarly we can prove that
$\omega_\tau$ is a closed current by the formula
\[
\omega_{\tau}=-\partial\bar\partial\log\det(g_{_{WP}}^\ast)
\]
and Corollary 4.1 and 4.2 of \cite{lsy1}.

Now we prove the second statement of the lemma. Since
$\omega_{_{KE}}$ is a closed current, to show it represents the
first Chern class of $\bar E$, we need to prove that for any
closed $(n-1,n-1)$-form $\widetilde\phi$ on $\bar\M_g$,
\begin{eqnarray}\label{c8}
\int_{\M_g}\omega_{_{KE}}\wedge\widetilde\phi=\int_{\M_g}c_1(\bar E)\wedge\widetilde\phi.
\end{eqnarray}
However, this can be easily proved by using the above argument
where we replace $d\phi$ by $\widetilde\phi$. The same argument
works for $\omega_\tau$. This finishes the proof.

\qed

Now we compute the degree of $\bar E$. In the following, by degree of a bundle over
$\bar\M_g$ we always mean the $\Phi$-degree.

\begin{lemma}\label{amdegree}
The degree of $\bar E$ is given by $\int_{\M_g}\omega_{_{KE}}^n$.
\end{lemma}

{\bf Proof.} Since the degree of $\bar E$ is given by
\[
\deg(\bar E)=\int_{\bar\M_g}c_1(\bar E)\wedge\omega_{_{KE}}^{n-1}
\]
we need to show that
\begin{eqnarray}\label{c9}
\int_{\bar\M_g}c_1(\bar E)\wedge\omega_{_{KE}}^{n-1}=\int_{\M_g}\omega_{_{KE}}^n.
\end{eqnarray}
By the property of the asymptotic Poincar\'e metric, we know that
\begin{eqnarray}\label{100}
 \int_{\M_g} \omega_p^n <\infty.
\end{eqnarray}
Since the K\"ahler-Einstein metric is equivalent to the asymptotic
Poincar\'e metric, we know that
$\int_{\M_g}\omega_{_{KE}}^n<\infty$. Also, since $c_1(\bar E)$ is
a closed $(1,1)$-form on $\bar\M_g$ which is compact, we know that
any representative of $c_1(\bar E)$ is bounded on $\bar\M_g$. This
implies there is a
constant $C_5$ such that $-C_5\omega_p\leq c_1(\bar E)\leq
C_5\omega_p$. This implies that $\left |\int_{\bar\M_g}c_1(\bar
E)\wedge\omega_{_{KE}}^{n-1}\right |<\infty$. By using the
notations as in the above lemma we have
\begin{align}\label{c10}
\begin{split}
&\int_{\bar\M_g}c_1(\bar E)\wedge\omega_{_{KE}}^{n-1}-\int_{\M_g}\omega_{_{KE}}^n
=\int_{\M_g}(c_1(\bar E)-\omega_{_{KE}})\wedge\omega_{_{KE}}^{n-1}\\
=&\lim_{\epsilon\to 0}\int_{\M_g}\rho_\epsilon
(c_1(\bar E)-\omega_{_{KE}})\wedge\omega_{_{KE}}^{n-1}
=\lim_{\epsilon\to 0}\int_{\M_g}\rho_\epsilon  \partial\bar\partial\log\left (
\frac{\det(g_{_{KE}}^\ast)}{\det(h)}\right )\wedge\omega_{_{KE}}^{n-1}\\
=&\lim_{\epsilon\to 0}\int_{\M_g}\log\left (
\frac{\det(g_{_{KE}}^\ast)}{\det(h)}\right )\partial\bar\partial\rho_\epsilon
\wedge\omega_{_{KE}}^{n-1}\\
=&\lim_{\epsilon\to 0}\int_{\M_g\cap supp(\partial\bar\partial\rho_\epsilon)}\log\left (
\frac{\det(g_{_{KE}}^\ast)}{\det(h)}\right )\partial\bar\partial\rho_\epsilon
\wedge\omega_{_{KE}}^{n-1}.
\end{split}
\end{align}
Now we show that
\begin{eqnarray}\label{200}
 \int_{\M_g}  \left |\log\frac{\det(g_{_{KE}}^\ast)}{\det(h)}\right | \omega_p^n
 <\infty.
\end{eqnarray}
Since
\begin{align*}
\begin{split}
\int_{\M_g}  \left |\log\frac{\det(g_{_{KE}}^\ast)}{\det(h)}\right | \omega_p^n
=\sum_{i=1}^p
\int_{U_i}\psi_i \left |\log\frac{\det(g_{_{KE}}^\ast)}{\det(h)}\right | \omega_p^n
+\sum_{j=1}^q
\int_{V_j}\psi_{p+j}
\left |\log\frac{(g_{_{KE}}^\ast)}{\det(h)}\right | \omega_p^n
\end{split}
\end{align*}
and $V_j$ lies in the compact set $\M_g\setminus V_0$ and $0\leq \phi_i\leq 1$,
we only need to show that
\begin{eqnarray}\label{300}
\int_{U_i}\left |\log\frac{(g_{_{KE}}^\ast)}{\det(h)}\right | \omega_p^n
<\infty.
\end{eqnarray}
We know that there is a constant $C_i$ such that, on $U_i$,
\[
\omega_p^n\leq C_i \prod_{j=1}^{\alpha_i}\frac{1}{|t_j^i|^2(\log |t_j^i|)^2}.
\]
Formula \eqref{300} follows from the above formula, inequality \eqref{c5} and a simple
computation.

Now we pick $\epsilon$ small such that $supp(\partial\bar\partial\rho_\epsilon)
\subset supp(1-\rho_\epsilon)\subset V_0$. By Lemma \ref{cutoff} we know that there is a
constant $C$ such that
\[
-C\omega_p\leq \partial\bar\partial\rho_\epsilon\leq C\omega_p
\]
and
\[
0\leq \omega_{_{KE}}\leq C\omega_p.
\]
So we have
\begin{align*}
\begin{split}
&\left | \int_{\M_g\cap supp(\partial\bar\partial\rho_\epsilon)}\log\left (
\frac{\det(g_{_{KE}}^\ast)}{\det(h)}\right )\partial\bar\partial\rho_\epsilon
\wedge\omega_{_{KE}}^{n-1}\right |\\
\leq & C^n\sum_{i=1}^p \int_{U_i\cap supp(1-\rho_\epsilon)}
\left |\log\left (
\frac{\det(g_{_{KE}}^\ast)}{\det(h)}\right )\right |\omega_{p}^{n}\to 0
\end{split}
\end{align*}
as $\epsilon$ goes to $0$ because of inequality \eqref{300} and the fact that the
Lebesgue measure of
\[
U_i\cap supp(1-\rho_\epsilon)
\]
goes to $0$. By combining with formula \eqref{c10} we have proved
this lemma.

\qed

Now we define the pointwise version of the degree. Let $F\subset \bar E$ be a
holomorphic subbundle of rank $k\leq n$. Let $g_{_{KE}}^\ast\mid_F$ and $h\mid_F$ be
the restriction to $F$ of the metrics induced by the K\"ahler-Einstein metric and
the metric $h$. Let
\begin{eqnarray}\label{c11}
d(F)=-\partial\bar\partial\log\det(g_{_{KE}}^\ast\mid_F)\wedge\omega_{_{KE}}^{n-1}.
\end{eqnarray}
The following result is well-known. Please see \cite{ko1} for
details.
\begin{lemma}\label{slope}
For any holomorphic subbundle $F$ of $\bar E$ with rank $k$, we have
\begin{eqnarray}\label{slopeineq}
\frac{d(F)}{k}\leq \frac{d(\bar E)}{n}.
\end{eqnarray}
\end{lemma}

Now we prove the main theorem.

{\bf Proof.} Let $F$ be a holomorphic subbundle of $\bar E$ of
rank $k$. We first check that $\int_{\M_g}d(F)$ is finite and
equal to the degree of $F$. To prove that $\int_{\M_g}d(F)$ is
finite, we need to show that
$-\partial\bar\partial\log\det(g_{_{KE}}^\ast\mid_F)$ has
Poincar\'e growth. This involves the estimate of the derivatives
of the K\"ahler-Einstein metric up to second order. Our method is
to use Lemma \ref{slope} together with the monotone convergence
theorem and integration by parts to reduce the $C^2$ estimates of
the K\"ahler-Einstein metric to $C^0$ estimates.

By Lemma \ref{cutoff} we know that $\rho_\epsilon$ is monotonically increasing when
$\epsilon$ is monotonically decreasing. Also by Lemma \ref{slope} we know that
\[
\frac{k}{n}\omega_{_{KE}}^n-d(F)=\frac{k}{n}d(\bar E)-d(F)\geq 0.
\]
By the monotone convergence theorem we have
\begin{eqnarray}\label{c12}
\lim_{\epsilon\searrow 0}\int_{\M_g}\rho_\epsilon \left (\frac{k}{n}\omega_{_{KE}}^n-d(F)
\right )
=\int_{\M_g}\left (\frac{k}{n}\omega_{_{KE}}^n-d(F)\right ).
\end{eqnarray}
By Lemma \ref{amdegree} we have
\begin{align}\label{c13}
\begin{split}
&\frac{k}{n}\deg(\bar E)-\deg(F)=\int_{\M_g}\left ( \frac{k}{n}\omega_{_{KE}}^n-
Ric(h\mid_F)\wedge\omega_{_{KE}}^{n-1}\right )\\
=& \lim_{\epsilon\to 0}\int_{\M_g}\rho_\epsilon \left (
\frac{k}{n}\omega_{_{KE}}^n-
Ric(h\mid_F)\wedge\omega_{_{KE}}^{n-1}\right ).
\end{split}
\end{align}
However,
\begin{align}\label{c14}
\begin{split}
&\int_{\M_g}\rho_\epsilon \left (\frac{k}{n}\omega_{_{KE}}^n-d(F)\right )
-\int_{\M_g}\rho_\epsilon \left ( \frac{k}{n}\omega_{_{KE}}^n-
Ric(h\mid_F)\wedge\omega_{_{KE}}^{n-1}\right )\\
=& \int_{\M_g}\rho_\epsilon \left (Ric(h\mid_F)\wedge\omega_{_{KE}}^{n-1}
- d(F)\right )
= \int_{\M_g}\rho_\epsilon \left (Ric(h\mid_F)
+\partial\bar\partial\log\det\left (g_{_{KE}}^\ast\mid_F\right )\right )
\wedge\omega_{_{KE}}^{n-1}\\
=& \int_{\M_g}\rho_\epsilon\partial\bar\partial\log
 \left (\frac{\det\left (g_{_{KE}}^\ast\mid_F\right )}{\det\left (h\mid_F\right)}
\right )\wedge\omega_{_{KE}}^{n-1}
= \int_{\M_g}\log
 \left (\frac{\det\left (g_{_{KE}}^\ast\mid_F\right )}{\det\left (h\mid_F\right)}
\right )\partial\bar\partial\rho_\epsilon\wedge\omega_{_{KE}}^{n-1}\\
=&\int_{\M_g\cap supp(1-\rho_\epsilon)}\log
 \left (\frac{\det\left (g_{_{KE}}^\ast\mid_F\right )}{\det\left (h\mid_F\right)}
\right )\partial\bar\partial\rho_\epsilon\wedge\omega_{_{KE}}^{n-1}.
\end{split}
\end{align}

Now we show that
\begin{eqnarray}\label{c15}
\lim_{\epsilon\to 0}\int_{\M_g\cap supp(1-\rho_\epsilon)}\log
 \left (\frac{\det\left (g_{_{KE}}^\ast\mid_F\right )}{\det\left (h\mid_F\right)}
\right )\partial\bar\partial\rho_\epsilon\wedge\omega_{_{KE}}^{n-1}=0.
\end{eqnarray}

By the proof of Lemma \ref{amdegree}, to prove formula \eqref{c15}
we only need to show that
\begin{eqnarray}\label{c16}
 \int_{\M_g}\left |\log
 \left (\frac{\det\left (g_{_{KE}}^\ast\mid_F\right )}{\det\left (h\mid_F\right)}
\right )\right |\omega_p^n<\infty
\end{eqnarray}
which is reduced to show that
\begin{eqnarray}\label{c17}
 \int_{U_i}\left |\log
 \left (\frac{\det\left (g_{_{KE}}^\ast\mid_F\right )}{\det\left (h\mid_F\right)}
\right )\right |\omega_p^n<\infty.
\end{eqnarray}

Since the Ricci metric is equivalent to the K\"ahler-Einstein
metric, we know that $\frac{\det\left (g_{_{KE}}^\ast\mid_F\right
)}{\det\left (\tau^\ast\mid_F\right)}$ is bounded from above and
below by positive constants.

We fix a $U_i$. Let $\{ f_1,\cdots,f_k\}$ be a local holomorphic
frame of $F$. We know that there exists a $k\times n$ matrix
$B=(b_{\alpha\beta})$ whose entries are holomorphic functions on
$U_i$ such that the rank of $B$ is $k$ and $f_\alpha=\sum_\beta
b_{\alpha\beta}e_\beta^i$ where $e_\beta^i$ is defined in
\eqref{logframe}. Now we have
\[
\frac{\det\left (\tau^\ast\mid_F\right )}{\det\left
(h\mid_F\right)}=\frac{\det\left (B\left (\tau^\ast_{i\bar
j}\right )\bar B^T\right )}{\det\left (B\left (h_{i\bar
j}\right )\bar B^T\right )}.
\]

Now we need the following linear algebraic lemma:
\begin{lemma}\label{easylemma}
For any positive Hermitian $n\times n$ matrix $A$, we denote its
eigenvalues by $\lambda_1,\cdots,\lambda_n$ where
$\lambda_i=\lambda_i(A)$ such that
$\lambda_1(A)\geq\cdots\geq\lambda_n(A)>0$. Let $A_1$ and $A_2$ be
two positive Hermitian $n\times n$ matrices. Let $B$ be an
$k\times n$ matrix with $k\leq n$ such that the rank of $B$ is
$k$. Then there are positive constants $c_1$ and $c_2$ only
depending on $n,k$ such that
\[
c_1\frac{\lambda_{n-k+1}(A_1)\cdots\lambda_n(A_1)}
{\lambda_1(A_2)\cdots\lambda_k(A_2)}\leq
\frac{\det\left (BA_1\bar B^T\right )}{\det\left (BA_2\bar B^T\right )}
\leq c_2\frac{\lambda_1(A_1)\cdots\lambda_k(A_1)}
{\lambda_{n-k+1}(A_2)\cdots\lambda_n(A_3)}.
\]
\end{lemma}
We briefly show the proof here.

{\bf Proof.} We fix $A_1$ and $A_2$ and let
\[
Q(B)=\frac{\det\left (BA_1\bar B^T\right )}{\det\left (BA_2\bar B^T\right )}.
\]
Let $B_i$ be the $k\times n$ matrix obtained by multiplying the
$i$-th row of $B$ by a non-zero constant $c$ and leave other rows
invariant and let $B_{ij}$ be the $k\times n$ matrix obtained by
adding a constant multiple of the $j$-th row to the $i$-th row and
leave other rows invariant. It is easy to check that
$Q(B_i)=Q(B_{ij})=Q(B)$.

Thus we can assume that the row vectors of $B$ form an orthonormal set
of $\mathbb C^n$. With this assumption, it is easy to see that there
are positive constants $c_3$ and $c_4$ only depending on $n,k$ such
that
\[
c_3\lambda_{n-k+1}(A_i)\cdots\lambda_n(A_i)\leq
\det\left (BA_i\bar B^T\right )\leq c_4\lambda_1(A_i)\cdots\lambda_k(A_i)
\]
for $i=1,2$. The lemma follows directly.

\qed

Now we go back to the proof of the theorem. By using Theorem 1.4
and corollary 4.2 of \cite{lsy1}, we know that, under the frame
\eqref{logframe},
\begin{enumerate}
\item $\tau^\ast_{i\bar i}=u_i^{-2}(1+O(u_0))$ if $i\leq m$;
\item $\tau^\ast_{i\bar j}=O(1)$ if $i,j\leq m$ with $i\ne j$ or
$i\leq m<j$ or $j\leq m<i$;
\item $\tau^\ast_{i\bar j}=\tau^{i\bar j}$ if $i,j\geq m+1$;
\item On $U_i$,
the submatrix $\left ( \tau^{i\bar j}\right )_{i,j\geq m+1}$ is
bounded from above and below by positive constants multiple of the
identity matrix where the constants depend on $U_i$.
\end{enumerate}

It is clear that, on $U_i$ the eigenvalues of matric matrix of $h$
with respect to the frame \eqref{logframe} are bounded from above
and below by positive constants which depend on $U_i$ and the
choice of the metric $h$.

By analyzing the eigenvalues of the matrix $\left (\tau^{\ast}\right
)$ and by using Lemma \ref{easylemma}, a simple computation shows
that there are positive constants $C_6$ and $C_7$ which depend on
$F$ and $U_i$ such that
\begin{eqnarray}\label{c18}
\left |\log
 \left (\frac{\det\left (g_{_{KE}}^\ast\mid_F\right )}{\det\left (h\mid_F\right)}
\right )\right |\leq C_6+C_7\sum_{j=1}^{\alpha_i}\log\log\frac{1}{\left |t_j^i\right |}.
\end{eqnarray}
Now by using a similar method to the proof of Lemma \ref{amdegree}
we know that formula \eqref{c17} and \eqref{c16} hold which imply
formula \eqref{c15} holds. Combining \eqref{c15} and \eqref{c14}
we have
\begin{eqnarray}\label{c19}
\lim_{\epsilon\searrow 0}\left (
\int_{\M_g}\rho_\epsilon \left (\frac{k}{n}\omega_{_{KE}}^n-d(F)\right )
-\int_{\M_g}\rho_\epsilon \left ( \frac{k}{n}\omega_{_{KE}}^n-
Ric(h\mid_F)\wedge\omega_{_{KE}}^{n-1}\right )\right )=0.
\end{eqnarray}
By combining \eqref{c19}, \eqref{c13} and \eqref{c12} we have
\begin{eqnarray}\label{c20}
\int_{\M_g}\left (\frac{k}{n}\omega_{_{KE}}^n-d(F)\right )
=\frac{k}{n}\deg(\bar E)-\deg(F).
\end{eqnarray}
By Lemma \ref{amdegree} we know that $\int_{\M_g}\omega_{_{KE}}^n=\deg(\bar E)$. From
formula \eqref{c20} we know that $\int_{\M_g}d(F)$ is finite and
\begin{eqnarray}\label{c21}
\deg(F)=\int_{\M_g}d(F).
\end{eqnarray}
Now by Lemma \ref{slope} we have
\[
\frac{\deg(F)}{k}-\frac{\deg(\bar E)}{n}=\int_{\M_g}\left ( \frac{d(F)}{k}
-\frac{d(\bar E)}{n}\right )\leq 0
\]
which implies
\[
\frac{\deg(F)}{k}\leq \frac{\deg(\bar E)}{n}.
\]
This proves that the bundle $\bar E$ is semi-stable in the sense
of Mumford.

To prove the strict stability of the logarithm cotangent bundle, we
need to show that this bundle cannot split. The following result about
the moduli group $\text{Mod}_g$ and its proof is due to F. Luo \cite{luo1}.
\begin{prop}\label{nonsplit}
Let $\text{Mod}_g$ be the moduli group of closed Riemann surfaces of genus
$g$ with $g\geq 2$. Then any finite index subgroup of $\text{Mod}_g$ is not
isomorphic to a product of groups.
\end{prop}

The proof of the above proposition is topological. For
completeness, we will include Luo's proof at the end of this
section.

Now we go back to the proof of stability.
If $\bar E$ is not stable, then it must split into a direct sum of
holomorphic subbundles $\bar E=\oplus_{i=1}^k E_i$ with $k\geq 2$.
Moreover, when restricted to the moduli space, both the connection
and the K\"ahler-Einstein metric split. It is well known that there
is a finite cover $\widetilde M_g$ of $\M_g$ which is smooth.
By the decomposition theorem of de Rham, the Teichm\"uller space,
as the universal covering space of $\widetilde M_g$ must split
isomorphically as a product of manifolds. Furthermore, the fundamental
group of $\widetilde M_g$ is isomorphic
to a product of groups. However, $\pi_1(\widetilde M_g)$ is a finite
index subgroup of the mapping class group. By the above proposition,
this is impossible. So we have proved the stability.

\qed

We remark that the positivity of $c_1({\bar E})$ implies that the
Deligne-Mumford compactification ${\overline \M}_g$ is of
logarithmic general type for $g>1$.

In the end of this section we give a proof by F. Luo of
Proposition \ref{nonsplit}.

{\bf Proof of Proposition \ref{nonsplit}.} The proof uses
Thurston's classification of elements in the mapping class group
\cite{th} and the solution of the Nielsen realization problem
\cite{ker1}. In the following we will use the words ``simple
loops'' and ``subsurfaces'' to denote the isotopy classes of
simple loops and subsurfaces. We fix a surface $X$.

Suppose there is a subgroup $G$ of $\text{Mod}_g$ with finite
index such that there are two nontrivial subgroups $A$ and $B$ of
$G$ so that $G=A\times B$. We will derive a contradiction. We need
to following lemma.
\begin{lemma}\label{applem1}
Let $A$, $B$ and $G$ be as above. Then
\begin{enumerate}
\item There are elements of infinite order in both $A$ and $B$;
\item There are no elements in $A$ or $B$ which is pseudo-Anosov.
\end{enumerate}
\end{lemma}

{\bf Proof.} If the first claim is not true, then we can assume
$A$ consists of torsions only. We claim that, in this case, $A$ is
a finite group.

Actually let $\pi:\text{Mod}_g\to\text{Aut}(H_1(X))$ be the
natural homomorphism. By the virtual of Theorem V.3.1 of
\cite{farkra1}, we know that the kernel of $\pi$ contains no
torsion elements. Thus $\pi(A)$ is isomorphic to $A$. Now $\pi(A)$
is a torsion subgroup of the general linear group
$\text{GL}(n,\mathbf R)$. By the well known solution of the
Burnside problem for the linear groups, we see that $\pi(A)$ must
be finite and so is $A$.

For the finite group $A$, by the solution of the Nielsen
realization problem of Kerckhoff \cite{ker1}, we know that there
is a point $d$ in the Teichm\"uller space of $X$ fixed by all
elements in $A$. Let $\text{Fix}(A)$ be the set of all points in
the Teichm\"uller space fixed by each element in $A$. Then we see
that $\text{Fix}(A)$ is a non-empty proper subset of the
Teichm\"uller space. Now for each $a\in A$, $b\in B$ and
$d\in\text{Fix}(A)$, since $ab=ba$, we have $ab(d)=ba(d)=b(d)$.
This implies $\text{Fix}(A)$ is invariant under the action of $B$
on the Teichm\"uller space. Thus we see that the finite index
subgroup $G=A\times B$ acts on the Teichm\"uller space leaving
$\text{Fix}(A)$ invariant. This contradicts the finite index
property of $G$ since $\text{Fix}(A)$ is the Teichm\"uller space
of the orbifold $X/A$. This proved the first claim.

Now we check the second claim. If it is not true, then we can
assume that there is an $a\in A$ which is pseudo-Anosov. Now we
consider the action of the mapping class group on the space of all
measured laminations in the surface. By Thurston's theory, there
are exactly two measured laminations $m,m'$ fixed by $a$. Now for
all $b\in B$, due to $ab=ba$, we see that $b$ leaves $\{ m,m'\}$
invariant. A result of McCarthy \cite{mcc} shows that the
stablizers of $\{ m,m'\}$ in the mapping class group is virtually
cyclic. Thus we see that each element $b\in B$ has some power
$b^n$ with $n\ne 0$ which is equal to $a^k$ with $k\ne 0$. By the
first claim we know that $B$ contains elements of infinite order.
This implies that some power $a^k$ with $k\ne 0$ is in $B$ which
is a contradiction. This proved the second claim.

\qed

Now we go back to the proof of the proposition. By the above
lemma, we can take $a\in A$ and $b\in B$, both are infinite order
and none of them is pseudo-Anosov. Thus by replacing $a$ and $b$
by a high power $a^n$ and $b^n$ with $n>0$, we may assume that for
$a$ there is a set of disjoint simple loops $c_1,\cdots,c_k$ in
$X$ so that
\begin{enumerate}
\item Each component of $X\setminus \cup_{i=1}^k c_i$ is invariant
under $a$;
\item The restriction of $a$ to each component of
$X\setminus \cup_{i=1}^k c_i$ is the identity map or a
pseudo-Anosov map.
\end{enumerate}
Let $N(c_i)$ be a regular neighborhood of $c_i$ in $X$ and let
$F(a)$ be the subsurface which is the union of all $N(c_i)$'s with
those components of $X\setminus \cup_{i=1}^k c_i$ on which $a$ is
the identity map. By Thurston's classification, the subsurface
$F(a)$ has the property that if $c$ is a curve system invariant
under $a$, then $c$ is in $F(a)$. We can construct $F(b)$ in a
similar way.

Now since $ab=ba$ we know that if a simple loop $c$ is invariant
under $b$, then $a(c)$ is also invariant under $b$ because
$ba(c)=ab(c)=a(c)$. Thus for all simple loops $c$ in $F(b)$,
$a(c)$ is still in $F(b)$. This implies that $a(F(b))=F(b)$. By
using the property of $F(a)$, we see $F(b)\subset F(a)$. Similar
reason implies $F(a)\subset F(b)$. Thus $F(a)=F(b)$.

Now we show that the subsurface $F(a)=F(b)$ is invariant under
each element in $A$ and $B$. We pick $x\in A$. Due to $xb=bx$, we
have $x(F(b))=b(x(F(b)))$ which implies $x(F(b))\subset F(b)$.
Since these two surfaces are homeomorphic, we have $x(F(b))=F(b)$.
Similarly, for all $y\in B$, $y(F(a))=F(a)$. Thus all elements in
$G=A\times B$ leave the subsurface $F(a)$ invariant.

Now there are two cases. In the first case, $F(a)$ is not
homeomorphic to $X$. In the second case, $F(a)=F(b)=X$.

It is clear that the first case $F(a)$ is not homeomorphic to $X$
cannot occur. Otherwise, the finite index subgroup $G$ leaves a
proper subsurface $F(a)$ invariant. This contradicts the known
properties of the mapping class group. As conclusion of this case,
we see that for any element $x\in A$, any power $x^n$ of $x$
cannot contain pseudo-Anosov components in Thurston's
classification.

In the second case, $F(a)=F(b)=X$. In this case, some power of $a$
(and $b$) is a composition of Dehn twists on disjoint simple
loops. By the argument in the first case, we see that for any two
elements $x\in A$ and $y\in B$, some powers $x^n$, $y^m$ with
$m,n\ne 0$ are either the identity map or compositions of Dehn
twists on disjoint simple loops.

let $\text{Fix}(x)$ be the union of the disjoint simple loops so
that $x$ is the composition of Dehn twists on these simple loops.
We define $\text{Fix}(x)$ to be the empty set if $x^n=id$ for some
non-zero integer $n$. We need the following claim:
\begin{claim}
For any $x\in A$ and $y\in B$, the geometric intersection number
between $\text{Fix}(x)$ and $\text{Fix}(y)$ is zero.
\end{claim}

To prove the claim, without loss of generality, we may take $x=a$
and $y=b$. We may assume that $a$ is the composition of Dehn
twists on curve system $\text{Fix}(a)=c$ and $b$ is the
composition of Dehn twists on curve system $\text{Fix}(b)=d$. Note
that in this case, if $s$ is a curve system invariant under $a$,
then the geometric intersection number between $s$ and $c$ is
zero. Namely $I(s,c)=0$. We also say that $s$ is \it disjoint \rm
from $c$. Now since $ab=ba$ and $a(c)=c$, we have that
$a(b(c))=b(c)$. Thus the curve system $b(c)$ is disjoint from $c$.
Since $b$ is a composition of Dehn twists on disjoint simple loops
$d$, this shows that $I(c,d)=0$. Namely $d$ and $c$ are disjoint
curve systems.

Now we finish the proof of the proposition as follows. First of
all, by the assumption $\text{Fix}(a)$ and $\text{Fix}(b)$ are
both non-empty. Let $S_1$ (and $S_2$) be the smallest subsurface
of $X$ which contains all curves in $\text{Fix}(x)$ for $x\in A$
(or $x\in B$). By the above claim, we have $S_1\cap
S_2=\emptyset$. Thus there is an essential simple loop $c$ which
is disjoint from both $S_1$ and $S_2$. By the construction,
$x(c)=c$ and $y(c)=c$ for all $x\in A$ and $y\in B$. Thus we see
that for each element $e\in G=A\times B$, there is a power $e^n$
with $n\ne 0$ so that $e^n$ leaves $c$ invariant. This contradicts
the fact that $A\times B$ is a finite index subgroup of the
mapping class group. This finishes the proof.

\qed

\section{The Bounded Geometry of the K\"ahler-Einstein Metric}\label{bddgeo}

In this section we show that the K\"ahler-Einstein metric has
bounded curvature and the injectivity radius of the Teichm\"uller
space equipped with the K\"ahler-Einstein metric is bounded from
below.

We begin with a metric that is equivalent to the K\"ahler-Einstein
metric and its curvature as well as the covariant derivatives of
the curvature are uniformly bounded. We deform the Ricci metric
whose curvature is bounded to obtain this metric by using
K\"ahler-Ricci flow. Then we establish the Monge-Amper\'e equation
from this new metric. Our work in \cite{lsy1} implies the new
metric is equivalent to the K\"ahler-Einstein metric which gave us
$C^2$ estimates. Based on this, we do the $C^3$ and $C^4$
estimates. This will give us the boundedness of the curvature of
the K\"ahler-Einstein metric. The same method can be used to show
that all the covariant derivatives of the K\"ahler-Einstein metric
are bounded.

We slightly change our notations. We will use
$g$ to denote the K\"ahler-Einstein metric and use $h$ to denote
an equivalent metric whose curvature and covariant derivatives of
curvature are bounded. We use $z_1,\cdots,z_n$ to denote local
holomorphic coordinates on the Teichm\"uller space. The main
result of this section is the following theorem:

\begin{theorem}\label{kebdd}
Let $g$ be the K\"ahler-Einstein metric on the Teichm\"uller space
$\mathcal T$. Then the curvature of $g$ and all of its covariant
derivatives are all bounded.
\end{theorem}
{\bf Proof.} We begin with the Ricci metric $\tau$. We first
obtain a new equivalent metric $h$ by deforming $\tau$ with the
Ricci flow. Consider the following K\"ahler-Ricci flow:
\begin{eqnarray}\label{bdd10}
\begin{cases}
\frac{\partial g_{i\bar j}}{\partial t}=-(R_{i\bar j}+g_{i\bar j})\\
g_{i\bar j}(0)=\tau_{i\bar j}
\end{cases}.
\end{eqnarray}
If we let $s=e^t-1$ and $\widetilde g=e^t g$, we have
\begin{eqnarray}\label{bdd20}
\begin{cases}
\frac{\partial\widetilde g_{i\bar j}}{\partial s}=-\widetilde R_{i\bar
j}\\
\widetilde g_{i\bar j}(0)=\tau_{i\bar j}
\end{cases}.
\end{eqnarray}
Since the initial metric has bounded curvature, by the work of Shi
\cite{shi1}, the flow \eqref{bdd20} has short time existence and
for small $s$, the metric $\widetilde g(s)$ is equivalent to the
initial metric $\tau$. Furthermore, the curvature and its
covariant derivatives of $\widetilde g(s)$ are bounded. Hence for
small $t$, the metric $g(t)$ is equivalent to the Ricci metric
$\tau$ and has bounded curvature as well as covariant derivatives
of the curvature.

Now we fix a small $t$ and denote the metric $g(t)$ by $h$. Since the
Teichm\"uller space is contractible, there are smooth functions $u$
and $F$ such that
\begin{eqnarray}\label{gh}
\omega_g=\omega_h+\partial\bar\partial u
\end{eqnarray}
and
\begin{eqnarray}\label{rich}
Ric(h)+\omega_h=\partial\bar\partial F.
\end{eqnarray}
Since the metrics $h$ and $\tau$ are equivalent, we know that $h$
and $g$ are equivalent which implies the tensor $u_{i\bar
j}dz_i\otimes d\bar z_j$ is bounded with respect to either metric.
Also, because the curvature and its covariant derivatives of the
metric $h$ are bounded, we know that a covariant derivative of $F$
with respect to $h$ is bounded if this derivative is at least
order $2$ and has at least one holomorphic direction and one
anti-holomorphic direction. So we have $C^2$ estimates.

By the K\"ahler-Einstein condition of $g$, we have the Monge-Amp\`ere
equation
\begin{eqnarray}\label{maeq}
\log\det(h_{i\bar j}+u_{i\bar j})-\log\det(h_{i\bar j})=u+F.
\end{eqnarray}

We use $\Delta$, $\Delta^{'}$, $\nabla$, $\nabla^{'}$,
$\Gamma_{ij}^k$, $\widetilde\Gamma_{ij}^k$, $R_{i\bar jk\bar l}$,
$P_{i\bar jk\bar l}$, $R_{i\bar j}$, $P_{i\bar j}$, $R$ and $P$ to denote the
Laplacian, gradient, Christoffell symbol, curvature tensor, Ricci
curvature and scalar curvature of the metrics $h$ and $g$ respectively.
In the following, all covariant derivatives of functions and tensors
are taken with respect to the background metric $h$.

Inspired by Yau's work in \cite{yau3}, we let
\begin{eqnarray}\label{bdd25}
\mathcal F=u+F,
\end{eqnarray}
\begin{eqnarray}\label{bdd30}
S=g^{i\bar j}g^{k\bar l}g^{p\bar q}u_{;i\bar qk}u_{;\bar jp\bar l}
\end{eqnarray}
and
\begin{eqnarray}\label{bdd40}
V=g^{i\bar j}g^{k\bar l}g^{p\bar q}g^{m\bar n}
u_{;i\bar qk\bar n}u_{;\bar jp\bar lm}+
g^{i\bar j}g^{k\bar l}g^{p\bar q}g^{m\bar n}
u_{;i\bar nkp}u_{;\bar jm\bar l\bar q}.
\end{eqnarray}
To simplify the notation, we define the following quantities:
\begin{align}\label{bdd4100}
\begin{split}
A_{ipkm\alpha}=& u_{;i\bar pk\bar m\alpha}-g^{\gamma\bar\delta}
u_{;i\bar\delta\alpha}u_{;\gamma\bar pk\bar m}-g^{\gamma\bar\delta}
u_{;k\bar\delta\alpha}u_{;i\bar p\gamma\bar m}-g^{\gamma\bar\delta}
u_{;\gamma\bar p\alpha}u_{;i\bar\delta k\bar m}\\
&-g^{\gamma\bar\delta}u_{;\gamma\bar m\alpha}u_{;i\bar pk\bar\delta}
+g^{\gamma\bar\delta}u_{;i\bar p\gamma}u_{;\bar\delta k\bar m\alpha}
+g^{\gamma\bar\delta}u_{;i\bar p\gamma}u_{;k\bar m\alpha\bar\delta}
+2g^{\gamma\bar\delta}u_{;\bar pi\bar\delta}u_{;k\bar m\alpha\gamma},
\end{split}
\end{align}
\begin{align}\label{bdd4200}
\begin{split}
B_{ipkm\alpha}=& u_{;\bar ip\bar km\alpha}-g^{\gamma\bar\delta}
u_{;\gamma\bar i\alpha}u_{;\bar\delta p\bar km}-g^{\gamma\bar\delta}
u_{;\gamma\bar k\alpha}u_{;\bar ip\bar\delta m}-g^{\gamma\bar\delta}
u_{;p\bar\delta\alpha}u_{;\bar i\gamma\bar km}\\
&-g^{\gamma\bar\delta}u_{;m\bar\delta\alpha}u_{;\bar ip\bar k\gamma}
+g^{\gamma\bar\delta}u_{;p\bar i\gamma}u_{;\bar km\bar\delta\alpha}
+g^{\gamma\bar\delta}u_{;\bar ip\bar\delta}u_{;\gamma\bar km\alpha},
\end{split}
\end{align}
\begin{align}\label{bdd4300}
\begin{split}
C_{ipkm\alpha}=& u_{;i\bar mkp\alpha}-g^{\gamma\bar\delta}
u_{;i\bar\delta\alpha}u_{;\gamma\bar mkp}-g^{\gamma\bar\delta}
u_{;k\bar\delta\alpha}u_{;i\bar m\gamma p}-g^{\gamma\bar\delta}
u_{;p\bar\delta\alpha}u_{;i\bar mk\gamma}\\
&-g^{\gamma\bar\delta}u_{;\gamma\bar m\alpha}u_{;i\bar\delta kp}
+g^{\gamma\bar\delta}u_{;i\bar m\gamma}u_{;k\bar\delta p\alpha},
\end{split}
\end{align}
\begin{align}\label{bdd4400}
\begin{split}
D_{ipkm\alpha}=& u_{;\bar im\bar k\bar p\alpha}-g^{\gamma\bar\delta}
u_{;\gamma\bar i\alpha}u_{;\bar\delta m\bar k\bar p}-g^{\gamma\bar\delta}
u_{;\gamma\bar k\alpha}u_{;\bar im\bar\delta\bar p}-g^{\gamma\bar\delta}
u_{;\gamma\bar p\alpha}u_{;\bar im\bar k\bar\delta}\\
&-g^{\gamma\bar\delta}u_{;m\bar\delta\alpha}u_{;\bar i\gamma\bar k\bar p}
+g^{\gamma\bar\delta}u_{;\bar im\bar\delta}u_{;\bar k\gamma\bar p\alpha},
\end{split}
\end{align}
\begin{align}\label{bdd50}
\begin{split}
W=& g^{i\bar j}g^{k\bar l}g^{p\bar q}g^{m\bar n}g^{\alpha\bar\beta}
u_{;i\bar qk\bar n\alpha}u_{;\bar jp\bar lm\bar\beta}
+g^{i\bar j}g^{k\bar l}g^{p\bar q}g^{m\bar n}g^{\alpha\bar\beta}
u_{;i\bar qk\bar n\bar\beta}u_{;\bar jp\bar lm\alpha}\\
&+g^{i\bar j}g^{k\bar l}g^{p\bar q}g^{m\bar n}g^{\alpha\bar\beta}
u_{;i\bar nkp\alpha}u_{;\bar jm\bar l\bar q\bar\beta}
+g^{i\bar j}g^{k\bar l}g^{p\bar q}g^{m\bar n}g^{\alpha\bar\beta}
u_{;i\bar nkp\bar\beta}u_{;\bar jm\bar l\bar q\alpha}
\end{split}
\end{align}
and
\begin{align}\label{bdd55}
\begin{split}
\widetilde W=& g^{i\bar j}g^{k\bar l}g^{p\bar q}g^{m\bar n}g^{\alpha\bar\beta}
\left (
A_{iqkn\alpha}\bar{A_{jplm\beta}}
+B_{jplm\alpha}\bar{B_{iqkn\beta}}\right )\\
&+g^{i\bar j}g^{k\bar l}g^{p\bar q}g^{m\bar n}g^{\alpha\bar\beta}\left (
C_{ipkn\alpha}\bar{C_{jqlm\beta}}
+B_{jqlm\alpha}\bar{B_{ipkn\beta}}
\right ).
\end{split}
\end{align}

Firstly, a simple computation shows that
\begin{eqnarray}\label{crist}
\widetilde\Gamma_{ik}^\alpha-\Gamma_{ik}^\alpha=g^{\alpha\bar\beta}
u_{;i\bar\beta k}.
\end{eqnarray}
Now we compute $P_{i\bar jk\bar l}$. We first note that
\[
u_{;i\bar jk\bar l}=u_{i\bar jk\bar l}-u_{\bar jp\bar l}\Gamma_{ik}^p
-u_{i\bar qk}\bar{\Gamma_{jl}^q}+u_{p\bar
  q}\Gamma_{ik}^p\bar{\Gamma_{jl}^q}
-u_{p\bar j}h^{p\bar q}R_{i\bar qk\bar l}.
\]
Since $g_{i\bar j}=h_{i\bar j}+u_{i\bar j}$ we have
\begin{align}\label{bdd60}
\begin{split}
P_{i\bar jk\bar l}=& \partial_k\partial_{\bar l}g_{i\bar j}-
g^{p\bar q}\partial_k g_{i\bar q}\partial_{\bar l}g_{p\bar j}
= \partial_k\partial_{\bar l}h_{i\bar j}+u_{i\bar jk\bar l}
-g_{p\bar q}\widetilde\Gamma_{ik}^p\bar{\widetilde\Gamma_{jl}^q} \\
=& \partial_k\partial_{\bar l}h_{i\bar j}+u_{i\bar jk\bar l}
-g_{p\bar q}\left (\Gamma_{ik}^p+g^{p\bar\beta}u_{;i\bar\beta k}
\right )\left ( \bar{\Gamma_{jl}^q}+g^{\alpha\bar q}u_{;\bar
  j\alpha\bar l}\right )\\
=& \left ( \partial_k\partial_{\bar l}h_{i\bar j}-h_{p\bar q}
\Gamma_{ik}^p\bar{\Gamma_{jl}^q}\right )-u_{p\bar q}
\Gamma_{ik}^p\bar{\Gamma_{jl}^q}\\
&+\left (u_{i\bar jk\bar l}-u_{\bar jp\bar l}\Gamma_{ik}^p
-u_{i\bar qk}\bar{\Gamma_{jl}^q}+2u_{p\bar
  q}\Gamma_{ik}^p\bar{\Gamma_{jl}^q} \right )-g^{p\bar q}u_{;i\bar qk}
u_{;\bar jp\bar l}\\
=& R_{i\bar jk\bar l}+u_{p\bar j}h^{p\bar q}R_{i\bar qk\bar l}
+u_{;i\bar jk\bar l}-g^{p\bar q}u_{;i\bar qk}
u_{;\bar jp\bar l}.
\end{split}
\end{align}
Since the curvature of the background metric $h$ and the tensor
$u_{i\bar j}dz_i\otimes d\bar z_j$ are both bounded,
to prove that the curvature of the metric $g$ is bounded, we only need
to show that both $S$ and $V$ are bounded.

We first consider the quantity $S$. We follow the idea of Yau in
\cite{yau3} and use the following notations:
\begin{definition}
Let $A$ and $B$ be two functions. We denote
\begin{enumerate}
\item $A \overset{3}{\simeq} B$ if $|A-B|\leq C_1\sqrt{S}+C_2$;
\item $A \overset{4}{\simeq} B$ if $|A-B|\leq C_1\sqrt{V}+C_2$;
\item $A \overset{3}{\cong} B$ if $|A-B|\leq C_1S+C_2\sqrt{S}+C_3$;
\item $A \overset{4}{\cong} B$ if $|A-B|\leq C_1V+C_2\sqrt{V}+C_3$
\end{enumerate}
where $C_1$, $C_2$ and $C_3$ are universal constants.

Also, by diagonalizing we mean to choose holomorphic coordinates
$z_1,\cdots,z_n$ such that
\[
h_{i\bar j}=\delta_{ij}
\]
and
\[
u_{i\bar j}=\delta_{ij}u_{i\bar i}.
\]
\end{definition}

Now we differentiate the equation \eqref{maeq} twice and
reorganize the terms. We have
\begin{eqnarray}\label{bdd70}
g^{i\bar j}u_{;i\bar jk\bar l}=\mathcal F_{k\bar l}+g^{i\bar j}
g^{p\bar q}u_{;i\bar qk}u_{;\bar jp\bar l}.
\end{eqnarray}
By differentiating this equation once more we have
\begin{align}\label{bdd80}
\begin{split}
g^{i\bar j}u_{;i\bar jk\bar l\alpha}=&\mathcal F_{;k\bar l\alpha}
+g^{i\bar j}g^{p\bar q}\left (u_{;i\bar q\alpha}u_{;p\bar jk\bar l}
+u_{;i\bar q k}u_{;\bar jp\bar l\alpha}+
u_{;i\bar q k\alpha}u_{;\bar jp\bar l}\right )\\
&-g^{i\bar j}g^{p\bar q}g^{m\bar n}\left ( u_{;\bar jm\bar l}
u_{;p\bar n\alpha}u_{;i\bar qk}+u_{;m\bar qk}u_{;i\bar n\alpha}
u_{;\bar jp\bar l}\right ).
\end{split}
\end{align}
Since
\[
\partial_k(\Delta u)=\partial_k\left (h^{i\bar j}u_{i\bar j}\right )
=h^{i\bar j}\left (u_{i\bar jk}-u_{p\bar j}\Gamma_{ik}^p\right )
=h^{i\bar j}u_{;i\bar jk},
\]
by diagonalizing and the Schwarz inequality we have
\begin{eqnarray}\label{bdd90}
\left |\nabla^{'}(\Delta u)\right |^2=g^{i\bar j}\left (h^{k\bar l}
u_{;k\bar li}\right )\left (h^{p\bar q}u_{;\bar qp\bar j}\right )=
\sum_{i}\frac{1}{1+u_{i\bar i}}\left |\sum_k u_{;k\bar ki}\right
|^2\leq C_1 S
\end{eqnarray}
since the metrics $h$ and $g$ are equivalent and
$S=\sum_{i,p,k}\frac{1}{1+u_{i\bar i}}\frac{1}{1+u_{p\bar p}}
\frac{1}{1+u_{k\bar k}}|u_{;i\bar pk}|^2$. We also have
\[
\Delta^{'}(\Delta u)= g^{k\bar l}\partial_k\partial_{\bar l}\left (
h^{i\bar j}u_{i\bar j}\right )=g^{k\bar l}h^{i\bar j}u_{;i\bar jk\bar l}
=g^{k\bar l}h^{i\bar j}\left ( u_{;k\bar li\bar j}-
u_{p\bar j}h^{p\bar q}R_{i\bar qk\bar l}+
u_{p\bar l}h^{p\bar q}R_{k\bar qi\bar j}\right ).
\]
By using equation \eqref{bdd70} we have
\begin{align}\label{bdd100}
\begin{split}
\Delta^{'}(\Delta u)=& h^{i\bar j} \left ( \mathcal F_{i\bar j}+
g^{k\bar l}g^{p\bar q}u_{;k\bar q i}u_{;\bar lp\bar j}\right )
-h^{i\bar j}g^{k\bar l}h^{p\bar q}u_{p\bar j}R_{i\bar qk\bar l}
+h^{i\bar j}g^{k\bar l}h^{p\bar q}u_{p\bar l}R_{k\bar qi\bar j}\\
=&\widetilde S+\Delta\mathcal F
-h^{i\bar j}g^{k\bar l}h^{p\bar q}u_{p\bar j}R_{i\bar qk\bar l}
+h^{i\bar j}g^{k\bar l}h^{p\bar q}u_{p\bar l}R_{k\bar qi\bar j}
\end{split}
\end{align}
where $\widetilde S= h^{i\bar j}g^{k\bar l}g^{p\bar q}
u_{;k\bar q i}u_{;\bar lp\bar j}$. Since the metrics $h$ and $g$ are
equivalent, we know that there is a constant $C_2$ such that
\[
\widetilde S\geq C_2 S.
\]
Now the term $\left |
h^{i\bar j}g^{k\bar l}h^{p\bar q}u_{p\bar j}R_{i\bar qk\bar l}\right
|$ is bounded since $h$ is equivalent to $g$, the curvature of $h$ is
bounded and we have $C^2$ estimates on $u$. Similarly, $\left |\Delta\mathcal
F\right |$ is bounded. Finally, since
\[
h^{i\bar j}g^{k\bar l}h^{p\bar q}u_{p\bar l}R_{k\bar qi\bar j}
=-g^{k\bar l}h^{p\bar q}u_{p\bar l}R_{k\bar q}
=-g^{k\bar l}h^{p\bar q}(g_{p\bar l}-h_{p\bar l})R_{k\bar q}=
g^{k\bar l}R_{k\bar q}-R,
\]
we know that $\left | h^{i\bar j}g^{k\bar l}h^{p\bar q}u_{p\bar
l}R_{k\bar qi\bar j}\right |$ is also bounded for similar reasons.
By combining the above argument we know that there is a constant
$C_3$ such that
\begin{eqnarray}\label{bdd110}
\Delta^{'}(\Delta u)\geq C_2 S-C_3.
\end{eqnarray}
Now by the computation in \cite{yau3}, we know there are positive
constants $C_4$, $C_5$ and $C_6$ such that
\begin{eqnarray}\label{bdd120}
\Delta^{'}(S+C_4\Delta u)\geq C_5S-C_6.
\end{eqnarray}
So for any positive $\lambda>0$, we can find a positive constant $C_7$
such that
\begin{eqnarray}\label{bdd130}
S+C_4\Delta u+C_7\geq 0
\end{eqnarray}
and
\begin{eqnarray}\label{bdd140}
\Delta^{'}(S+C_4\Delta u+C_7)\geq -\lambda(S+C_4\Delta u+C_7)
\end{eqnarray}
since $\Delta u$ is bounded. We fix $\lambda$ and let
$f=S+C_4\Delta u+C_7$. Now we know that the Ricci curvature of $g$
is $-1$ and $g$ is equivalent to the Ricci metric $\tau$ whose
injectivity radius has a lower bound, by the work of Li and Schoen
\cite{lischoen1} and Li \cite{lip1} we can find a positive $r_0$
such that the mean value inequality
\begin{eqnarray}\label{bdd150}
f(p)\leq C V_p^{-1}(r_0)\int_{B_p(r_0)} f \ dV
\end{eqnarray}
hold for any $p$ in the Teichm\"uller space. Here $V_p(r_0)$ is the
volume of the K\"ahler-Einstein ball centered at $p$ with radius
$r_0$, $dV=\omega_g^n$ is the volume element of the metric $g$
and $C$ is a constant depending on $r_0$ and $\lambda$ but is
independent of $p$. Let $r(z)$ be the function on $B_p(2r_0)$
measuring the $g$-distance between $z$ and $p$.
We fix a small $r_0$ and let $\rho=\rho(r)$ be a cutoff
function such that $0\leq \rho\leq 1$,
 $\rho(r)=1$ for $r\leq r_0$ and $\rho(r)=0$ for
$r\geq 2r_0$. Since $\Delta^{'}(\Delta u)+C_3\geq C_2S\geq 0$, we have
\begin{align*}
\begin{split}
C_2\int_{B_p(2r_0)}\rho^2 S\ dV-C_3V_p(2r_0)\leq &
\int_{B_p(2r_0)}\rho^2\Delta^{'}(\Delta u)\ dV=-2
\int_{B_p(2r_0)}\nabla^{'}\rho\cdot(\rho\nabla^{'}(\Delta u))\ dV\\
\leq & C_8\left (\int_{B_p(2r_0)}\left |\nabla^{'}\rho\right |^2\ dV
\right )^{\frac{1}{2}} \left (\int_{B_p(2r_0)}\rho^2\left |
\nabla^{'}(\Delta u)\right |^2\ dV\right )^{\frac{1}{2}}.
\end{split}
\end{align*}
Since $\left |\nabla^{'}\rho\right |$ is bounded and $\left |
\nabla^{'}(\Delta u)\right |^2\leq C_1S$, we have
\[
C_2\int_{B_p(2r_0)}\rho^2 S\ dV-C_3V_p(2r_0)\leq C_9
\left (\int_{B_p(2r_0)}\rho^2 S\ dV\right )^{\frac{1}{2}}
\left (V_p(2r_0)\right )^{\frac{1}{2}}
\]
which implies
\begin{eqnarray}\label{bdd160}
\int_{B_p(2r_0)}\rho^2 S\ dV\leq C_{10}V_p(2r_0).
\end{eqnarray}
By using inequalities \eqref{bdd150}and \eqref{bdd160} we have
\begin{align}\label{bdd170}
\begin{split}
f(p)\leq & C V_p^{-1}(r_0)\int_{B_p(r_0)}(S+C_4\Delta u+C_7)\ dV
\leq  C V_p^{-1}(r_0)\int_{B_p(r_0)}S\ dV+C_{11}\\
\leq & C V_p^{-1}(r_0)\int_{B_p(2r_0)}\rho^2 S\ dV+C_{11}
\leq C_{12}\frac{V_p(2r_0)}{V_p(r_0)}+C_{11}.
\end{split}
\end{align}

For each point $p\in\mathcal T_g$, let $f_p:\ \mathcal T_g\to
\mathbb C^{3g-3}$ be the Bers' embedding map such that $f_p(p)=0$
and $B_2\subset f_p(\mathcal T_g)\subset B_6$ where $B_r\subset
\mathbb C^{3g-3}$ is the open Euclidean ball of radius $r$. Since
both metrics $h$ and $g$ are equivalent to the Ricci metric which
is equivalent to the Euclidean metric on the unit Euclidean ball
$B_1$, we know that $\frac{V_p(2r_0)}{V_p(r_0)}$ is uniformly
bounded since both balls have Euclidean volume growth. Thus
$f(p)\leq C_{13}$. Since $\Delta u$ is bounded, we conclude that
$S$ is uniformly bounded.

Now we do the $C^4$ estimate. Let $\kappa$ be a large positive
constant. We first compute $\Delta^{'}\left [(S+\kappa)V\right ]$. We have
\begin{align}\label{bdd180}
\begin{split}
\Delta^{'}\left [(S+\kappa)V\right ]=& g^{i\bar j}
\partial_i\partial_{\bar j}\left [(S+\kappa)V\right ]=
g^{i\bar j}\partial_{\bar j}\left [V \partial_i S+
(S+\kappa)\partial_i V\right ]\\
=& g^{i\bar j}\left [(S+\kappa)
\partial_i\partial_{\bar j}V+V\partial_i\partial_{\bar j}S+
\partial_i V\partial_{\bar j}S+\partial_i S\partial_{\bar j}V
\right ]\\
\geq & (S+\kappa)\Delta^{'}V+V\Delta^{'}S-2\left |\nabla^{'}S\right |
\left |\nabla^{'}V\right |.
\end{split}
\end{align}
The computations of the first and second derivatives of $V$ are
very long. We only list the results here. For the first derivative
of $V$, we have
\begin{align}\label{bdd190}
\begin{split}
\partial_\alpha V=& g^{i\bar j}g^{k\bar l}g^{p\bar q}g^{m\bar n}
\left [ u_{;i\bar qk\bar n\alpha}u_{;\bar jp\bar lm}+
u_{;i\bar qk\bar n}u_{;\bar jp\bar lm\alpha}\right ]\\
&-g^{i\bar j}g^{k\bar l}g^{p\bar q}g^{m\bar n}
g^{\gamma\bar\delta}\left [u_{;i\bar\delta\alpha}
u_{;\gamma\bar qk\bar n}u_{;\bar jp\bar lm}
+u_{;k\bar\delta\alpha}u_{;i\bar q\gamma\bar n}u_{;\bar jp\bar lm}
\right ]\\
&-g^{i\bar j}g^{k\bar l}g^{p\bar q}g^{m\bar n}
g^{\gamma\bar\delta}\left [u_{;p\bar\delta\alpha}
u_{;\bar j\gamma\bar lm}u_{;i\bar qk\bar n}
+u_{;m\bar\delta\alpha}
u_{;\bar jp\bar l\gamma}u_{;i\bar qk\bar n}
\right ]\\
&+ g^{i\bar j}g^{k\bar l}g^{p\bar q}g^{m\bar n}
\left [ u_{;i\bar nkp\alpha}u_{;\bar jm\bar l\bar q}
+ u_{;i\bar nkp}u_{;\bar jm\bar l\bar q\alpha}
\right ]\\
&-g^{i\bar j}g^{k\bar l}g^{p\bar q}g^{m\bar n}
g^{\gamma\bar\delta}\left [u_{;i\bar\delta\alpha}
u_{;\gamma\bar nkp}u_{;\bar jm\bar l\bar q}
+u_{;k\bar\delta\alpha}
u_{;i\bar n\gamma p}u_{;\bar jm\bar l\bar q}
\right ]\\
&-g^{i\bar j}g^{k\bar l}g^{p\bar q}g^{m\bar n}
g^{\gamma\bar\delta}\left [u_{;p\bar\delta\alpha}
u_{;i\bar nk\gamma}u_{;\bar jm\bar l\bar q}
+u_{;m\bar\delta\alpha}u_{;\bar j\gamma\bar l\bar q}
u_{;i\bar nkp}\right ].
\end{split}
\end{align}
By differentiating the above formula we have
\begin{eqnarray}\label{bdd200}
\Delta^{'}V=g^{\alpha\bar\beta}
\partial_\alpha\partial_{\bar\beta}V=A_1+A_2+A_3+A_4+A_5+A_6+A_7
\end{eqnarray}
where
\begin{align}\label{bdd210}
\begin{split}
A_1=& g^{i\bar j}g^{k\bar l}g^{p\bar q}g^{m\bar n}g^{\alpha\bar\beta}
\left [ u_{;i\bar qk\bar n\alpha}u_{;\bar jp\bar lm\bar\beta}
+u_{;i\bar qk\bar n\bar\beta}u_{;\bar jp\bar lm\alpha}\right ]\\
&+g^{i\bar j}g^{k\bar l}g^{p\bar q}g^{m\bar n}g^{\alpha\bar\beta}
\left [  u_{;i\bar nkp\alpha}u_{;\bar jm\bar l\bar q\bar\beta}
+ u_{;i\bar nkp\bar\beta}u_{;\bar jm\bar l\bar q\alpha}\right ]\\
=&W,
\end{split}
\end{align}
\begin{align}\label{bdd220}
\begin{split}
A_2=& g^{i\bar j}g^{k\bar l}g^{p\bar q}g^{m\bar n}g^{\alpha\bar\beta}
\left [  u_{;i\bar qk\bar n\alpha\bar\beta}u_{;\bar jp\bar lm}
+ u_{;i\bar qk\bar n}u_{;\bar jp\bar lm\alpha\bar\beta}\right ] \\
&+g^{i\bar j}g^{k\bar l}g^{p\bar q}g^{m\bar n}g^{\alpha\bar\beta}
\left [  u_{;i\bar nkp\alpha\bar\beta}u_{;\bar jm\bar l\bar q}
+u_{;i\bar nkp}u_{;\bar jm\bar l\bar q\alpha\bar\beta}\right ],
\end{split}
\end{align}
\begin{align}\label{bdd230}
\begin{split}
A_3=&-g^{i\bar j}g^{k\bar l}g^{p\bar q}g^{m\bar n}g^{\alpha\bar\beta}
g^{\gamma\bar\delta}\left [
u_{;i\bar\delta k\bar n\alpha}u_{;\bar q\gamma\bar\beta}
u_{;\bar jp\bar lm}
+u_{;i\bar q k\bar\delta\alpha}u_{;\bar n\gamma\bar\beta}
u_{;\bar jp\bar lm}\right ]\\
&-g^{i\bar j}g^{k\bar l}g^{p\bar q}g^{m\bar n}g^{\alpha\bar\beta}
g^{\gamma\bar\delta}\left [
u_{;i\bar q k\bar n\alpha}u_{;\bar j\gamma\bar\beta}
u_{;\bar\delta p\bar lm}
+u_{;i\bar q k\bar n\alpha}u_{;\bar l\gamma\bar\beta}
u_{;\bar j p\bar\delta m}\right ]\\
&-g^{i\bar j}g^{k\bar l}g^{p\bar q}g^{m\bar n}g^{\alpha\bar\beta}
g^{\gamma\bar\delta}\left [
u_{;i\bar\delta k\bar n}u_{;\bar q\gamma\bar\beta}
u_{;\bar jp\bar lm\alpha}
+u_{;i\bar q k\bar\delta}u_{;\bar n\gamma\bar\beta}
u_{;\bar jp\bar lm\alpha}\right ]\\
&-g^{i\bar j}g^{k\bar l}g^{p\bar q}g^{m\bar n}g^{\alpha\bar\beta}
g^{\gamma\bar\delta}\left [
u_{;i\bar q k\bar n}u_{;\bar j\gamma\bar\beta}
u_{;\bar\delta p\bar lm\alpha}
+u_{;i\bar q k\bar n}u_{;\bar l\gamma\bar\beta}
u_{;\bar jp\bar\delta m\alpha}\right ]\\
&-g^{i\bar j}g^{k\bar l}g^{p\bar q}g^{m\bar n}g^{\alpha\bar\beta}
g^{\gamma\bar\delta}\left [
u_{;\gamma\bar q k\bar n\bar\beta}u_{;i\bar\delta\alpha}
u_{;\bar jp\bar lm}
+u_{;\bar jp\bar lm\bar\beta}u_{;i\bar\delta\alpha}
u_{;\gamma\bar q k\bar n}\right ]\\
&-g^{i\bar j}g^{k\bar l}g^{p\bar q}g^{m\bar n}g^{\alpha\bar\beta}
g^{\gamma\bar\delta}\left [
u_{;i\bar q\gamma\bar n\bar\beta}u_{;k\bar\delta\alpha}
u_{;\bar jp\bar lm}
+u_{;\bar jp\bar lm\bar\beta}u_{;k\bar\delta\alpha}
u_{;i\bar q\gamma\bar n}\right ]\\
&-g^{i\bar j}g^{k\bar l}g^{p\bar q}g^{m\bar n}g^{\alpha\bar\beta}
g^{\gamma\bar\delta}\left [
u_{;i\bar q k\bar n\bar\beta}u_{;p\bar\delta\alpha}
u_{;\bar j\gamma\bar lm}
+u_{;\bar j\gamma\bar lm\bar\beta}u_{;p\bar\delta\alpha}
u_{;i\bar q k\bar n}\right ]\\
&-g^{i\bar j}g^{k\bar l}g^{p\bar q}g^{m\bar n}g^{\alpha\bar\beta}
g^{\gamma\bar\delta}\left [
u_{;i\bar q k\bar n\bar\beta}u_{;m\bar\delta\alpha}
u_{;\bar j p\bar l\gamma}
+u_{;\bar j p\bar l\gamma\bar\beta}u_{;m\bar\delta\alpha}
u_{;i\bar q k\bar n}\right ],
\end{split}
\end{align}
\begin{align}\label{bdd240}
\begin{split}
A_4=&-g^{i\bar j}g^{k\bar l}g^{p\bar q}g^{m\bar n}g^{\alpha\bar\beta}
g^{\gamma\bar\delta}\left [
u_{;i\bar\delta kp\alpha}u_{;\bar n\gamma\bar\beta}
u_{;\bar jm\bar l\bar q}
+u_{;i\bar n kp\alpha}u_{;\bar j\gamma\bar\beta}
u_{;\bar\delta m\bar l\bar q}\right ]\\
&-g^{i\bar j}g^{k\bar l}g^{p\bar q}g^{m\bar n}g^{\alpha\bar\beta}
g^{\gamma\bar\delta}\left [
u_{;i\bar n kp\alpha}u_{;\bar l\gamma\bar\beta}
u_{;\bar j m\bar\delta\bar q}
+u_{;i\bar n kp\alpha}u_{;\bar q\gamma\bar\beta}
u_{;\bar j m\bar l\bar\delta}\right ]\\
&-g^{i\bar j}g^{k\bar l}g^{p\bar q}g^{m\bar n}g^{\alpha\bar\beta}
g^{\gamma\bar\delta}\left [
u_{;i\bar\delta kp}u_{;\bar n\gamma\bar\beta}
u_{;\bar jm\bar l\bar q\alpha}
+u_{;i\bar n kp}u_{;\bar j\gamma\bar\beta}
u_{;\bar\delta m\bar l\bar q\alpha}\right ]\\
&-g^{i\bar j}g^{k\bar l}g^{p\bar q}g^{m\bar n}g^{\alpha\bar\beta}
g^{\gamma\bar\delta}\left [
u_{;i\bar n kp}u_{;\bar l\gamma\bar\beta}
u_{;\bar j m\bar\delta\bar q\alpha}
+u_{;i\bar n kp}u_{;\bar q\gamma\bar\beta}
u_{;\bar j m\bar l\bar\delta\alpha}\right ]\\
&-g^{i\bar j}g^{k\bar l}g^{p\bar q}g^{m\bar n}g^{\alpha\bar\beta}
g^{\gamma\bar\delta}\left [
u_{;\gamma\bar n kp\bar\beta}u_{;i\bar\delta\alpha}
u_{;\bar j m\bar l\bar q}
+u_{;\bar j m\bar l\bar q\bar\beta}u_{;i\bar\delta\alpha}
u_{;\gamma\bar n kp}\right ]\\
&-g^{i\bar j}g^{k\bar l}g^{p\bar q}g^{m\bar n}g^{\alpha\bar\beta}
g^{\gamma\bar\delta}\left [
u_{;i\bar n \gamma p\bar\beta}u_{;k\bar\delta\alpha}
u_{;\bar j m\bar l\bar q}
+u_{;\bar j m\bar l\bar q\bar\beta}u_{;k\bar\delta\alpha}
u_{;i\bar n\gamma p}\right ]\\
&-g^{i\bar j}g^{k\bar l}g^{p\bar q}g^{m\bar n}g^{\alpha\bar\beta}
g^{\gamma\bar\delta}\left [
u_{;i\bar nk \gamma\bar\beta}u_{;p\bar\delta\alpha}
u_{;\bar j m\bar l\bar q}
+u_{;\bar j m\bar l\bar q\bar\beta}u_{;p\bar\delta\alpha}
u_{;i\bar n k\gamma}\right ]\\
&-g^{i\bar j}g^{k\bar l}g^{p\bar q}g^{m\bar n}g^{\alpha\bar\beta}
g^{\gamma\bar\delta}\left [
u_{;i\bar nkp\bar\beta}u_{;m\bar\delta\alpha}
u_{;\bar j\gamma\bar l\bar q}
+u_{;\bar j\gamma\bar l\bar q\bar\beta}u_{;m\bar\delta\alpha}
u_{;i\bar n kp}\right ],
\end{split}
\end{align}
\begin{align}\label{bdd250}
\begin{split}
A_5=&-g^{i\bar j}g^{k\bar l}g^{p\bar q}g^{m\bar n}g^{\alpha\bar\beta}
g^{\gamma\bar\delta}\left [
u_{;i\bar\delta\alpha\bar\beta}u_{;\gamma\bar qk\bar n}
u_{;\bar jp\bar lm}
+u_{;k\bar\delta\alpha\bar\beta}u_{;i\bar q\gamma\bar n}
u_{;\bar jp\bar lm}\right ]\\
&-g^{i\bar j}g^{k\bar l}g^{p\bar q}g^{m\bar n}g^{\alpha\bar\beta}
g^{\gamma\bar\delta}\left [
u_{;p\bar\delta\alpha\bar\beta}u_{;i\bar qk\bar n}
u_{;\bar j\gamma\bar lm}
+u_{;m\bar\delta\alpha\bar\beta}u_{;i\bar qk\bar n}
u_{;\bar jp\bar l\gamma}\right ]\\
&-g^{i\bar j}g^{k\bar l}g^{p\bar q}g^{m\bar n}g^{\alpha\bar\beta}
g^{\gamma\bar\delta}\left [
u_{;i\bar\delta\alpha\bar\beta}u_{;\gamma\bar nkp}
u_{;\bar jm\bar l\bar q}
+u_{;k\bar\delta\alpha\bar\beta}u_{;i\bar n\gamma p}
u_{;\bar jm\bar l\bar q}\right ]\\
&-g^{i\bar j}g^{k\bar l}g^{p\bar q}g^{m\bar n}g^{\alpha\bar\beta}
g^{\gamma\bar\delta}\left [
u_{;p\bar\delta\alpha\bar\beta}u_{;i\bar nk\gamma}
u_{;\bar jm\bar l\bar q}
+u_{;m\bar\delta\alpha\bar\beta}u_{;i\bar nkp}
u_{;\bar j\gamma\bar l\bar q}\right ],
\end{split}
\end{align}
\begin{align}\label{bdd260}
\begin{split}
A_6=&g^{i\bar j}g^{k\bar l}g^{p\bar q}g^{m\bar n}g^{\alpha\bar\beta}
g^{\gamma\bar\delta}g^{s\bar t}\left [
u_{;i\bar t\alpha}u_{;\bar\delta s\bar\beta}
u_{;\gamma\bar qk\bar n}u_{;\bar jp\bar lm}
+u_{;i\bar\delta\alpha}u_{;\bar q s\bar\beta}
u_{;\gamma\bar tk\bar n}u_{;\bar jp\bar lm}\right ]\\
&+g^{i\bar j}g^{k\bar l}g^{p\bar q}g^{m\bar n}g^{\alpha\bar\beta}
g^{\gamma\bar\delta}g^{s\bar t}\left [
u_{;i\bar\delta\alpha}u_{;\bar n s\bar\beta}
u_{;\gamma\bar qk\bar t}u_{;\bar jp\bar lm}
+u_{;i\bar\delta\alpha}u_{;\bar j s\bar\beta}
u_{;\gamma\bar qk\bar n}u_{;\bar tp\bar lm}\right ]\\
&+g^{i\bar j}g^{k\bar l}g^{p\bar q}g^{m\bar n}g^{\alpha\bar\beta}
g^{\gamma\bar\delta}g^{s\bar t}\left [
u_{;i\bar\delta\alpha}u_{;\bar l s\bar\beta}
u_{;\gamma\bar qk\bar n}u_{;\bar jp\bar tm}
+u_{;k\bar t\alpha}u_{;\bar\delta s\bar\beta}
u_{;i\bar q\gamma\bar n}u_{;\bar jp\bar lm}\right ]\\
&+g^{i\bar j}g^{k\bar l}g^{p\bar q}g^{m\bar n}g^{\alpha\bar\beta}
g^{\gamma\bar\delta}g^{s\bar t}\left [
u_{;k\bar\delta\alpha}u_{;\bar q s\bar\beta}
u_{;i\bar t\gamma\bar n}u_{;\bar jp\bar lm}
+u_{;k\bar\delta\alpha}u_{;\bar n s\bar\beta}
u_{;i\bar q\gamma\bar t}u_{;\bar jp\bar lm}\right ]\\
&+g^{i\bar j}g^{k\bar l}g^{p\bar q}g^{m\bar n}g^{\alpha\bar\beta}
g^{\gamma\bar\delta}g^{s\bar t}\left [
u_{;k\bar\delta\alpha}u_{;\bar j s\bar\beta}
u_{;i\bar q\gamma\bar n}u_{;\bar tp\bar lm}
+u_{;k\bar\delta\alpha}u_{;\bar l s\bar\beta}
u_{;i\bar q\gamma\bar n}u_{;\bar jp\bar tm}\right ]\\
&+g^{i\bar j}g^{k\bar l}g^{p\bar q}g^{m\bar n}g^{\alpha\bar\beta}
g^{\gamma\bar\delta}g^{s\bar t}\left [
u_{;p\bar t\alpha}u_{;\bar\delta s\bar\beta}
u_{;i\bar q k\bar n}u_{;\bar j\gamma\bar lm}
+u_{;p\bar\delta\alpha}u_{;\bar j s\bar\beta}
u_{;i\bar q k\bar n}u_{;\bar t\gamma\bar lm}\right ]\\
&+g^{i\bar j}g^{k\bar l}g^{p\bar q}g^{m\bar n}g^{\alpha\bar\beta}
g^{\gamma\bar\delta}g^{s\bar t}\left [
u_{;p\bar\delta\alpha}u_{;\bar l s\bar\beta}
u_{;i\bar q k\bar n}u_{;\bar j\gamma\bar tm}
+u_{;p\bar\delta\alpha}u_{;\bar q s\bar\beta}
u_{;i\bar t k\bar n}u_{;\bar j\gamma\bar lm}\right ]\\
&+g^{i\bar j}g^{k\bar l}g^{p\bar q}g^{m\bar n}g^{\alpha\bar\beta}
g^{\gamma\bar\delta}g^{s\bar t}\left [
u_{;p\bar\delta\alpha}u_{;\bar n s\bar\beta}
u_{;i\bar q k\bar t}u_{;\bar j\gamma\bar lm}
+u_{;m\bar t\alpha}u_{;\bar\delta s\bar\beta}
u_{;i\bar q k\bar n}u_{;\bar j p\bar l\gamma}\right ]\\
&+g^{i\bar j}g^{k\bar l}g^{p\bar q}g^{m\bar n}g^{\alpha\bar\beta}
g^{\gamma\bar\delta}g^{s\bar t}\left [
u_{;m\bar\delta\alpha}u_{;\bar j s\bar\beta}
u_{;i\bar q k\bar n}u_{;\bar t p\bar l\gamma}
+u_{;m\bar\delta\alpha}u_{;\bar l s\bar\beta}
u_{;i\bar q k\bar n}u_{;\bar j p\bar t\gamma}\right ]\\
&+g^{i\bar j}g^{k\bar l}g^{p\bar q}g^{m\bar n}g^{\alpha\bar\beta}
g^{\gamma\bar\delta}g^{s\bar t}\left [
u_{;m\bar\delta\alpha}u_{;\bar q s\bar\beta}
u_{;i\bar t k\bar n}u_{;\bar t j\bar l\gamma}
+u_{;m\bar\delta\alpha}u_{;\bar n s\bar\beta}
u_{;i\bar q k\bar t}u_{;\bar j p\bar l\gamma}\right ]
\end{split}
\end{align}
and
\begin{align}\label{bdd270}
\begin{split}
A_7=&g^{i\bar j}g^{k\bar l}g^{p\bar q}g^{m\bar n}g^{\alpha\bar\beta}
g^{\gamma\bar\delta}g^{s\bar t}\left [
u_{;i\bar t\alpha}u_{;\bar\delta s\bar\beta}
u_{;\gamma\bar n kp}u_{;\bar j m\bar l\bar q}
+u_{;i\bar\delta\alpha}u_{;\bar n s\bar\beta}
u_{;\gamma\bar t kp}u_{;\bar j m\bar l\bar q}\right ]\\
&+g^{i\bar j}g^{k\bar l}g^{p\bar q}g^{m\bar n}g^{\alpha\bar\beta}
g^{\gamma\bar\delta}g^{s\bar t}\left [
u_{;i\bar\delta\alpha}u_{;\bar j s\bar\beta}
u_{;\gamma\bar n kp}u_{;\bar t m\bar l\bar q}
+u_{;i\bar\delta\alpha}u_{;\bar l s\bar\beta}
u_{;\gamma\bar n kp}u_{;\bar j m\bar t\bar q}\right ]\\
&+g^{i\bar j}g^{k\bar l}g^{p\bar q}g^{m\bar n}g^{\alpha\bar\beta}
g^{\gamma\bar\delta}g^{s\bar t}\left [
u_{;i\bar\delta\alpha}u_{;\bar q s\bar\beta}
u_{;\gamma\bar n kp}u_{;\bar j m\bar l\bar t}
+u_{;k\bar t\alpha}u_{;\bar\delta s\bar\beta}
u_{;i\bar n\gamma p}u_{;\bar j m\bar l\bar q}\right ]\\
&+g^{i\bar j}g^{k\bar l}g^{p\bar q}g^{m\bar n}g^{\alpha\bar\beta}
g^{\gamma\bar\delta}g^{s\bar t}\left [
u_{;k\bar\delta\alpha}u_{;\bar n s\bar\beta}
u_{;i\bar t\gamma p}u_{;\bar j m\bar l\bar q}
+u_{;k\bar\delta\alpha}u_{;\bar j s\bar\beta}
u_{;i\bar n\gamma p}u_{;\bar t m\bar l\bar q}\right ]\\
&+g^{i\bar j}g^{k\bar l}g^{p\bar q}g^{m\bar n}g^{\alpha\bar\beta}
g^{\gamma\bar\delta}g^{s\bar t}\left [
u_{;k\bar\delta\alpha}u_{;\bar l s\bar\beta}
u_{;i\bar n\gamma p}u_{;\bar j m\bar t\bar q}
+u_{;k\bar\delta\alpha}u_{;\bar q s\bar\beta}
u_{;i\bar n\gamma p}u_{;\bar j m\bar l\bar t}\right ]\\
&+g^{i\bar j}g^{k\bar l}g^{p\bar q}g^{m\bar n}g^{\alpha\bar\beta}
g^{\gamma\bar\delta}g^{s\bar t}\left [
u_{;p\bar t\alpha}u_{;\bar\delta s\bar\beta}
u_{;i\bar n k\gamma}u_{;\bar j m\bar l\bar q}
+u_{;p\bar\delta\alpha}u_{;\bar n s\bar\beta}
u_{;i\bar tk\gamma}u_{;\bar j m\bar l\bar q}\right ]\\
&+g^{i\bar j}g^{k\bar l}g^{p\bar q}g^{m\bar n}g^{\alpha\bar\beta}
g^{\gamma\bar\delta}g^{s\bar t}\left [
u_{;p\bar\delta\alpha}u_{;\bar j s\bar\beta}
u_{;i\bar n k\gamma}u_{;\bar t m\bar l\bar q}
+u_{;p\bar\delta\alpha}u_{;\bar l s\bar\beta}
u_{;i\bar nk\gamma}u_{;\bar j m\bar t\bar q}\right ]\\
&+g^{i\bar j}g^{k\bar l}g^{p\bar q}g^{m\bar n}g^{\alpha\bar\beta}
g^{\gamma\bar\delta}g^{s\bar t}\left [
u_{;p\bar\delta\alpha}u_{;\bar q s\bar\beta}
u_{;i\bar n k\gamma}u_{;\bar j m\bar l\bar t}
+u_{;m\bar t\alpha}u_{;\bar\delta s\bar\beta}
u_{;i\bar n kp}u_{;\bar j\gamma\bar l\bar q}\right ]\\
&+g^{i\bar j}g^{k\bar l}g^{p\bar q}g^{m\bar n}g^{\alpha\bar\beta}
g^{\gamma\bar\delta}g^{s\bar t}\left [
u_{;m\bar\delta\alpha}u_{;\bar j s\bar\beta}
u_{;i\bar n kp}u_{;\bar t\gamma\bar l\bar q}
+u_{;m\bar\delta\alpha}u_{;\bar l s\bar\beta}
u_{;i\bar n kp}u_{;\bar j\gamma\bar t\bar q}\right ]\\
&+g^{i\bar j}g^{k\bar l}g^{p\bar q}g^{m\bar n}g^{\alpha\bar\beta}
g^{\gamma\bar\delta}g^{s\bar t}\left [
u_{;m\bar\delta\alpha}u_{;\bar q s\bar\beta}
u_{;i\bar n kp}u_{;\bar j\gamma\bar l\bar t}
+u_{;m\bar\delta\alpha}u_{;\bar n s\bar\beta}
u_{;i\bar t kp}u_{;\bar j\gamma\bar l\bar q}\right ].
\end{split}
\end{align}

Now we estimate each $A_i$ in the sum of $\Delta^{'}V$. Since we
have $C^3$ estimate, that is, $S$ is bounded, by diagonalizing,
we know that each term in the sum $A_6$ and $A_7$ is bounded by
a constant multiple of $V$. So we have
\begin{eqnarray}\label{bdd290}
\left | A_6\right |+\left | A_7\right |
\leq C_{14}V.
\end{eqnarray}

Now we estimate terms in $A_5$. We have
\[
u_{;i\bar\delta\alpha\bar\beta}=u_{;\alpha\bar\beta i\bar\delta}
+(u_{;i\bar\delta\alpha\bar\beta}-u_{;\alpha\bar\beta i\bar\delta})
=u_{;\alpha\bar\beta i\bar\delta}+h^{s\bar t}
(u_{s\bar\beta}R_{\alpha\bar ti\bar\delta}-
u_{s\bar\delta}R_{i\bar t\alpha\bar\beta}).
\]
By using equation \eqref{bdd70} we have
\begin{align*}
\begin{split}
g^{\alpha\bar\beta}u_{;i\bar\delta\alpha\bar\beta}=&
g^{\alpha\bar\beta}u_{;\alpha\bar\beta i\bar\delta}+
g^{\alpha\bar\beta}h^{s\bar t}
(u_{s\bar\beta}R_{\alpha\bar ti\bar\delta}-
u_{s\bar\delta}R_{i\bar t\alpha\bar\beta})\\
=&\mathcal F_{i\bar\delta}+g^{\alpha\bar\beta}g^{s\bar t}
u_{;\alpha\bar ti}u_{;\bar\beta s\bar\delta}
+g^{\alpha\bar\beta}h^{s\bar t}
(u_{s\bar\beta}R_{\alpha\bar ti\bar\delta}-
u_{s\bar\delta}R_{i\bar t\alpha\bar\beta}).
\end{split}
\end{align*}
Since the curvature of the metric $h$ is bounded and we have
$C^2$ and $C^3$ estimate, we know that
\[
\left |-g^{i\bar j}g^{k\bar l}g^{p\bar q}
g^{m\bar n}g^{\alpha\bar\beta}g^{\gamma\bar\delta}
u_{;i\bar\delta\alpha\bar\beta}u_{;\gamma\bar qk\bar n}
u_{;\bar jp\bar lm}\right |\leq CV.
\]
Similarly we can compute other terms in the sum $A_5$. So
we have
\begin{eqnarray}\label{bdd300}
\left |A_5\right |\leq C_{15}V.
\end{eqnarray}
By combining formulas \eqref{bdd200}, \eqref{bdd290} and
\eqref{bdd300} we have
\begin{eqnarray}\label{bdd310}
\Delta^{'}V\geq A_1+A_2+A_3+A_4-C_{16}V.
\end{eqnarray}

Now we deal with terms in $A_2$. By differentiating formula
\eqref{bdd80} in a holomorphic direction or an anti-holomorphic
direction we have
\begin{align}\label{bdd320}
\begin{split}
g^{i\bar j}u_{;i\bar jk\bar l\alpha\bar\beta}=&
g^{i\bar j}g^{p\bar q}\left [
u_{;i\bar qk\bar l\alpha}u_{;\bar jp\bar\beta}
+u_{;i\bar qk\alpha\bar\beta}u_{;\bar jp\bar l}
+u_{;p\bar jk\bar l\bar\beta}u_{;i\bar q\alpha}
+u_{;\bar jp\bar l\alpha\bar\beta}u_{;i\bar qk}
\right ]\\
&+g^{i\bar j}g^{p\bar q}\left [
u_{;p\bar jk\bar l}u_{;i\bar q\alpha\bar\beta}
+u_{;\bar jp\bar l\alpha}u_{;i\bar qk\bar\beta}
+u_{;i\bar qk\alpha}u_{;\bar jp\bar l\bar\beta}
\right ]\\
&-g^{i\bar j}g^{p\bar q}g^{m\bar n}\left [
u_{;p\bar nk\bar l}u_{;\bar jm\bar\beta}u_{;i\bar q\alpha}
+u_{;p\bar jk\bar l}u_{;\bar qm\bar\beta}u_{;i\bar n\alpha}
+u_{;\bar np\bar l\alpha}u_{;\bar jm\bar\beta}u_{;i\bar qk}
\right ]\\
&-g^{i\bar j}g^{p\bar q}g^{m\bar n}\left [
u_{;\bar jp\bar l\alpha}u_{;\bar qm\bar\beta}u_{;i\bar nk}
+u_{;i\bar nk\alpha}u_{;\bar qm\bar\beta}u_{;\bar jp\bar l}
+u_{;i\bar qk\alpha}u_{;\bar jm\bar\beta}u_{;\bar np\bar l}
\right ]\\
&-g^{i\bar j}g^{p\bar q}g^{m\bar n}\left [
u_{;\bar jm\bar l\bar\beta}u_{;p\bar n\alpha}u_{;i\bar qk}
+u_{;\bar jp\bar l\bar\beta}u_{;i\bar n\alpha}u_{;m\bar qk}
+u_{;p\bar n\alpha\bar\beta}u_{;\bar jm\bar l}u_{;i\bar qk}
\right ]\\
&-g^{i\bar j}g^{p\bar q}g^{m\bar n}\left [
u_{;i\bar qk\bar\beta}u_{;\bar jm\bar l}u_{;p\bar n\alpha}
+u_{;m\bar qk\bar\beta}u_{;\bar jp\bar l}u_{;i\bar n\alpha}
+u_{;i\bar n\alpha\bar\beta}u_{;\bar jp\bar l}u_{;m\bar qk}
\right ]\\
&+g^{i\bar j}g^{p\bar q}g^{m\bar n}g^{s\bar t}\left [
u_{;i\bar qk}u_{;\bar js\bar\beta}
u_{;p\bar n\alpha}u_{;\bar tm\bar l}
+u_{;i\bar qk}u_{;\bar ns\bar\beta}
u_{;p\bar t\alpha}u_{;\bar jm\bar l}
\right ]\\
&+g^{i\bar j}g^{p\bar q}g^{m\bar n}g^{s\bar t}\left [
u_{;i\bar tk}u_{;\bar qs\bar\beta}
u_{;p\bar n\alpha}u_{;\bar jm\bar l}
+u_{;m\bar tk}u_{;\bar qs\bar\beta}
u_{;i\bar n\alpha}u_{;\bar jp\bar l}
\right ]\\
&+g^{i\bar j}g^{p\bar q}g^{m\bar n}g^{s\bar t}\left [
u_{;m\bar qk}u_{;\bar ns\bar\beta}
u_{;i\bar t\alpha}u_{;\bar jp\bar l}
+u_{;m\bar qk}u_{;\bar js\bar\beta}
u_{;i\bar n\alpha}u_{;\bar tp\bar l}
\right ]\\
&+u_{;k\bar l\alpha\bar\beta}+F_{;k\bar l\alpha\bar\beta}
\end{split}
\end{align}
and
\begin{align}\label{bdd330}
\begin{split}
g^{i\bar j}u_{;i\bar jk\bar l\alpha\gamma}=&
g^{i\bar j}g^{p\bar q}\left [
u_{;p\bar jk\bar l\alpha}u_{;i\bar q\gamma}
+u_{;p\bar jk\bar l\gamma}u_{;i\bar q\alpha}
+u_{;i\bar qk\alpha\gamma}u_{;\bar jp\bar l}
+u_{;\bar jp\bar l\alpha\gamma}u_{;i\bar qk}\right ]\\
&+g^{i\bar j}g^{p\bar q}\left [
u_{;p\bar jk\bar l}u_{;i\bar q\alpha\gamma}
+u_{;\bar jp\bar l\alpha}u_{;i\bar qk\gamma}
+u_{;i\bar qk\alpha}u_{;\bar jp\bar l\gamma}\right ]\\
&-g^{i\bar j}g^{p\bar q}g^{m\bar n}\left [
u_{;m\bar jk\bar l}u_{;p\bar n\gamma}u_{;i\bar q\alpha}
+u_{;p\bar jk\bar l}u_{;i\bar n\gamma}u_{;m\bar q\alpha}
+u_{;\bar jm\bar l\alpha}u_{;p\bar n\gamma}u_{;i\bar qk}
\right ]\\
&-g^{i\bar j}g^{p\bar q}g^{m\bar n}\left [
u_{;\bar jp\bar l\alpha}u_{;i\bar n\gamma}u_{;m\bar qk}
+u_{;m\bar qk\alpha}u_{;i\bar n\gamma}u_{;\bar jp\bar l}
+u_{;i\bar qk\alpha}u_{;p\bar n\gamma}u_{;\bar jm\bar l}
\right ]\\
&-g^{i\bar j}g^{p\bar q}g^{m\bar n}\left [
u_{;i\bar qk\gamma}u_{;p\bar n\alpha}u_{;\bar jm\bar l}
+u_{;p\bar n\alpha\gamma}u_{;i\bar qk}u_{;\bar jm\bar l}
+u_{;\bar jm\bar l\gamma}u_{;i\bar qk}u_{;p\bar n\alpha}
\right ]\\
&-g^{i\bar j}g^{p\bar q}g^{m\bar n}\left [
u_{;m\bar qk\gamma}u_{;i\bar n\alpha}u_{;\bar jp\bar l}
+u_{;i\bar n\alpha\gamma}u_{;m\bar qk}u_{;\bar jp\bar l}
+u_{;\bar jp\bar l\gamma}u_{;m\bar qk}u_{;i\bar n\alpha}
\right ]\\
&+g^{i\bar j}g^{p\bar q}g^{m\bar n}g^{s\bar t}\left [
u_{;i\bar t\gamma}u_{;s\bar qk}
u_{;p\bar n\alpha}u_{;\bar jm\bar l}
+u_{;p\bar t\gamma}u_{;i\bar qk}
u_{;s\bar n\alpha}u_{;\bar jm\bar l}\right ]\\
&+g^{i\bar j}g^{p\bar q}g^{m\bar n}g^{s\bar t}\left [
u_{;m\bar t\gamma}u_{;i\bar qk}
u_{;p\bar n\alpha}u_{;\bar js\bar l}
+u_{;m\bar t\gamma}u_{;s\bar qk}
u_{;i\bar n\alpha}u_{;\bar jp\bar l}\right ]\\
&+g^{i\bar j}g^{p\bar q}g^{m\bar n}g^{s\bar t}\left [
u_{;i\bar t\gamma}u_{;m\bar qk}
u_{;s\bar n\alpha}u_{;\bar jp\bar l}
+u_{;p\bar t\gamma}u_{;m\bar qk}
u_{;i\bar n\alpha}u_{;\bar js\bar l}\right ]\\
&+u_{;k\bar l\alpha\gamma}+F_{;k\bar l\alpha\gamma}.
\end{split}
\end{align}
By using a similar computation as in \cite{yau3} we have
\begin{eqnarray}\label{bdd340}
u_{;i\bar qk\bar n\alpha\bar\beta}\overset{4}{\simeq}
u_{;\alpha\bar\beta i\bar qk\bar n},
\end{eqnarray}
\begin{eqnarray}\label{bdd345}
u_{;\bar jp\bar lm\alpha\bar\beta}\overset{4}{\simeq}
u_{;\bar\beta\alpha\bar jp\bar lm},
\end{eqnarray}
\begin{eqnarray}\label{bdd350}
u_{;i\bar nkp\alpha\bar\beta}\overset{4}{\simeq}
u_{;\alpha\bar\beta i\bar nkp}
\end{eqnarray}
and
\begin{eqnarray}\label{bdd355}
u_{;\bar jm\bar l\bar q\alpha\bar\beta}\overset{4}{\simeq}
u_{;\bar\beta\alpha\bar jm\bar l\bar q}
\end{eqnarray}
which imply that
\begin{eqnarray}\label{bdd360}
g^{i\bar j}g^{k\bar l}g^{p\bar q}g^{m\bar n}g^{\alpha\bar\beta}
u_{;i\bar qk\bar n\alpha\bar\beta}u_{;\bar jp\bar lm}
\overset{4}{\cong}
g^{i\bar j}g^{k\bar l}g^{p\bar q}g^{m\bar n}g^{\alpha\bar\beta}
u_{;\alpha\bar\beta i\bar qk\bar n}u_{;\bar jp\bar lm},
\end{eqnarray}
\begin{eqnarray}\label{bdd365}
g^{i\bar j}g^{k\bar l}g^{p\bar q}g^{m\bar n}g^{\alpha\bar\beta}
u_{;\bar jp\bar lm\alpha\bar\beta}u_{;i\bar qk\bar n}
\overset{4}{\cong}
g^{i\bar j}g^{k\bar l}g^{p\bar q}g^{m\bar n}g^{\alpha\bar\beta}
u_{;\bar\beta\alpha\bar jp\bar lm}u_{;i\bar qk\bar n},
\end{eqnarray}
\begin{eqnarray}\label{bdd370}
g^{i\bar j}g^{k\bar l}g^{p\bar q}g^{m\bar n}g^{\alpha\bar\beta}
u_{;i\bar nkp\alpha\bar\beta}u_{;\bar jm\bar l\bar q}
\overset{4}{\cong}
g^{i\bar j}g^{k\bar l}g^{p\bar q}g^{m\bar n}g^{\alpha\bar\beta}
u_{;\alpha\bar\beta i\bar nkp}u_{;\bar jm\bar l\bar q}
\end{eqnarray}
and
\begin{eqnarray}\label{bdd372}
g^{i\bar j}g^{k\bar l}g^{p\bar q}g^{m\bar n}g^{\alpha\bar\beta}
u_{;\bar jm\bar l\bar q\alpha\bar\beta}u_{;i\bar nkp}
\overset{4}{\cong}
g^{i\bar j}g^{k\bar l}g^{p\bar q}g^{m\bar n}g^{\alpha\bar\beta}
u_{;\bar\beta\alpha\bar jm\bar l\bar q}u_{;i\bar nkp}.
\end{eqnarray}

By using equations \eqref{bdd320}, \eqref{bdd330} and their
conjugations, we have
\begin{eqnarray}\label{bdd374}
\left |
g^{i\bar j}g^{k\bar l}g^{p\bar q}g^{m\bar n}g^{\alpha\bar\beta}
u_{;\alpha\bar\beta i\bar qk\bar n}u_{;\bar jp\bar lm}-T_1\right |
\leq C_{17}V^{\frac{3}{2}}+C_{18}V+C_{19}V^{\frac{1}{2}},
\end{eqnarray}
\begin{eqnarray}\label{bdd3745}
\left |
g^{i\bar j}g^{k\bar l}g^{p\bar q}g^{m\bar n}g^{\alpha\bar\beta}
u_{;\bar\beta\alpha\bar jp\bar lm}u_{;i\bar qk\bar n}-T_2\right |
\leq C_{17}V^{\frac{3}{2}}+C_{18}V+C_{19}V^{\frac{1}{2}},
\end{eqnarray}
\begin{eqnarray}\label{bdd376}
\left |
g^{i\bar j}g^{k\bar l}g^{p\bar q}g^{m\bar n}g^{\alpha\bar\beta}
u_{;\alpha\bar\beta i\bar nkp}u_{;\bar jm\bar l\bar q}-T_3\right |
\leq  C_{17}V^{\frac{3}{2}}+C_{18}V+C_{19}V^{\frac{1}{2}}
\end{eqnarray}
and
\begin{eqnarray}\label{bdd3765}
\left |
g^{i\bar j}g^{k\bar l}g^{p\bar q}g^{m\bar n}g^{\alpha\bar\beta}
u_{;\bar\beta\alpha\bar jm\bar l\bar q}u_{;i\bar nkp}-T_4\right |
\leq  C_{17}V^{\frac{3}{2}}+C_{18}V+C_{19}V^{\frac{1}{2}}
\end{eqnarray}
where
\begin{align}\label{bdd3782}
\begin{split}
T_1=& g^{i\bar j}g^{k\bar l}g^{p\bar q}g^{m\bar n}
g^{\alpha\bar\beta}g^{\gamma\bar\delta}\left [
u_{;i\bar qk\bar n\alpha}u_{;\bar jp\bar\delta}
u_{;\bar lm\bar\beta\gamma}
+u_{;i\bar nkp\bar\beta}u_{;\bar jm\bar\delta}
u_{;\bar l\gamma\bar q\alpha}
\right ]\\
&+g^{i\bar j}g^{k\bar l}g^{p\bar q}g^{m\bar n}
g^{\alpha\bar\beta}g^{\gamma\bar\delta}\left [
u_{;i\bar qk\bar n\bar\beta}u_{;p\bar j\gamma}
u_{;\bar lm\bar\delta\alpha}
+u_{;\bar jp\bar lm\bar\beta}u_{;i\bar q\gamma}
u_{;\bar\delta k\bar n\alpha}
\right ],
\end{split}
\end{align}
\begin{align}\label{bdd3784}
\begin{split}
T_2=& g^{i\bar j}g^{k\bar l}g^{p\bar q}g^{m\bar n}
g^{\alpha\bar\beta}g^{\gamma\bar\delta}\left [
u_{;\bar jp\bar lm\bar\beta}u_{;i\bar q\gamma}
u_{;k\bar n\alpha\bar\delta}
+u_{;\bar jm\bar l\bar q\alpha}u_{;i\bar n\gamma}
u_{;k\bar\delta p\bar\beta}
\right ]\\
&+g^{i\bar j}g^{k\bar l}g^{p\bar q}g^{m\bar n}
g^{\alpha\bar\beta}g^{\gamma\bar\delta}\left [
u_{;\bar jp\bar lm\alpha}u_{;\bar qi\bar\delta}
u_{;k\bar n\gamma\bar\beta}
+u_{;i\bar qk\bar n\alpha}u_{;\bar jp\bar\delta}
u_{;\gamma\bar lm\bar\beta}
\right ],
\end{split}
\end{align}
\begin{align}\label{bdd3786}
\begin{split}
T_3=& g^{i\bar j}g^{k\bar l}g^{p\bar q}g^{m\bar n}
g^{\alpha\bar\beta}g^{\gamma\bar\delta}\left [
u_{;i\bar qk\bar n\alpha}u_{;p\bar j\gamma}
u_{;\bar lm\bar\beta\bar\delta}
+u_{;i\bar nkp\alpha}u_{;\bar jm\bar\delta}
u_{;\bar l\gamma\bar q\bar\beta}
\right ]\\
&+g^{i\bar j}g^{k\bar l}g^{p\bar q}g^{m\bar n}
g^{\alpha\bar\beta}g^{\gamma\bar\delta}\left [
u_{;i\bar qk\bar n\alpha}u_{;p\bar j\gamma}
u_{;\bar lm\bar\delta\bar\beta}
+u_{;\bar jp\bar lm\alpha}u_{;i\bar q\gamma}
u_{;\bar\delta k\bar n\bar\beta}
\right ]
\end{split}
\end{align}
and
\begin{align}\label{bdd3788}
\begin{split}
T_4=& g^{i\bar j}g^{k\bar l}g^{p\bar q}g^{m\bar n}
g^{\alpha\bar\beta}g^{\gamma\bar\delta}\left [
u_{;\bar jp\bar lm\bar\beta}u_{;\bar qi\bar\delta}
u_{;k\bar n\alpha\gamma}
+u_{;\bar jm\bar l\bar q\bar\beta}u_{;i\bar n\gamma}
u_{;k\bar\delta p\alpha}
\right ]\\
&+g^{i\bar j}g^{k\bar l}g^{p\bar q}g^{m\bar n}
g^{\alpha\bar\beta}g^{\gamma\bar\delta}\left [
u_{;\bar jp\bar lm\bar\beta}u_{;\bar qi\bar\delta}
u_{;k\bar n\gamma\alpha}
+u_{;i\bar qk\bar n\bar\beta}u_{;\bar jp\bar\delta}
u_{;\gamma\bar lm\alpha}
\right ].
\end{split}
\end{align}

Now we choose local coordinates such that $g_{i\bar
j}=\delta_{ij}$. By combining formulas \eqref{bdd310},
\eqref{bdd230},  \eqref{bdd240} and
\eqref{bdd360}$-$\eqref{bdd3788} we have
\begin{align}\label{bdd390}
\begin{split}
\Delta^{'}V\geq & \sum_{i,p,k,m,\alpha}
\left [ \left |A_{ipkm\alpha}\right |^2
+\left |B_{ipkm\alpha}\right |^2+\left |C_{ipkm\alpha}\right |^2
+\left |D_{ipkm\alpha}\right |^2\right ]\\
& -C_{20}V^{\frac{3}{2}}-C_{21}V-C_{22}V^{\frac{1}{2}}\\
=&\widetilde W -C_{20}V^{\frac{3}{2}}-C_{21}V-C_{22}V^{\frac{1}{2}}.
\end{split}
\end{align}
Now we estimate $\left |\nabla^{'}V\right |$. For each fixed
$\alpha$, by \eqref{bdd190} we have
\begin{align}\label{bdd400}
\begin{split}
\partial_\alpha V=\sum_{i,p,k,m}\left [
A_{ipkm\alpha}u_{;\bar ip\bar km}
+B_{ipkm\alpha}u_{;i\bar pk\bar m}
+C_{ipkm\alpha}u_{;\bar im\bar k\bar p}
+D_{ipkm\alpha}u_{;i\bar mkp}\right ]+X_{\alpha}
\end{split}
\end{align}
where
\[
\left |X_\alpha\right |\leq C V.
\]
By using the Schwarz inequality, it is easy to see that
\[
\left | \partial_\alpha V\right | \leq \sqrt 2V^{\frac{1}{2}}
\widetilde W^{\frac{1}{2}}+CV
\]
which implies
\begin{eqnarray}\label{bdd410}
\left |\nabla^{'}V\right |\leq C_{23}\widetilde
W^{\frac{1}{2}}V^{\frac{1}{2}}+C_{24}V.
\end{eqnarray}

Now we estimate the derivatives of $S$. By using similar
computation as above, we have
\begin{align}\label{bdd420}
\begin{split}
\partial_\alpha S=&g^{i\bar j}g^{k\bar l}g^{p\bar q}\left (
u_{;i\bar qk\alpha}u_{;\bar jp\bar l}+u_{;\bar jp\bar l\alpha}
u_{;i\bar qk}\right )\\
&-g^{i\bar j}g^{k\bar l}g^{p\bar q}g^{m\bar n}\left (
u_{;m\bar qk}u_{;i\bar n\alpha}u_{;\bar jp\bar l}
+u_{;i\bar qm}u_{;k\bar n\alpha}u_{;\bar jp\bar l}
+u_{;i\bar qk}u_{;p\bar n\alpha}u_{;\bar jm\bar l}\right ).
\end{split}
\end{align}
Since $S$ is bounded, we have
\[
\left |\partial_\alpha S\right |\leq C V^{\frac{1}{2}}+\widetilde C
\]
which implies
\begin{eqnarray}\label{bdd430}
\left |\nabla^{'}S\right |\leq C_{25}V^{\frac{1}{2}}+C_{26}.
\end{eqnarray}
By Yau's work in \cite{yau3} we have
\begin{eqnarray}\label{bdd440}
\Delta^{'}S\geq V-C_{27}V^{\frac{1}{2}}-C_{28}.
\end{eqnarray}
From \eqref{bdd390}, we have
\begin{align}\label{bdd450}
\begin{split}
(S+\kappa)\Delta^{'}V\geq & (S+\kappa)\left (
\widetilde W -C_{20}V^{\frac{3}{2}}-C_{21}V-C_{22}V^{\frac{1}{2}}
\right )\\
\geq & \kappa\widetilde W
-C_{29}V^{\frac{3}{2}}-C_{30}V-C_{31}V^{\frac{1}{2}}
\end{split}
\end{align}
where the constants $C_{29}$,  $C_{30}$ and  $C_{31}$ depend on
$\kappa$. By \eqref{bdd440}, \eqref{bdd410} and \eqref{bdd440} we
also have
\begin{eqnarray}\label{bdd460}
V\Delta^{'}S\geq V^2-C_{27}V^{\frac{3}{2}}-C_{28}V
\end{eqnarray}
and
\begin{eqnarray}\label{bdd470}
2\left |\nabla^{'}S\right |\left |\nabla^{'}V\right |\leq
C_{32}\widetilde W^{\frac{1}{2}}V
+C_{33}\widetilde W^{\frac{1}{2}}V^{\frac{1}{2}}
+C_{34}V^{\frac{3}{2}}+C_{35}V.
\end{eqnarray}

Combining \eqref{bdd180}, \eqref{bdd450}, \eqref{bdd460} and
\eqref{bdd470} we have
\begin{align}\label{bdd480}
\begin{split}
\Delta^{'}\left [(S+\kappa)V\right ]\geq \kappa\widetilde W+V^2
-C_{32}\widetilde W^{\frac{1}{2}}V
-C_{33}\widetilde W^{\frac{1}{2}}V^{\frac{1}{2}}
-C_{36}V^{\frac{3}{2}}-C_{37}V-C_{38}V^{\frac{1}{2}}
\end{split}
\end{align}
where constants $C_{36}$,  $C_{37}$ and  $C_{38}$ depend on $\kappa$.
Now we fix a $\kappa$  such that
\[
\kappa\geq \max\{3C_{32}^2,3C_{33}^2,1\}.
\]
We have
\[
\frac{\kappa}{3}\widetilde W-C_{32}\widetilde W^{\frac{1}{2}}V
+\frac{1}{4}V^2\geq 0
\]
and
\[
\frac{\kappa}{3}\widetilde W
-C_{33}\widetilde W^{\frac{1}{2}}V^{\frac{1}{2}}\geq -\frac{1}{4}V.
\]
With this choice of $\kappa$, by formula \eqref{bdd480} we have
\begin{align}\label{bdd490}
\begin{split}
\Delta^{'}\left [(S+\kappa)V\right ]\geq & \frac{\kappa}{3}\widetilde W
+\frac{3}{4}V^2-C_{39}V^{\frac{3}{2}}-C_{40}V-C_{41}V^{\frac{1}{2}}\\
\geq & C_{42}\left [(S+\kappa)V\right ]^2-C_{43}
\left [(S+\kappa)V\right ]^{\frac{3}{2}}-C_{44}
\left [(S+\kappa)V\right ]-C_{45}\left [(S+\kappa)V\right ]^{\frac{1}{2}}
\end{split}
\end{align}
since $S+\kappa\geq\kappa\geq 1$ and $S+\kappa$ is bounded from above
uniformly.

By the work of Cheng and Yau in \cite{cy2}, we know inequality
\eqref{bdd490} implies that
$(S+\kappa)V$ is bounded. This implies $V$ is bounded since $V\leq
(S+\kappa)V$. Thus we obtain the $C^4$ estimate. By formula \eqref{bdd60}
we know that the curvature of the K\"ahler-Einstein metric $g$
is bounded. This finishes the proof of Theorem \ref{kebdd}.

Now we briefly describe how to control the covariant derivatives of
the curvature of the K\"ahler-Einstein metric.

Firstly, by
differentiating equation \eqref{bdd60}, we see that the boundedness
of the derivatives of the $P_{i\bar jk\bar l}$ is equivalent to the
boundedness of the covariant derivatives of $u$ with respect to the
background metric. Furthermore,  the derivatives involved are at least
order $2$ and were taken in at least one holomorphic direction and
one anti-holomorphic direction.

To bound such $k$-th order derivatives of $u$, we form the
quantity $S_k$ such that $S_3=S$,  $S_4=V$ and $S_5=W$. In
general, if we fix normal coordinates with respect to the
K\"ahler-Einstein metric at one point, then $S_k$ is a sum of
square of terms where each term is a covariant derivative of $u$
and the derivative is described above. All terms are obtained in
the following way:

For each covariant derivative of $u$ whose square appeared in the
sum $S_{k-1}$, we take covariant derivative of this term with
respect to the background metric in $z_\alpha$ and $\bar z_\beta$
respectively. Then we obtain two terms whose square appear in the
sum $S_k$. It is easy to see that $S_k$ is a sum of $2^{n-3}$
squares of certain covariant derivatives of $u$ where the derivatives
are of the type described above.

It is clear that the covariant derivatives of the curvature of the
K\"ahler-Einstein metric is bounded is equivalent to the fact that
the quantities $S_k$ is bounded.

We now estimate $S_k$ inductively. Assume $S_l$ is bounded for any
$l\leq k-1$, we compute
\[
\Delta^{'}\left ((S_{k-1}+\kappa)S_k
\right )
\]
where $\kappa$ is a large constant. Similar to
inequality \eqref{bdd180} we have
\[
\Delta^{'}\left ((S_{k-1}+\kappa)S_k\right )\geq (S_{k-1}+\kappa)
\Delta^{'}S_k+S_k\Delta^{'}S_{k-1}-2\left |\nabla^{'}S_{k-1}\right |
\left |\nabla^{'}S_{k}\right |.
\]
In the above formula, the leading term of $\Delta^{'}S_{k-1}$ is
$S_k$ as we did in formula \eqref{bdd390} and the term
$\nabla^{'}S_{k-1}$ is of order $S_k^{\frac{1}{2}}$ as we know in
formula \eqref{bdd430}. Similarly, the term  $\nabla^{'}S_{k}$ is
of order $S_{k+1}^{\frac{1}{2}}S_k^{\frac{1}{2}}$ as we did in
formula \eqref{bdd410}. When we compute  $\Delta^{'}S_{k}$, the
leading term is $S_{k+1}$. However, there will be products of
$(k+2)$-th order derivatives of $u$ and $k$-th order derivatives
of $u$. We can reduce the $(k+2)$-th order derivatives of $u$ by
using the Monge-Amp\'ere equation as we did in formulas
\eqref{bdd320}-\eqref{bdd3788}. That is, by differentiating
equation \eqref{bdd70} successively and by switching the order of
derivatives, we see that these products are of order at most
$S_{k+1}^{\frac{1}{2}}S_{k}^{\frac{3}{2}}$.

By using similar argument as above, finally we can derive an
inequality of form \eqref{bdd490} when $\kappa$ is large enough.
By using Cheng-Yau's work, we conclude that $S_k$ is bounded. The
computation is very long but straightforward. We omit it here for
simplicity .

\qed

As a direct corollary, we have
\begin{cor}
The injectivity radius of the K\"ahler-Einstein metric on the
Teichm\"uller space is bounded from below. Thus the Teichm\"uller
space equipped with this metric has bounded geometry.
\end{cor}

This corollary can be proved in the same way as Corollary
\ref{injectivity} by using the above theorem.


\begin{thebibliography}{10}

\bibitem{cy2}
S.~Y. Cheng and S.~T. Yau.
\newblock Differential equations on {R}iemannian manifolds and their geometric
  applications.
\newblock {\em Comm. Pure Appl. Math.}, 28(3):333--354, 1975.

\bibitem{fa}
G.~Faltings.
\newblock Arithmetic varieties and rigidity.
\newblock In {\em Seminar on number theory, Paris 1982--83 (Paris, 1982/1983)},
  volume~51 of {\em Progr. Math.}, pages 63--77. Birkh\"auser Boston, Boston,
  MA, 1984.

\bibitem{farkra1}
H.~M. Farkas and I.~Kra.
\newblock {\em Riemann surfaces}, volume~71 of {\em Graduate Texts in
  Mathematics}.
\newblock Springer-Verlag, New York, second edition, 1992.

\bibitem{ker1}
S.~P. Kerckhoff.
\newblock The {N}ielsen realization problem.
\newblock {\em Ann. of Math. (2)}, 117(2):235--265, 1983.

\bibitem{ko1}
S.~Kobayashi.
\newblock {\em Differential geometry of complex vector bundles}, volume~15 of
  {\em Publications of the Mathematical Society of Japan}.
\newblock Princeton University Press, Princeton, NJ, 1987.

\bibitem{ko2}
S.~Kobayashi.
\newblock {\em Hyperbolic complex spaces}, volume 318 of {\em Grundlehren der
  Mathematischen Wissenschaften [Fundamental Principles of Mathematical
  Sciences]}.
\newblock Springer-Verlag, Berlin, 1998.

\bibitem{lip1}
P.~Li.
\newblock Lecture notes.
\newblock 2003.

\bibitem{lischoen1}
P.~Li and R.~Schoen.
\newblock {$L\sp p$} and mean value properties of subharmonic functions on
  {R}iemannian manifolds.
\newblock {\em Acta Math.}, 153(3-4):279--301, 1984.

\bibitem{lsy1}
K.~Liu, X.~Sun, and S.-T. Yau.
\newblock Canonical metrics on the moduli space of riemann surface {I}.
\newblock Preprint, 2004.

\bibitem{ls2}
Z.~Lu and X.~Sun.
\newblock On the {W}eil-{P}etersson volume and the first chern class o of the
  {M}oduli {S}paces of {C}alabi -{Y}au {M}anifolds.
\newblock submitted to Communications in Mathematical Physics, 2003.

\bibitem{luo1}
F.~Luo.
\newblock Oral communication.
\newblock 2004.

\bibitem{ma1}
H.~Masur.
\newblock Extension of the {W}eil-{P}etersson metric to the boundary of
  {T}eichmuller space.
\newblock {\em Duke Math. J.}, 43(3):623--635, 1976.

\bibitem{mcc}
J.~McCarthy.
\newblock A ``{T}its-alternative'' for subgroups of surface mapping class
  groups.
\newblock {\em Trans. Amer. Math. Soc.}, 291(2):583--612, 1985.

\bibitem{yau4}
R.~Schoen and S.-T. Yau.
\newblock {\em Lectures on differential geometry}.
\newblock Conference Proceedings and Lecture Notes in Geometry and Topology, I.
  International Press, Cambridge, MA, 1994.
\newblock Lecture notes prepared by Wei Yue Ding, Kung Ching Chang [Gong Qing
  Zhang], Jia Qing Zhong and Yi Chao Xu, Translated from the Chinese by Ding
  and S. Y. Cheng, Preface translated from the Chinese by Kaising Tso.

\bibitem{shi1}
W.-X. Shi.
\newblock Ricci flow and the uniformization on complete noncompact {K}\"ahler
  manifolds.
\newblock {\em J. Differential Geom.}, 45(1):94--220, 1997.

\bibitem{th}
W.~P. Thurston.
\newblock On the geometry and dynamics of diffeomorphisms of surfaces.
\newblock {\em Bull. Amer. Math. Soc. (N.S.)}, 19(2):417--431, 1988.

\bibitem{tr1}
S.~Trapani.
\newblock On the determinant of the bundle of meromorphic quadratic
  differentials on the {D}eligne-{M}umford compactification of the moduli space
  of {R}iemann surfaces.
\newblock {\em Math. Ann.}, 293(4):681--705, 1992.

\bibitem{wol1}
S.~A. Wolpert.
\newblock The hyperbolic metric and the geometry of the universal curve.
\newblock {\em J. Differential Geom.}, 31(2):417--472, 1990.

\bibitem{yau3}
S.~T. Yau.
\newblock On the {R}icci curvature of a compact {K}\"ahler manifold and the
  complex {M}onge-{A}mp\`ere equation. {I}.
\newblock {\em Comm. Pure Appl. Math.}, 31(3):339--411, 1978.

\bibitem{yau2}
S.-T. Yau.
\newblock {\em Nonlinear analysis in geometry}, volume~33 of {\em Monographies
  de L'Enseignement Math\'ematique [Monographs of L'Enseignement
  Math\'ematique]}.
\newblock L'Enseignement Math\'ematique, Geneva, 1986.
\newblock S\'erie des Conf\'erences de l'Union Math\'ematique Internationale
  [Lecture Series of the International Mathematics Union], 8.

\end{thebibliography}
\end{document}